\documentclass[10pt]{article}
\usepackage{amsmath}
\usepackage{amssymb}
\usepackage{amscd}

\def\oplusinf{\mathop{\oplus}}
\def\otimesinf{\mathop{\otimes}}

\def\dim{{\mathrm{dim}}}
\def\ker{{\mathrm{Ker}}}

\def\aut{{\mathrm{Aut}}}

\def\Talpha#1{\vbox{\ialign{##\crcr
      $\alpha$\crcr\noalign{\kern2pt\nointerlineskip}
	   $\hfil\displaystyle{#1}\hfil$\crcr}}}
\def\vol{{\mathrm{vol}}}

\def\alg{{\mathrm{alg}}}

\def\lie{{\mathbf{Lie}}}

\def\calo{{\mathcal O}}

\def\cala{{\mathcal A}}

\def\call{{\mathcal L}}
\def\calm{{\mathcal M}}

\def\calh{{\mathcal H}}
\def\cals{{\mathcal S}}

\def\cals{{\mathcal S}}

\def\bbbone{\mbox{\rm 1\hspace {-.6em} l}}
\def\us{\mathbf{u}}

\def\gr{\mathrm{\scriptsize gr}}

\def\diagth{\mathrm{\scriptsize Diag}}

\def\detf{\mathrm{det}}

\def\pol{{C_{\mathrm{alg}}}}
\def\calg{\pol}
\def\Halg {{H_{\mathrm{alg}}}}
\def\Oalg {{\Omega_{\mathrm{alg}}}}

\def\clif{{\mathrm{Cliff}}}
\def\ch{{\mathrm{ch}}}
\def\tr{{\mathrm{tr}}}

\def\rnt{{\mathbb R^{2n}_{\theta}}}
\def\rntu{{\mathbb R^{2n+1}_{\theta}}}
\def\snt{{S^{2n}_{\theta}}}

\def\Cb{{\mathbb  C}}
\def\Zb{{\mathbb  Z}}
\def\build#1_#2^#3{\mathrel{
\mathop{\kern 0pt#1}\limits_{#2}^{#3}}}

\newtheorem{lemma}{LEMMA}
\newtheorem{proposition}{PROPOSITION}
\newtheorem{corol}{COROLLARY}
\newtheorem{theo}{THEOREM}

\begin{document}
\thispagestyle{empty}
\enlargethispage{3cm}
\begin{center}
{\bf NONCOMMUTATIVE FINITE-DIMENSIONAL MANIFOLDS}
\end{center}
\begin{center}
{\bf I. SPHERICAL MANIFOLDS AND RELATED EXAMPLES}
\end{center}
\vspace{0,5cm}

\begin{center} Alain CONNES\footnote{Coll\`ege de France, 3 rue
d'Ulm, 75  005 Paris, and \\ I.H.E.S., 35 route de Chartres, 91440
Bures-sur-Yvette\\
connes$@$ihes.fr}
and Michel DUBOIS-VIOLETTE
\footnote{Laboratoire de Physique Th\'eorique, UMR 8627\\
Universit\'e Paris XI,
B\^atiment 210\\ F-91 405 Orsay Cedex, France\\
Michel.Dubois-Violette$@$th.u-psud.fr}\\
\end{center}
\vspace{0,5cm}
\begin{abstract}

       We exhibit large classes of examples of noncommutative 
finite-dimen\-sional
      manifolds which are (non-formal) deformations of classical
      manifolds. The main result of this paper is a
      complete description of noncommutative three-dimensional
spherical manifolds, a noncommutative version of the sphere $S^3$
defined by basic K-theoretic equations. We find a
3-parameter family of deformations $S^3_{\mathbf{u}}$ of the standard 
3-sphere $S^3$
and a corresponding 3-parameter deformation of the 4-dimensional
Euclidean space $\mathbb R^4$. For generic values of the deformation
parameters we show that the obtained algebras of polynomials on the
deformed $\mathbb R^4_{\mathbf{u}}$ only depend on two parameters and 
are isomorphic to the algebras introduced
by Sklyanin in connection with the Yang-Baxter equation.
It follows that different $S^3_{\mathbf{u}}$ can span the same 
$\mathbb R^4_{\mathbf{u}}$. This equivalence generates a foliation of 
the parameter space $\Sigma $. This foliation admits singular leaves 
reduced to a point. These critical points are either isolated or fall 
in two 1-parameter families $ C_{\pm}\subset \Sigma $. Up to the 
simple operation of taking the fixed algebra by an involution, these 
two families are identical and we concentrate here on $ C_{+}$. For 
$u \in C_{+}$ the above isomorphism with the Sklyanin algebra breaks 
down and the corresponding algebras are special cases of 
$\theta$-deformations, a notion which we generalize in
      any dimension and various contexts, and  study in some
      details.
Here, and this point is crucial,
      the dimension is not an artifact, i.e. the dimension of
      the classical model, but is the Hochschild dimension of the
      corresponding algebra which remains constant during the
      deformation.
      Besides the standard noncommutative tori,
      examples of $\theta$-deformations include the recently defined
      noncommutative 4-sphere $S^4_{\theta}$ as well as $
      m$-dimensional generalizations, noncommutative versions of
      spaces $\mathbb R^m$ and quantum groups which are deformations of
      various classical groups. We develop general tools such as the 
twisting of the Clifford algebras in order to exhibit the spherical 
property of the hermitian projections
      corresponding to the noncommutative $2n$-dimensional spherical
      manifolds $S^{2n}_{\theta}$. A key result is the differential self-duality
      properties of these projections which generalize the 
self-duality of the round
      instanton.

      \end{abstract}

\section{Introduction}

Our aim in this paper is to describe large classes of tractable
concrete examples of {\sl noncommutative manifolds}. Our original
motivation is the problem of classification of {\sl spherical}
noncommutative manifolds which arose from the basic discussion of
Poincar\'e duality in $K$-homology \cite{connes:62}, \cite{connes:08}.\\

\noindent The algebra $\cala$ of functions on a spherical
noncommutative manifold $S$ of dimension $n$ is generated by the
matrix components of a cycle $x$ of the $K$ theory of $\cala$, whose
dimension is the same as $n=\dim\ (S)$.\\
More specifically, for $n$ even, $n=2m$, the algebra $\cala$ is
generated by the matrix elements $e^i_{j}$ of a self-adjoint idempotent
\begin{equation}
e=[e^i_{j}]\in M_{q}(\cala),\> \> e= e^2=e^\ast,
\label{al1}
\end{equation}
and one assumes that all the components $\ch_{k}(e)$ of the Chern
character of $e$ in cyclic homology satisfy,
\begin{equation}
\ch_{k}(e)=0\> \> \forall k=0,1,\dots,m-1
\label{al2}
\end{equation}
while $\ch_{m}(e)$ defines a non zero Hochschild cycle playing the
role of the volume form of $S$.\\
For $n$ odd the algebra $\cala$ is similarly generated by the matrix
components $U^i_{j}$ of a unitary
\begin{equation}
U=[U^i_{j}] \in M_{q}(\cala),\> \> UU^\ast=U^\ast U=1
\label{al3}
\end{equation}
and, with $n=2m+1$, the vanishing condition (\ref{al2}) becomes
\begin{equation}
\ch_{k+\frac{1}{2}}(U)=0\> \> \> \forall k=0,1,\dots,m-1.
\label{al4}
\end{equation}
The components $\ch_{k}$ of the Chern character in cyclic homology
are the following explicit elements of the tensor product
\begin{equation}
\cala\otimes (\tilde\cala)^{\otimes 2k}
\label{al5}
\end{equation}
where $\tilde \cala$ is the quotient of $\cala$ by the subspace
$\mathbb C 1$,
\begin{equation}
\ch_{k}(e)= \left(e^{i_{0}}_{ i_{1}}-\frac{1}{2}\delta^{i_{0}}_{i_{1}}
\right)\otimes e^{i_{1}}_{i_{2}}\otimes e^{i_{2}}_{i_{3}}\otimes \dots \otimes
e^{i_{2k}}_{i_{0}}
\label{al6}
\end{equation}
and
\begin{equation}
\ch_{k+\frac{1}{2}}(U)=U^{i_{0}}_{i_{1}} \otimes U^{\ast
i_{1}}_{i_{2}} \otimes  U^{i_{2}}_{i_{3}}\otimes \dots \otimes
U^{\ast i_{2k+1}}_{i_{0}}-U^{\ast i_{0}}_{i_{1}}\otimes \dots \otimes
U^{i_{2k+1}}_{i_{0}}
\label{al7}
\end{equation}
up to an irrelevant normalization constant.\\

\noindent It was shown in \cite{connes:62} that the Bott generator on the
classical spheres
$S^n$ give solutions to the above equations (\ref{al2}), (\ref{al4})
and in \cite{connes:08} that non trivial noncommutative solutions exist for
$n=3$, $q=2$ and $n=4$, $q=4$.\\

\noindent In fact, as will be explained in our next paper (Part II), 
consistency
with the suspension functor requires a coupling between the dimension
$n$ of $S$ and $q$. Namely $q$ must be the same for $n=2m$ and $n=2m+1$
whereas it must be doubled when going from $n=2m-1$ to $n=2m$. This
implies that for dimensions $n=2m$ and $n=2m+1$, one has $q=2^m
q_{0}$ for some $q_{0}\in \mathbb N$. Furthermore the normalization
$q_{0}=1$ is induced by the identification of $S^2$ with
one-dimensional projective space $P_{1}(\mathbb C)$ (which means
$q=2$ for $n=2$). We shall take this convention (i.e. $q=2^m$ for
$n=2m$ and $n=2m+1$) in the following.\\

\noindent The main result of the present paper is the complete description of
the noncommutative solutions for $n=3$ ($q=2$). We find a
three-parameter family of deformations of the standard three-sphere $S^3$
and a corresponding 3-parameter deformation of the 4-dimensional
Euclidean space $\mathbb R^4$.  For generic values of the deformation
parameters we show that the obtained algebras of polynomials on the
deformed $\mathbb R^4_{\mathbf{u}}$ only depend on two parameters and 
are isomorphic to the algebras introduced
by Sklyanin in connection with the Yang-Baxter equation.
It follows that different $S^3_{\mathbf{u}}$ can span the same 
$\mathbb R^4_{\mathbf{u}}$. This equivalence relation generates a 
foliation of the parameter space $\Sigma $. This foliation admits 
singular leaves reduced to a point. These critical points are either 
isolated or fall in two 1-parameter families $ C_{\pm}\subset \Sigma 
$. Up to the simple operation of taking the fixed algebra by an 
involution, these two families are identical and we concentrate here 
on $ C_{+}$. For $u \in C_{+}$ the above isomorphism with the 
Sklyanin algebra breaks down and the corresponding algebras are 
special cases of $\theta$-deformations. It gives rise to a
one-parameter deformation $\mathbb C^2_{\theta}$ of $\mathbb C^2$
(identified with $\mathbb R^4$) which
is well suited for simple higher dimensional generalizations (i.e.
$\mathbb C^2$ replaced by $\mathbb C^n\simeq \mathbb R^{2n}$). We
shall describe and analyse them in details to understand  this 
particular critical case, while the general case (of generic values 
of the parameters) will be treated in Part II.\\

  First we shall show that, unlike
most deformations used to produce noncommutative spaces from classical
ones, the above deformations do not alter the Hochschild dimension. The
latter is the natural generalization of the notion of dimension to the
noncommutative case and is the smallest integer $m$ such that the
Hochschild homology of $\cala$ with values in a bimodule $\calm$
vanishes for $k>m$ ($H^k(\cala,\calm)=0$ $\forall k>m$). Second we
shall describe the natural notion of differential forms on the above
noncommutative spaces and obtain the natural quantum groups of
symmetries as ``$\theta$-deformations" of the classical groups
$GL(m,\mathbb R)$, $SL(m,\mathbb R)$ and $GL(n,\mathbb C)$.\\

\noindent The algebraic versions of differential forms on the above quantum
groups turn out to be graded involutive differential Hopf algebras,
which implies that the corresponding differential
calculi are bicovariant in the sense of \cite{slw}. It is worth
noticing here that conversely as shown in \cite{brz}, a bicovariant
differential calculus on a quantum group always comes from a graded
differential Hopf algebra as above.\\

\noindent Finally we shall come back to the metric aspect of the construction
which was the original motivation for the definition of spherical
manifolds from the polynomial operator equation
fulfilled by the Dirac operator.\\

\noindent We shall check in detail that $\theta$-deformations of Riemannian
spin geometries fulfill all axioms of noncommutative geometry,
thus completing the path, in the special case of $\theta$-deformations,
from the crudest level of the algebra $\calg(S)$ of polynomial functions on $S$
  to the full-fledged structure
of noncommutative geometry \cite{connes:61}.  \\

\noindent In the course of the paper it will be shown that the self-duality
property of the round instanton on $S^4$ extends directly to the
self-adjoint idempotent identifying $S^4_{\theta}$ as a
noncommutative 4-dimensional
spherical manifold and that, more generally, the self-adjoint idempotents
corresponding to the noncommutative $2n$-dimensional spherical
manifolds $S^{2n}_{\theta}$ defined below satisfy a differential
self-duality property which is a direct extension of the one
satisfied by their classical counterparts as explained in
\cite{mdv2}.\\

\noindent In conclusion the above examples appear as an interesting point of
contact between various approaches to noncommutative geometry. The
original motivation came from the operator equation of degree $n$
fulfilled by the Dirac operator of an $n$-dimensional spin
manifold \cite{connes:61}. The simplest equation ``quantizing" the 
corresponding
Hochschild cycle $c$, namely $c=\ch (e)$ (\cite{connes:62}) led to 
the definition of
spherical manifolds. What we show here is that in the simplest non
trivial case ($n=3$, $q=2$) the answer is intimately related to the
Sklyanin algebras which play a basic role in noncommutative algebraic
geometry.\\

\noindent Many algebras occuring in this
paper are  finitely
generated and finitely presented.
  These algebras are viewed as algebras of polynomials on the
corresponding noncommutative space $S$ and we denote them by $\calg(S)$.
  With these notations $\calg(S)$ has to be
distinguished from $C^\infty(S)$, the algebra of smooth
functions on $S$ obtained as a suitable completion of $\calg(S)$.
Basic algebraic properties such as Hochschild dimension are
not necessarily preserved under the transition from $\calg(S)$
to $C^\infty(S)$. The topology of $S$ is specified by the $C^*$
completion of $C^\infty(S)$.

\noindent The plan of the paper is the following. After this introduction, in
section~\ref{sec02}, we give a complete description of
noncommutative spherical manifolds for the lowest non trivial
dimension : Namely for dimension $n=3$ and for $q=2$. These
form a 3-parameter family $S^3_{\us}$ of deformations of the standard
3-sphere as explained above and correspondingly one has a homogeneous
version which is a three-parameter family $\mathbb R^4_{\us}$ of
deformations of the standard 4-dimensional Euclidean space $\mathbb
R^4$.  We then consider
their suspensions and show that the suspension $S^4_{\us}$ of
$S^3_{\us}$ is a four-dimensional noncommutative spherical manifold
(with $q=4=2^2$). In Section~\ref{sec03}, we show that for generic
values of the parameters, the algebra $\calg(\mathbb R^4_{\us})$ of polynomial
functions on the noncommutative  $\mathbb R^4_{\us}$ reduces to a Sklyanin
algebra \cite{Skly:1}, \cite{Skly:2}.These
Sklyanin algebras have been intensively studied \cite{Od:fe},
\cite{sm:st}, from the point of view of noncommutative algebraic geometry
but we postpone their analysis to Part II of this paper.
We concentrate instead on the determination of the scaling foliation 
of the parameter space $\Sigma$ for 3-dimensional
spherical manifolds $S^3_{\mathbf{u}}$. Different $S^3_{\mathbf{u}}$
can span isomorphic 4-dimensional $\mathbb R^4_{\mathbf{u}}$ and we
shall determine here the corresponding foliation of $\Sigma$
 using the geometric data
associated \cite{Od:fe}
\cite{art:1}\cite{art:2}
to such algebras. This will allow us to isolate the critical points in the 
parameter space $\Sigma$ and we devote the end of the paper to the study
of the corresponding algebras. The simplest way to analyse them
is to view them as a special case of a general procedure of $\theta$-deformation.
In  Section~\ref{sec2} we define a noncommutative deformation $\rnt$
of $\mathbb
R^{2n}$ for $n\geq 2$ which is coherent with the identification
$\mathbb C^n=\mathbb R^{2n}$ as real spaces and is also consequently
a noncommutative deformation $\mathbb C^n_{\theta}$ of $\mathbb C^n$.
For $n=2$, $\mathbb R^4_{\theta}$ reduces to the above one-parameter family of
deformations of $\mathbb R^4$ which is included in the multiparameter
deformation $\mathbb R^4_{\us}$ of Section \ref{sec02}.
We introduce in this section a deformation of the generators of the
Clifford algebra of $\mathbb R^{2n}$ which will be very useful for
the computations of Sections~\ref{sec3} and ~\ref{sec10}. In
Section~\ref{sec3} we define noncommutative versions
$\rntu$, $S^{2n}_{\theta}$ and $S^{2n-1}_{\theta}$ of $\mathbb
R^{2n+1}$, $S^{2n}$ and $S^{2n-1}$ for $n\geq 2$.  For  $n=2$,
$S^{2n}_{\theta}$ reduces to the noncommutative $4$-sphere
$S^4_{\theta}$ of \cite{connes:08} whereas $S^{2n-1}_{\theta}$
reduces to the one-parameter family $S^3_{\theta}$ of deformation of
$S^3$ associated to the non-generic values of $\us$.  We generalize
the results of
\cite{connes:08} on $S^4_{\theta}$ to $S^{2n}_{\theta}$ for arbitrary $n\geq 2$
and we describe their counterpart for the odd-dimensional cases
$S^{2n-1}_{\theta}$ showing thereby that these $S^m_{\theta}$ ($m\geq
3$) are noncommutative spherical manifolds.  Furthermore, it will be
shown later (in Section \ref{sec10}) that the defining hermitian
projections of
$S^{2n}_{\theta}$ possess differential self-duality properties which
generalize the ones of their classical counterpart (i.e. for $S^{2n}$)
as explained in \cite{mdv2}.
In Section~\ref{sec4}, we define
algebraic versions of differential forms on the above noncommutative
spaces. These definitions, which are essentially unique, provide
dense subalgebras of the canonical algebras of smooth differential
forms defined in
Sections~\ref{sec9}, \ref{sec10}, \ref{sec11} for these particular cases.
These differential calculi are diagonal  \cite{mdv:pm1} which implies
that they are quotients of the corresponding universal diagonal
differential calculi \cite{mdv3}. In Section~\ref{sec5} we construct
quantum groups which are deformations  (called $\theta$-deformations)
of the classical
groups $GL(m,\mathbb R)$, $SL(m,\mathbb R)$ and $GL(n,\mathbb C)$
for $m\geq 4$ and $n\geq 2$. The point of view for this construction
is close to the one of \cite{malt} which is itself a generalization of
a construction described in  \cite{man:0}, \cite{manin}.  It is
pointed out that
there is no such $\theta$-deformation of $SL(n,\mathbb C)$ although there
is a $\theta$-deformation of the subgroup of $GL(n,\mathbb C)$
consisting of matrices with determinants of modulus one ($\vert
\det_{\mathbb C}(M)\vert^2=1$).
In Section~\ref{sec6} we define the
corresponding deformations of the groups $O(m)$, $SO(m)$ and $U(n)$.
As above there is no $\theta$-deformation of
$SU(n)$ which is the counterpart of the same statement for
$SL(n,\mathbb C)$. All the quantum groups $G_{\theta}$ considered in
Section~\ref{sec5}  and in Section~\ref{sec6} are matrix quantum
groups \cite{slw:0}
and in fact as coalgebra $\calg(G_{\theta})$ is undeformed i.e.
isomorphic to the classical coalgebra $\calg(G)$ of representative
functions on $G$ \cite{mdv2b}, (only the associative product is
deformed). In Section~\ref{sec7}, we analyse the structure of the
algebraic version of differential forms on the above quantum groups.
These graded-involutive differential algebras turn out to be
{\sl graded-involutive differential Hopf algebras} (with
coproducts and counits extending the original ones) which, in view of
\cite{brz}, means that the corresponding differential calculi are
bicovariant in the sense of \cite{slw}. It is worth noticing that the
above $\theta$-deformations of $\mathbb R^m$, of the differential
calculus on $\mathbb R^m$ and of some classical groups have been
already considered for instance in \cite{asch}. The quantum group
setting analysis of \cite{asch} is clearly very
interesting: There,
$\mathbb R^m_{\theta}$ appears (with other notations) as a quantum
space on which some quantum group acts (or more precisely as a quantum
homogeneous space) and the differential calculus on $\mathbb
R^m_{\theta}$ is the covariant one. Another powerful approach
to the above quantum group aspects is to make use of the notion
  of Drinfeld twist \cite{Drin} since it is clear that the
$\theta$-deformed quantum group of Sections 7 and 8 can be
obtained by twisting (see e.g. in \cite{sit} for a particular case);
  thus many results of Sections 7 to 9 can be obtained by using for
instance Proposition 2.3.8 of \cite{maj2},
its graded counterpart and the result of \cite{M.Oe} for
the differential calculus in this case. Here the emphasis is rather
different. The noncommutative $\mathbb R^m_{\theta}$ appears as a
solution of the $K$-theoretic equations (1.2) or (1.4) appropriate to
the dimension $m$ and the differential calculus which is essentially
unique is used  to produce the projective resolution of
$C^\infty(\mathbb R^m_{\theta})$ which ensures that the Hochschild
dimension of $C^\infty(\mathbb R^m_{\theta})$ is $m$ (i.e. that
$R^m_{\theta}$ is $m$-dimensional). It turns out that the
differential calculus on $\mathbb R^m_{\theta}$ is covariant for some
quantum group actions and that these quantum groups are again
$\theta$-deformations. However, our interest in $\theta$-deformation
is connected to the fact that it preserves the Hochschild dimension.
Furthermore the analysis of Section~12 shows that in general the
differential calculi over $\theta$-deformations do not rely on the
existence of quantum group actions, (see below).
In Section~\ref{sec8}, we define the splitting
homomorphisms mapping the polynomial algebras $\calg$ of the various
$\theta$-deformations introduced previously into the polynomial
algebras on the product of the corresponding classical spaces with the
noncommutative $n$-torus $T^n_{\theta}$. In Section~\ref{sec9} we use the
splitting homomorphisms to produce smooth structures on the previously
defined noncommutative spaces, that is the algebras of smooth
functions and of smooth differential forms. In Section~\ref{sec10}  we describe
in general the construction which associates to each
finite-dimensional manifold $M$ endowed with a smooth
action $\sigma$ of the $n$-torus $T^n$ a noncommutative deformation
$C^\infty(M_{\theta})$ of the
algebra $C^\infty(M)$  of smooth functions on $M$ (and of the algebra of smooth
differential forms) which defines the noncommutative manifold $M_{\theta}$
and we explain why the Hochschild dimension of the deformed algebra
remains constant
and equal to the dimension of $M$. The construction of differential
forms given in this section shows that the $\theta$-deformation of
differential forms does not rely on a quantum group action since
generically there is no such an action on $M_{\theta}$ (beside the
action of the $n$-torus). The deformation
$C^\infty(M_{\theta})$ of the algebra
$C^\infty(M)$ is a special case of Rieffel's deformation quantization
\cite{rief:3} and close to the form
adopted in \cite{rief:4}. It is worth noticing here that at the
formal level deformations of algebras for actions of $\mathbb R^n$
have been also analysed in \cite{mourre}. It is however crucial
  to consider (non formal) actions of $T^n$; our results
would be generically wrong for actions of $\mathbb R^n$.

\noindent In Section~\ref{sec11} we analyse the metric aspect of the 
construction showing
that the deformation is isospectral in the sense of \cite{connes:08}
and that our construction gives an alternate setting for results like
Theorem 6 of \cite{connes:08}.  We use the
splitting homomorphism to show
that when $M$ is a compact riemannian spin manifold endowed with an
isometric action of $T^n$ the corresponding spectral triple (\cite{connes:08})
$(C^\infty(M_{\theta}), \calh_{\theta},D_{\theta})$ satisfies the
axioms of noncommutative geometry of \cite{connes:61}.
We show moreover (theorem 9) that any $T^n$-invariant metric
on $S^m$, ($m=2n$, $2n-1$), whose volume form is rotation invariant
yields a solution of the original polynomial equation for the
Dirac operator on $S_{\theta}^m$.  Section
\ref{seconclu} is our conclusion.\\

\noindent Throughout this paper $n$ denotes an integer such that $n\geq 2$,
$\theta\in M_n(\mathbb R)$ is an antisymmetric real ($n,n$)-matrix
with matrix elements denoted by $\theta_{\mu\nu}$ ($\mu,
\nu=1,2,\dots,n$) and we set $\lambda^{\mu\nu}=e^{
i\theta_{\mu\nu}}=\lambda_{\mu\nu}$.  The reason for this double
notation $\lambda^{\mu\nu},\lambda_{\mu\nu}$ for the same object
$e^{ i\theta_{\mu\nu}}$ is to avoid ambiguities connected with the
Einstein  summation convention (of repeated up down indices) which is
used throughout.  The
symbol $\otimes$ without other specification will always denote the
tensor product over the field $\mathbb C$.  A {\sl  self-adjoint
idempotent} or a {\sl hermitian projection} in a $\ast$-algebra is an
element $e$ satisfying $e^2=e=e^\ast$.
By a {\sl graded-involutive algebra} we here mean a graded $\mathbb
C$-algebra endowed with an antilinear involution $\omega\mapsto\bar
\omega$ such that $\overline{\omega\omega'}=(-1)^{pp'}\bar \omega'\bar
\omega$ for $\omega$ of degree $p$ and $\omega'$ of degree $p'$.  A
{\sl graded-involutive differential algebra} will be a
graded-involutive algebra endowed with a real differential
  $d$ such that $d(\bar\omega)=\overline{d(\omega)}$ for
any $\omega$. Given a graded vector space $V=\oplus_{n}V^n$, we denote
by $(-I)^\gr$ the endomorphism of $V$ which is the identity mapping
on $\oplus_{k} V^{2k}$ and minus the identity mapping on $\oplus_{k}
V^{2k+1}$. If $\Omega'$ and $\Omega''$ are graded algebras one can
endow $\Omega'\otimes \Omega''$ with  the usual product $(x'\otimes
x'')(y'\otimes y'')=x'y'\otimes x''y''$  or with the graded
twisted one $(x'\otimes x'')(y'\otimes y'')=(-1)^{\vert x''\vert\vert
y'\vert} x'y'\otimes x''y''$ where $\vert x''\vert$ is the degree of
$x''$ and $\vert y'\vert$ is the degree of $y'$; in the latter case we
denote by $\Omega'\otimes_{\gr}\Omega''$ the corresponding graded
algebra. If furthermore $\Omega'$ and $\Omega''$ are graded
differential algebras $\Omega'\otimes_{\gr}\Omega''$ will denote the
corresponding graded algebra endowed with the differential
$d=d'\otimes I+(-I)^\gr \otimes d''$. A bimodule over an algebra $A$
is said to be {\sl diagonal} if it is a subbimodule of $A^I$ for some
set $I$.  Concerning locally convex
algebras, topological modules, bimodules and resolutions we use the
conventions of \cite{connes:02}. All our locally convex algebras and
locally convex modules will be nuclear and complete. Finally we shall
need some notations  concerning matrix algebras $M_{n}(A)=M_{n}(\mathbb
C)\otimes A$ with entries in an algebra $A$.  For $M\in M_{n} (A)$, we
denote by $\tr(M)$ the element $\sum^n_{\alpha=1}M^\alpha_{\alpha}$ of $A$.
If $M$ and $N$ are in $M_{n}(A)$, we denote by $M\circledcirc N$ the
element of $M_{n}(A\otimes A)$ defined by $(M\circledcirc
N)^\alpha_{\beta}=M^\alpha_{\gamma}\otimes M^\gamma_{\beta}$.

\section{Noncommutative 3-spheres and 4-planes}\label{sec02}
\setcounter{equation}{0}
Our aim in this section is to give a complete description
  of
noncommutative spherical three-manifolds. More specifically we give
here a complete description of the class
of complex unital $\ast$-algebras $\cala^{(1)}$ satisfying the
following conditions (I$_{1}$) and (II):\\

(I$_{1}$) $\cala^{(1)}$ is generated as unital $\ast$-algebra by the
matrix elements of a unitary $U\in M_{2}(\cala^{(1)})=M_{2}(\mathbb
C)\otimes \cala^{(1)}$, \\

(II) $U$ satisfies $\ch_{\frac{1}{2}}(U)=U^i_{j}\otimes U^{\ast
j}_{i}-U^{\ast i}_{j}\otimes U^j_{i}=0$\\

\noindent (i.e. with the notations explained above
$\tr(U\circledcirc  U^\ast-U^\ast
\circledcirc U)=0$).\\

\noindent It is convenient to consider the corresponding  homogeneous
problem, i.e. the
class of unital $\ast$-algebras $\cala$ such that \\

(I) $\cala$ is generated by the matrix elements of a $U\in
M_{2}(\cala)=M_{2}(\mathbb C)\otimes \cala$ satisfying $U^\ast
U=UU^\ast\in \bbbone_{2} \otimes \cala$ where $\bbbone_{2}$ is the
unit of $M_{2}(\mathbb C)$ and,\\

(II) $U$ satisfies
$\tilde\ch_{\frac{1}{2}}(U)=U^i_{j}\otimes U^{\ast
j}_{i}-U^{\ast i}_{j}\otimes U^j_{i}=0$\\

\noindent i.e. $\tr(U\circledcirc  U^\ast-U^\ast
\circledcirc U)=0$.

\noindent Notice that if $\cala^{(1)}$ satisfies Conditions (I$_{1}$) and
(II) or if $\cala$ satisfies Conditions (I) and (II) with $U$ as
above, nothing changes if one makes the replacement
\begin{equation}
U\mapsto U'=uV_{1}UV_{2}
\label{tr1}
\end{equation}
with $u=e^{i\varphi}\in U(1)$ and $V_{1},V_{2}\in SU(2)$. This
corresponds to a linear change in  generators, $(\cala^{(1)},U')$
satisfies (I$_{1}$) and (II) whenever $(\cala^{(1)},U)$ satisfies
(I$_{1}$) and (II) and $(\cala,U')$ satisfies (I) and (II)
whenever $(\cala,U)$ satisfies (I) and (II).\\

\noindent Let $\cala$ be a unital $\ast$-algebra and $U\in 
M_{2}(\cala)$. We use
the standard Pauli matrices $\sigma_{k}$ to write $U$ as
\begin{equation}
U=\bbbone_{2}z^0 + i\sigma_{k}z^k
\label{eqU}
\end{equation}
where $z^\mu$ are elements of $\cala$ for $\mu=0,1,2,3$. In terms of
the $z^\mu$, the transformations (\ref{tr1}) reads
\begin{equation}
z^\mu\mapsto u S^\mu_{\nu}z^\nu
\label{tr2}
\end{equation}
with $u\in U(1)$ as above and where $S^\mu_{\nu}$ are the matrix
elements of the real rotation $S\in SO(4)$ corresponding to
$(V_{1},V_{2})\in SU(2)\times SU(2)=\mbox{Spin}(4)$.
The pair
$(\cala,U)$ fulfills (I) if and only if $\cala$ is generated by
the $z^\mu$ as unital $\ast$-algebra and the $z^\mu$ satisfy
\begin{eqnarray}
   z^k z^{0\ast}-z^0z^{k\ast}+\epsilon_{k\ell m}z^\ell z^{m\ast}=0\label{aaa}\\
   z^{0\ast}z^k-z^{k\ast}z^0 + \epsilon_{k\ell m} z^{\ell \ast} z^m =0
   \label{aab}\\
   \sum^3_{\mu=0}(z^\mu z^{\mu\ast}-z^{\mu\ast}z^\mu)=0 \label{aac}
   \end{eqnarray}
   for $k=1, 2, 3$, where $\epsilon_{k\ell m}$ is completely
   antisymmetric in $k,\ell ,m\in \{1,2,3\}$ with $\epsilon_{1 2 3}=1$.
   Condition  (I$_{1}$) is satisfied if and only if one
   has in addition $\sum_{\mu}z^{\mu \ast} z^\mu=\bbbone$. The following
   lemma shows that there is no problem to pass from (I) to (I$_{1}$)
   just imposing the relation $\sum_{\mu}z^{\mu\ast}z^\mu
   -\bbbone=0$.

   \begin{lemma}\label{lemu}

    Let $\cala$, U satisfy $\mathrm{(I)}$ as above. Then $\sum^3_{\mu=0}$
    $z^{\mu\ast} z^\mu$ is in the center of $\cala$.
    \end{lemma}
    This result is easily verified using relations (\ref{aaa}),
    (\ref{aab}), (\ref{aac}) above. \\

\noindent Let us now investigate condition (II). In terms of the
    representation (\ref{eqU}), for $U$, condition (II) reads
    \begin{equation}
    \sum^3_{\mu=0}(z^{\mu\ast}\otimes z^\mu-z^\mu\otimes z^{\mu\ast})=0
    \label{aad}
    \end{equation}
    for the $z^\mu\in \cala$. One has the following result.

    \begin{lemma} \label{lemlambda}
     Condition $\mathrm{(II)}$, i.e. equation (\ref{aad}), is satisfied
if and only if     there is a symmetric unitary matrix $\Lambda \in 
M_{4}(\mathbb
     C)$ such that $z^{\mu\ast}=\Lambda^\mu_{\nu} z^\nu$ for
     $\mu\in\{0,\cdots,3\}$.
     \end{lemma}

   \noindent  The existence of $\Lambda\in M_{4}(\mathbb C)$ as above clearly
     implies Equation (\ref{aad}). Conversely assume that (\ref{aad}) is
     satisfied. If the $(z^\mu)$ are linearly independent, the existence
     and uniqueness of a matrix $\Lambda$ such that
     $z^{\mu\ast}=\Lambda^\mu_{\nu}z^\nu$ is immediate, and the symmetry
     and unitarity of $\Lambda$ follow from its uniqueness. Thus the
     only difficulty is to take care in general of the linear dependence
     between the $(z^\mu)$. We let $I\subset\{0,\dots,3\}$ be a maximal
     subset of $\{0,\dots,3\}$ such that the $(z^i)_{i\in I}$ are
     linearly independent. Let $I'$ be its complements ; one has
$z^{i'}=\bar L^{i'}_{i}z^i$ for some
     $L^{i'}_{i}\in \mathbb C$.  On the other hand one has
     $z^{i\ast}=C^i_{j}z^j + E^i_{A}y^A$ where the $y^A$ are linearily
     independent elements of $\cala$ which are independent of the $z^i$
     and $C^i_{j}, E^i_{A}$ are complex numbers. This implies in
     particular that $z^{i'\ast}=L^{i'}_{i}C^i_{j}z^j+L^{i'}_{i}
     E^i_{A}y^A$. By expanding Equation (\ref{aad}) in terms of the
     linearily independent elements $z^i\otimes z^j$, $z^i\otimes
     y^A$, $y^A\otimes z^i$ of $\cala\otimes \cala$ one obtains
     \begin{equation}
     (\bbbone + L^\ast L)C=((\bbbone +L^\ast L)C)^t
     \label{aae}
     \end{equation}
     for the complex matrices $L=(L^{i'}_{j})$ and $C=(C^i_{j})$ ($C$ is
     a square matrix whereas $L$ is generally rectangular) and
     \[
     (\bbbone + L^\ast L)E_{A}=0
     \]
     for the $E^i_{A}$ which implies $E^i_{A}=0$ (since $\bbbone+L^\ast
     L>0$). Thus one has $z^{i\ast}=C^i_{j}z^j$ which implies $\bar C
     C=\bbbone$ for the matrix $C$, $z^{i'}=\bar L^{i'}_{i}z^i$,
     $z^{i'\ast}=L^{i'}_{i}C^i_{j}z^j$ together with Equation (\ref{aae}).
     This implies $z^{\mu\ast}=\Lambda^\mu_{\nu}z^\nu$ together with
     $\Lambda^\mu_{\nu}=\Lambda^\nu_{\mu}$ for $\Lambda\in M_{4}(\mathbb
     C)$ given by
     \[
     \begin{array}{lll}
      \Lambda^i_{j} & = &C^i_{j}-\sum_{n'}C^m_{i}L^{n'}_{m}\bar
      L^{n'}_{j}\\
      \\
      \Lambda^{i'}_{j}& = & L^{i'}_{m} C^m_{j} =\Lambda^j_{i'}\\
      \\
      \Lambda^{i'}_{j'} & = & 0
      \end{array}
      \]
      With an obvious relabelling of the $z^\mu$, one can write $\Lambda$ in
      block from
      \[
      \Lambda = \left (
      \begin{array}{ccc}
       C-C^t L ^t\bar L & & C^t L^t\\
       \\
       LC & &0
       \end{array}
       \right)
       \]

\noindent The equality $\Lambda z=z^\ast$ and the symmetry of $\Lambda$ show
     that $\Lambda^\ast z^\ast=z$ so that $\Lambda^\ast\Lambda z=z$.
     Let $\Lambda=U\vert \Lambda \vert$ be the polar decomposition of
     $\Lambda$. Since $\Lambda$ is symmetric, the matrix $U$ is also
     symmetric (symmetry means $\Lambda^\ast=J\Lambda J^{-1}$ where $J$
     is the antilinear involution defining the complex structure, one
     has $\Lambda=\vert \Lambda^\ast\vert U$ so that
     $\Lambda^\ast=U^\ast \vert \Lambda^\ast\vert$ and the uniqueness of
     the polar decomposition gives $U^\ast = JUJ^{-1}$). Moreover the
     equality $\Lambda^\ast \Lambda z=z$ shows that \\

     (1) $\Lambda z=Uz$, $Pz=0$ where $P=(1-U^\ast U)$\\

\noindent One has $(1-UU^\ast)=JPJ^{-1}$ and with $e_{j}$ an orthonormal
      basis of $P\mathbb C^4$, $f_{j}=Je_{j}$ the corresponding
      orthornormal basis of $JP\mathbb C^4$ one checks that the
      following matrix is symmetric,\\

      (2) $S = \sum \vert f_{j}\rangle\langle e_{j}\vert$\\

\noindent Let now $\tilde \Lambda=U+S$. By (1) one has $\tilde \Lambda
      z=z^\ast$ since $Sz=0$ and $Uz=\Lambda z = z^\ast$. Since
      $\tilde \Lambda$ is symmetric and unitary  we get the
      conclusion.\\

\noindent Under the transformation (\ref{tr2}), $\Lambda$ transforms as
        \[
        \Lambda\mapsto u^{2 }\  ^tS\Lambda S
        \]
        so one can diagonalize the symmetric unitary $\Lambda$ by a
        real rotation $S$ and fix its first  eigenvalue to be 1 by
        chosing the appropriate $u\in U(1)$ which
      shows that one can take $\Lambda$ in diagonal form
   \begin{equation}
    \Lambda=\left (
    \begin{array}{cccc}
     1 & & &\\
     & e^{-2i\varphi_{1}} &&\\
     && e^{-2i\varphi_{2}}&\\
     &&& e^{-2i\varphi_{k}}
     \end{array}
     \right)
     \label{aaf}
     \end{equation}
     i. e. one can assume that $z^0=x^0$ and
     $z^k=e^{i\varphi_{k}}x^k$ with $e^{i\varphi_{k}}\in U(1)\subset
     \mathbb C$ for $k\in\{1, 2, 3\}$ and $x^{\mu\ast}=x^\mu\> (\in
     \cala)$ for $\mu\in\{0,\cdots,3\}$.\\

\noindent Setting $z^0=x^0=x^{0\ast}$ and $z^k=e^{i\varphi_{k}}x^k$,
$x^k=x^{k\ast}$ relations (\ref{aaa}) and (\ref{aab}) read
\begin{eqnarray}
   \cos(\varphi_{k})[x^0,x^k]_{-}& = &i\sin
   (\varphi_{\ell}-\varphi_{m})[x^\ell,x^m]_{+}\label{aag}\\
   \cos(\varphi_{\ell}-\varphi_{m})[x^\ell,x^m]_{-}&= &-i\sin
   (\varphi_{k})[x^0,x^k]_{+}\label{aah}
    \end{eqnarray}
    for $k=1, 2, 3$ where $(k,\ell,m)$ is the cyclic permutation of ($1, 2,
    3$) starting with $k$ and where $[x,y]_{\pm}=xy\pm yx$. Let
    $\mathbf{u}$ be the element $(e^{i\varphi_{1}}, e^{i\varphi_{2}},
    e^{i\varphi_{3}})$ of $T^3$, we denote by $\cala_{\mathbf{u}}$ the
    complex unital $\ast$-algebra generated by four hermitian elements
    $x^\mu$, $\mu\in\{0,\cdots, 3\}$, with relations (\ref{aag}),
    (\ref{aah}) above. It follows from the above discussion that all  $\cala$
    satisfying (I) and (II) are quotient of $\cala_{\mathbf{u}}$ for
    some $\mathbf{u}$. However it is straightforward that the pair
    $(\cala_{\mathbf{u}},U_{\mathbf{u}})$ with $U_{\mathbf{u}}
    =\bbbone_{2}x^0+i\sum^3_{k=1}e^{i\varphi_{k}}\sigma_{k}x^k$
    satisfies (I) and (II) so the $\cala_{\mathbf{u}}$ are the maximal
    solutions of (I), (II) and any other solution is a quotient of some
$\cala_{\mathbf{u}}$. In particular each maximal solution of
(I$_{1}$), (II) is the quotient $\cala^{(1)}_{\us}$ of
$\cala_{\mathbf{u}}$ by the ideal
generated by $\sum_{\mu}(x^\mu)^2-\bbbone$ for some $\mathbf{u}$.
This quotient does not contain other relations since
$\sum_{\mu}(x^\mu)^2$ is in the
center of $\cala_{\mathbf{u}}$ (Lemma \ref{lemu}). In summary one has
the following theorem.

   \begin{theo}\label{theobasic}
    $(i)$ For any $\us\in T^3$ the complex unital $\ast$-algebra
    $\cala_{\us}$ satisfies conditions $(I)$ et $(II)$.  Moreover,
    if $\cala$ is a complex unital $\ast$-algebra satisfying
    conditions $(I)$ and $(II)$ then $\cala$ is a quotient of
    $\cala_{\us}$ (i.e. a homomorphic image of $\cala_{\us}$) for
    some $\us\in T^3$.\\
   \phantom{THEOREM 1 } $(ii)$ For any $\us\in T^3$, the complex unital
$\ast$-algebra
   $\cala^{(1)}_{\us}$ satisfies  conditions $(I_{1})$ and $(II)$.
   Moreover, if $\cala^{(1)}$ is a complex unital $\ast$-algebra
   satisfying conditions $(I_{1})$ and $(II)$ then $\cala^{(1)}$ is
   a quotient of $\cala^{(1)}_{\us}$ for some $\us\in T^3$.
   \end{theo}

\noindent By construction the algebras  $\cala^{(1)}_{\us}$ are all
quotients of the universal Grassmannian ${\cal A}$ generated by $(I_{1})$ i. e.
by the matrix components  $x_1 , \ldots , x_4$
of a two by two unitary matrix.

\noindent One can show that the intersection $\cal J$ of the
kernels of the representations $\rho$ of ${\cal A}$ such that
${\rm ch}_{\frac{1}{2}} (\rho (U)) = 0 $ is a non-trivial two sided ideal
of ${\cal A}$. More precisely  let
$\mu = [x_1 , \ldots , x_4]$, be the multiple commutator
$ \Sigma \ \varepsilon (\sigma) \, x_{\sigma
(1)} \ldots x_{\sigma (4)}$ (where the sum is over all permutations 
and $\varepsilon (\sigma)$ is
the signature of the permutation)
then $[\mu , \mu^*] \ne 0 \, $ in ${\cal A}$
and $[\mu , \mu^*]$ belongs to $\cal J$
(see the Appendix for the detailed proof).
Thus the odd Grassmannian $\cal B$
which was introduced in \cite{connes:08} is a nontrivial quotient of
${\cal A}$.\\

\noindent There is another way to write relations (\ref{aag}) and
   (\ref{aah}) which will be useful for the description of the
   suspension below, it is given in the following lemma.

\begin{lemma}\label{lemaa}

   Let $\gamma_{\mu}=\gamma^\ast_{\mu}\in M_{4}(\mathbb C)$ be the
   generators of the Clifford algebra of $\mathbb R^4$, that is
   $\gamma_{\mu }\gamma_{\nu}+\gamma_{\nu}\gamma_{\mu}=2
   \delta_{\mu\nu}\bbbone$, and let $\tilde \gamma_{\mu}$ be defined
   by $\tilde \gamma_{0}=\gamma_{0}$, $\tilde \gamma_{k}=e^{i \frac{1}{2}
   \varphi_{k}\gamma}\gamma_{k}e^{-i\frac{1}{2}\varphi_{k}\gamma}$ for
   $k\in\{1,2,3\}$ with
   $\gamma=\gamma_{0}\gamma_{1}\gamma_{2}\gamma_{3}$ ($=\gamma_{5}$).
   Then the relations (\ref{aag}) and (\ref{aah}) defining $\cala_{\us}$
   are equivalent to the relation
   \[
   (\tilde \gamma_{\mu}x^\mu)^2=\bbbone \otimes \sum_{\mu}(x^\mu)^2
   \]
   in $M_{4}(\cala_{\us})=M_{4}(\mathbb C)\otimes \cala_{\us}$.
   \end{lemma}

\noindent This is easy to check using $\gamma\gamma_{\mu}=-\gamma_{\mu}\gamma$
   and $\gamma^2=\bbbone$. On the right-hand side of the above equality
   appears the central element $\sum_{\mu}(x^\mu)^2$ of $\cala_{\us}$;
   the algebra $\cala_{\us}$ has another central element described in
   the following lemma.

\begin{lemma}\label{lema}
   The element $\sum^3_{k=1}\cos(\varphi_{k}-\varphi_{\ell}-\varphi_{m})
    \cos(\varphi_{k})\sin(\varphi_{k})(x^k)^2$ is in the center of
    $\cala_{\mathbf{u}}$, where in the summation $(k, \ell, m)$ is the
    cyclic permutation of $(1, 2, 3)$ starting with $k$ for $k\in\{1, 2,
    3\}$.
    \end{lemma}
    This can be checked directly  using (\ref{aag}), (\ref{aah}). So
   one has two quadratic elements in the $x^\mu$ which belong to the
   center $Z( \cala_{\mathbf{u}})$ of $\cala_{\mathbf{u}}$. In fact,
   for generic $\us$, the center is generated by these two quadratic
   elements.\\

\noindent By changing $x_{k}$ in $-x_{k}$ one can replace $\varphi_{k}$ by
   $\varphi_{k}+\pi$ and by a rotation of $SO(3)$ one can permute the
   $\varphi_{k}$ without changing the algebra $\cala_{\us}$ nor the
   algebra $\cala^{(1)}_{\us}$. It follows that it is sufficient to
   take $\us$ in the 3-cell defined by
   \begin{equation}
   \{ (e^{i\varphi_{k}})\in T^3 \vert \pi > \varphi_{1}\geq
   \varphi_{2}\geq \varphi_{3}\geq 0\}
   \label{eqregi}
   \end{equation}
   to cover all the $\cala_{\us}$ and $\cala^{(1)}_{\us}$.\\

\noindent It is apparent that $\cala_{\mathbf{u}}$ is a deformation of the
   commutative $\ast$-algebra $\calg(\mathbb R^4)$ of complex
   polynomial functions  on $\mathbb R^4$; it reduces to the latter for
   $\varphi_{1}=\varphi_{2}=\varphi_{3}=0$ that is for
$\mathbf{u}=\mathbf{e}$ where
   $\mathbf{e}=(1,1,1)$ is the unit of $T^3$. We shall denote
   $\cala_{\mathbf{u}}$ by $\calg(\mathbb R^4_{\mathbf{u}})$ defining
   thereby the noncommutative 4-plane $\mathbb R^4_{\mathbf{u}}$ as
   dual object. Similarily, the quotient $\cala^{(1)}_{\us}$ of
$\cala_{\mathbf{u}}$ by the
   ideal generated by $\sum_{\mu}(x^\mu)^2-\bbbone$ is a deformation of
   the $\ast$-algebra $\calg(S^3)$ of polynomial functions on $S^3$ that
   is of functions on $S^3$ which are restrictions to $S^3\subset
   \mathbb R^4$ of elements of $\calg(\mathbb R^4)$; we shall denote
   this quotient $\cala^{(1)}_{\us}$ by $\calg(S^3_{\mathbf{u}})$
defining thereby the
   noncommutative 3-sphere $S^3_{\mathbf{u}}$ by duality. \\

\noindent Let $\calg(\mathbb R^5_{\us})$ be the unital $\ast$-algebra obtained
   by adding a central hermitian generator $x^4$ to $\calg(\mathbb
   R^4_{\us})=\cala_{\us}$, i.e. $\calg(\mathbb R^5_{\us})$ is the
   unital $\ast$-algebra generated by hermitian elements $x^\mu$,
   $\mu\in\{ 0,\dots,3\}$, and $x^4$ such that the $x^\mu$ satisfy
   (\ref{aag}), (\ref{aah}) and that one has $x^\mu x^4=x^4x^\mu$
   for $\mu\in\{0,\dots,3\}$; the noncommutative 5-plane $\mathbb
   R^5_{\us}$ being defined by duality. Let $\calg(S^4_{\us})$
   be the unital $\ast$-algebra quotient of $\calg(\mathbb R^5_{\us})$
   by two-sided ideal generated by the hermitian central element
   $\sum^3_{\mu=0}(x^\mu)^2+(x^4)^2-\bbbone$. The noncommutative
   4-sphere $S^4_{\us}$ defined as dual object is in the obvious sense
   the {\sl suspension} of $S^3_{\us}$. This is a 3-parameter deformation of
   the sphere $S^4$ which reduces to $S^4_{\theta}$ for
   $\varphi_{1}=\varphi_{2}=-\frac{1}{2}\theta$ and $\varphi_{3}=0$,
   (see below). We denote by $u^\mu,u$ the canonical images of
   $x^\mu,x^4\in  \calg(\mathbb R^5_{\us})$ in $\calg(S^4_{\us})$ and by
   $v^\mu$ the canonical images of $x^\mu\in\calg(\mathbb R^4_{\us})$
   in $\calg(S^3_{\us})$, i.e. one has $\sum(u^\mu)^2+u^2=\bbbone$ and
   $\sum(v^\mu)^2=\bbbone$, etc.. It will be convenient for further
   purpose to summarize some important points discussed above by the
   following theorem.

   \begin{theo} \label{theo01}
    $(i)$ One obtains a hermitian projection
   $e\in M_{4}(\calg(S^4_{\us} ))$ by setting $e=\frac{1}{2}(\bbbone +
   \tilde\gamma_{\mu}u^\mu+\gamma u)$. Furthermore one has
   $\ch_{0}(e)=0$ and $\ch_{1}(e)=0$.\\
   \phantom {THEOREM 2 } $(ii)$  One obtains a unitary $U\in
   M_{2}(\calg(S^3_{\us}))$ by setting $U=\bbbone v^0 + i\tilde\sigma_{k }
   v^k$ where $\tilde\sigma_{k}=\sigma_{k} e^{i\varphi_{k}}$.
   Furthermore one has $\ch_{\frac{1}{2}}(U)=0$.
   \end{theo}

\noindent Statement $(ii)$  is just a reformulation of what has be
   done previously. Concerning Statement $(i)$, the fact that $e$
   is a hermitian projection with $\ch_{0}(e)=0$  follows directly from
   the definition and Lemma  \ref{lemaa}  whereas $\ch_{1}(e)=0$ is a
   consequence of $\ch_{\frac{1}{2}} (U)=0$ in $(ii)$.\\

\noindent We shall now compute
   $\tilde{\ch}_{\frac{3}{2}}(U)$ and check that,
except for exceptional values of $\us$ for which
   $\tilde\ch_{\frac{3}{2}}(U)=0$,
  it is a non trivial Hochschild
   cycle on $\cala_{\mathbf{u}}$.\\

\noindent One has by construction
   \[
   \tilde{\ch}_{\frac{3}{2}}(U_{\mathbf{u}})=\tr(U_{\mathbf{u}}
   \circledcirc U_{\mathbf{u}}^\ast \circledcirc
   U_{\mathbf{u}}\circledcirc
   U_{\mathbf{u}}^\ast- U_{\mathbf{u}}^\ast\circledcirc
   U_{\mathbf{u}} \circledcirc U_{\mathbf{u}}^\ast \circledcirc U_{\mathbf{u}})
   \]
   which is an element of $\cala_{\mathbf{u}} \otimes
   \cala_{\mathbf{u}}\otimes \cala_{\mathbf{u}}\otimes \cala_{\mathbf{u}}$
   and can be considered as a $\cala_{\mathbf{u}}$-valued Hochschild
   3-chain. One obtains using (\ref{aag}), (\ref{aah})
   \begin{eqnarray}
    \tilde{\ch}_{\frac{3}{2}}(U_{\mathbf{u}}) & = & - \sum_{3\geq \alpha,
    \beta,\gamma,\delta\geq 0}\epsilon_{\alpha\beta\gamma\delta}
    \cos(\varphi_{\alpha}-\varphi_{\beta}+\varphi_{\gamma}-\varphi_{\delta})
    x^\alpha\otimes x^\beta \otimes x^\gamma \otimes x^\delta\nonumber\\
    & & +\  i\sum_{3\geq\mu,\nu\geq 0}\sin(2(\varphi_{\mu}-\varphi_{\nu}))
    x^\mu\otimes x^\nu \otimes x^\mu \otimes x^\nu\label{aan}
    \end{eqnarray}
    where $\epsilon_{\alpha\beta\gamma\delta}$ is completely antisymmetric
    with $\epsilon_{0123}=1$ and where we have set $\varphi_{0}=0$. Using
    (\ref{aan}), (\ref{aag}), (\ref{aah}) one checks that
    $\tilde{\ch}_{\frac{3}{2}}(U_{\mathbf{u}})$ is in fact a Hochschild
    cycle, i.e. $b(\tilde{\ch}_{\frac{3}{2}}(U_{\mathbf{u}}))=0$.
    Actually, this follows on general grounds from the fact that
    $\tilde{\ch}_{\frac{1}{2}}(U_{\mathbf{u}})=0$ and that
    $U_{\mathbf{u}}^\ast U_{\mathbf{u}}=U_{\mathbf{u}}U_{\mathbf{u}}^\ast$
    is an element of the center $\bbbone_{2}\otimes
    Z(\cala_{\mathbf{u}})$ of $M_{2}(\cala_{\mathbf{u}})$ in view of
   Lemma \ref{lemu}. In fact the $\cala_{\mathbf{u}}$-valued Hochschild 3-cycle
    $\tilde{\ch}_{\frac{3}{2}}(U_{\mathbf{u}})$ is trivial (i.e. is a
    boundary) if and only if it vanishes (which means that all
    coefficients vanish in formula (\ref{aan})).
    Indeed $\cala_{\mathbf{u}}$ is a $\mathbb N$-graded algebra with
    $\cala_{\mathbf{u}}^0=\mathbb C\bbbone$ and $\cala^1_{\mathbf{u}}$=
    linear span of the $\{x^\mu\vert\mu\in\{0,\cdots,3\}\}$ and the
    Hochschild boundary preserves the degree. It follows that
    $\tilde{\ch}_{\frac{3}{2}}(U_{\mathbf{u}})$ can only be the boundary of
    linear combinations of terms which are in
    $\otimes^5\cala_{\mathbf{u}}$ of total degree 4 and contain
    therefore at least one tensor factor equal to $\bbbone$. Among
    these terms, the $\bbbone\otimes
    x^\alpha\otimes x^\beta \otimes x^\gamma\otimes x^\delta$ are the only
    ones which contain in
    their boundaries tensor products of four $x^\mu$. One has for these
    terms
    \[
    \begin{array}{l}
     b(\bbbone\otimes x^\alpha\otimes x^\beta\otimes x^\gamma \otimes
     x^\delta)  =  x^\alpha\otimes x^\beta\otimes x^\gamma \otimes
     x^\delta + x^\delta\otimes x^\alpha \otimes x^\beta \otimes
     x^\gamma\\
     \\
   - \bbbone\otimes (x^\alpha x^\beta \otimes x^\gamma \otimes
     x^\delta - x^\alpha\otimes x^\beta x^\gamma \otimes
     x^\delta + x^\alpha\otimes x^\beta \otimes x^\gamma x^\delta)
     \end{array}
     \]
     however the $x^\alpha\otimes x^\beta\otimes x^\gamma\otimes
     x^\delta + x^\delta \otimes x^\alpha \otimes x^\beta \otimes
     x^\gamma$ cannot produce by linear combination a term with the
     kind of generalized antisymmetry of
     $\tilde{\ch}_{\frac{3}{2}}(U_{\mathbf{u}})$ excepted of course if
    $\tilde{\ch}_{\frac{3}{2}}(U_{\mathbf{u}})=0$. Thus
    $\tilde{\ch}_{\frac{3}{2}}(U_{\mathbf{u}})$ is non trivial if not
    zero.

\noindent The $\cala^{(1)}_{\us}$-valued Hochschild 3-cycle
    $\ch_{\frac{3}{2}}(U)$ on $\cala^{(1)}_{\us}$ is the image of
    $\tilde\ch_{\frac{3}{2}}(U_{\us})$ by the projection of
    $\cala_{\us}$ onto $\cala^{(1)}_{\us}$. In particular
    $\ch_{\frac{3}{2}}(U)$ vanishes if
    $\tilde\ch_{\frac{3}{2}}(U_{\us})$ vanishes which occurs on
    $\Sigma^3$ for $\varphi_{1}=\varphi_{2}=\varphi_{3}=\frac{\pi}{2}$
    and for $\varphi_{1}=\frac{\pi}{2}$, $\varphi_{2}=\varphi_{3}=0$.
    For these two values of $\us$, the algebras $\cala_{\us}$ are
    isomorphic, one passes from
    $\varphi_{1}=\varphi_{2}=\varphi_{3}=\frac{\pi}{2}$ to
    $\varphi_{1}=\frac{\pi}{2}$, $\varphi_{2}=\varphi_{3}=0$ by the
    exchange of $x^0$ and $x^1$; this is of course the same for
    $\cala^{(1)}_{\us}$. One can furthermore check that the Hochschild
    dimension of $\cala^{(1)}_{\us}$ for these values of $\us$ is one.

\noindent To obtain the Hochschild 4-cycle on  $\cala_{\us}$ corresponding to
the volume form on the noncommutative 4-plane $\mathbb R^4_{\mathbf{u}}$,
we shall just apply   to $\tilde
    \ch_{\frac{3}{2}}(U_{\us})$ the natural extension of the
   de Rham coboundary in the noncommutative case, namely the
    operator $B:\cala_{\us}\otimes \tilde
    \cala^{\otimes^3}_{\us}\rightarrow \cala_{\us}\otimes
    \tilde\cala^{\otimes^4}_{\us}$ (\cite{connes:02} \cite{jll}). Since $\tilde
    \ch_{\frac{3}{2}}(U_{\us})$ is not only a Hochschild cycle but
    also fulfills the cyclicity condition, it follows that, up to an
    irrelevant normalization $B$ reduces there to the tensor product
    by $\bbbone$, thus
    \[
    B \tilde \ch_{\frac{3}{2}}(U_{\us})=\bbbone \otimes \tilde
    \ch_{\frac{3}{2}}(U_{\us})
    \]
    and the Hochschild 4-cycle $B\tilde \ch_{\frac{3}{2}}(U_{\us})$
     which plays the role of the volume form of $\mathbb R^4_{\us}$
    is thus given by
\begin{eqnarray}
    v & = & - \sum_{3\geq \alpha,
    \beta,\gamma,\delta\geq 0}\epsilon_{\alpha\beta\gamma\delta}
    \cos(\varphi_{\alpha}-\varphi_{\beta}+\varphi_{\gamma}-\varphi_{\delta})
    \bbbone \otimes x^\alpha\otimes x^\beta \otimes x^\gamma \otimes
x^\delta\nonumber\\
    & & +\  i\sum_{3\geq\mu,\nu\geq 0}\sin(2(\varphi_{\mu}-\varphi_{\nu}))
    \bbbone \otimes x^\mu\otimes x^\nu \otimes x^\mu \otimes x^\nu\label{aan2}
    \end{eqnarray}

\noindent It turns out
    that this 4-cycle is non trivial whenever it does not vanish as
    can be verified by evaluation at the origin which is the classical
    point of $\mathbb R^4_{\us}$. The nontriviality of
    $\ch_{\frac{3}{2}}(U)$ follows since $B\tilde
    \ch_{\frac{3}{2}}(U_{\us})$ is its suspension.

   \section{The Scaling Foliation and relation to Sklyanin algebras} 
\label{sec03}
  \setcounter{equation}{0}

We let as above $\Sigma=T^3$ be the parameter space for 3-dimensional
spherical manifolds $S^3_{\mathbf{u}}$. Different $S^3_{\mathbf{u}}$
can span isomorphic 4-dimensional $\mathbb R^4_{\mathbf{u}}$ and we
shall analyse here the corresponding foliation of $\Sigma$.\\

\noindent More precisely, let us say that $S^3_{\mathbf{u}}$ is 
"scale-equivalent"
to $S^3_{\mathbf{v}}$ and write $u \sim v $ when the quadratic algebras
corresponding to $\mathbb R^4_{\mathbf{u}}$ and $\mathbb 
R^4_{\mathbf{v}}$ are isomorphic.
  This generates a foliation of $\Sigma$ which
is completely described by the orbits of the flow of the following 
vector field:\\
\begin{equation}
Z \,=\,\sum^3_{k=1}\, \sin(2 \varphi_{k}) \, 
\sin(\varphi_{\ell}+\varphi_{m}-\varphi_{k})
  \,\,  \frac{\partial}{\partial \varphi_{k}}
\label{3a}
\end{equation}

\noindent as shown by,

\begin{theo}\label{theofeuille}
Let $\mathbf{u} \in \Sigma$. There exists a neighborhood $V$ of 
$\mathbf{u}$ such that
$\mathbf{v}\in V$ is scale-equivalent to $\mathbf{u}$ if and only if 
it belongs to
the orbit of $\mathbf{u}$ under the flow of $Z$.
   \end{theo}

   \noindent Let us first show that if $\mathbf{v}$ belongs to the
orbit of $\mathbf{u}$ under $Z$ then the corresponding quadratic algebras
are isomorphic. To the action of the group of permutations ${\cals}_4$ of the
4 generators of the quadratic algebra there corresponds an action of $\cals_4$
on the parameter space $\Sigma$.
This action is the obvious one on the subgroup $\cals_3$ of permutations 
fixing $0$
and the action of the permutation $(1,0,3,2)$ of $(0,1,2,3)$ is given by the
following transformation,

\begin{equation}
w( \varphi_{1},\varphi_{2},\varphi_{3}) \,=
(- \varphi_{1},\varphi_{3}-\varphi_{1},\varphi_{2}-\varphi_{1})
\label{3b}
\end{equation}

\noindent The transformation $w$ and its conjugates under the action of
$\cals_3$ by permutations of the $\varphi_{j}$ generate an abelian group $K$
of order 4 which is a normal subgroup of the group $W=\cals_4$ 
generated by $w$ and
$\cals_3$. By construction $g(\mathbf{u})$ is scale-equivalent to 
$\mathbf{u}$ for
any $g \in W$. At a more conceptual level the group
$W$ is the Weyl group of the symmetric space used in lemma 2, of 
symmetric unitary
(unimodular) 4 by 4 matrices. Moreover the flow of $Z$ is invariant 
under the action
of $W$. This is obvious for $g \in  \cals_3$ and can be checked 
directly for $w$.

\noindent Let $C$ be the set of critical points for $Z$, i.e. $C=\{ 
\mathbf{u}, Z_{\mathbf{u}}=0 \}$.
For $\mathbf{u} \in C$ the orbit of $\mathbf{u}$ is reduced to 
$\mathbf{u}$ and the
required equivalence is trivial.
To handle the case $\mathbf{u} \notin C$ we let $D \subset \Sigma$ be 
the zero set of the function,

\begin{equation}
\delta(\mathbf{u})= \prod^3_{k=1}\,
\sin(\varphi_{k})\,\cos(\varphi_{l}-\varphi_{m})
\label{3c}
\end{equation}

\noindent The inclusion $\cap \, gD \, \subset C$ where $g$ varies in $K$
shows that we can assume that $\mathbf{u} \notin D$. We can then find 4
non-zero scalars $s^\mu$, $\mu\in\{0,\cdots,3\}$ such that,

\begin{eqnarray}
      s^0s^1\cos(\varphi_{2}-\varphi_{3})+ s^2s^3\sin(\varphi_{1}) & = &
      0\nonumber\\
      s^0s^2\cos(\varphi_{3}-\varphi_{1})+ s^3s^1\sin(\varphi_{2}) & = &
      0\label{aal}\\
      s^0s^3\cos(\varphi_{1}-\varphi_{2})+ s^1s^2\sin(\varphi_{3}) & = &
      0\nonumber
      \end{eqnarray}

\noindent The solution is unique
(up to an overall normalization and choices of sign) and can be
written in the form,
\[
\begin{array}{lll}
  s^0& = & (\prod_{j} \sin \varphi_{j})^{1/2}\\
  \\
  s^k & = &(\sin \varphi_{k}\prod_{\ell\not=
  k}\cos(\varphi_{k}-\varphi_{\ell}))^{1/2}
  \end{array}
  \]
 where the square roots are chosen so that $\prod s^\mu =-\delta(\mathbf{u})$.
\noindent Then, provided that $\cos(\varphi_{j}) \neq 0 \, \,\forall j$,
the relations (\ref{aag}),  (\ref{aah}) can be written
  \begin{eqnarray}
    [S_{0},S_{k}]_{-} & = & iJ_{\ell m}[S_{\ell},S_{m}]_{+}\label{aai}\\
    {[S_{\ell},S_{m}]}_{-} & = & i[S_{0},S_{k}]_{+}\label{aaj}
    \end{eqnarray}
    where $J_{\ell m}= -\tan(\varphi_{\ell}-\varphi_{m}) \tan(\varphi_{k})$
    for any cyclic permutation $(k,\ell, m)$ of $(1,2,3)$ and where
    \begin{equation}
     S_{\mu}=s^\mu x^\mu
     \label{aak}
     \end{equation}

   \noindent    So defined the three real numbers $J_{k\ell}$ satisfy 
the relation
     \begin{equation}
      J_{12}+J_{23}+J_{31}+J_{12}J_{23}J_{31}=0
      \label{aam}
      \end{equation}
      as easily verified. The relations (\ref{aai}), (\ref{aaj})
      together with (\ref{aam}) for the constants $J_{k\ell}$
      characterize the algebra introduced by Sklyanin in connection
      with the Yang-Baxter equation \cite{Skly:1}, \cite{Skly:2}.\\

   \noindent  In the case when the $s^\mu$ are real, the 
transformation (\ref{aak})
  preserves the involution which on the Sklyanin algebra $S(J_{k\ell})
  $ is given by
  \[
  S^\ast_{\mu}=S_{\mu}\qquad  \mu=0,1,2,3.
  \]
  In general, however, one cannot choose the $s^{\mu}$'s to be real and
  the involutive algebra $\cala_{\us}$ gives a different real form of
  the Sklyanin algebra.\\

  \noindent The invariance of the $J_{k\ell}$ under the flow $Z$, 
$Z(J_{k\ell})\, =\, 0$,
thus gives the required scale-equivalence on the orbit of 
$\mathbf{u}$ provided
$\cos(\varphi_{j}) \neq 0 \, \,\forall j$.
\noindent The condition $\varphi_{j}=\pi /2$
is invariant under the flow $Z$ and this special
case is handled in the same way
(note that if moreover $\varphi_{l}=\varphi_{m}$
  one of the relations
becomes trivial, the corresponding algebra is not
  a Sklyanin algebra but is constant
on the orbit of $Z$).

\noindent We have thus shown that two points on the
  same orbit of $Z$ are scale-equivalent.
Let us now prove the converse in the form stated in
  theorem 3. In order to distinguish the
quadratic algebras $\cala_{\us}$ we shall use an
  invariant called the associated geometric data.

\noindent The Sklyanin algebras $S(J_{k\ell})$
have been extensively studied
  from the point of view of noncommutative algebraic geometry. An
  important role is played by the associated geometric data
  \[
  \{ E,\sigma, \call \}
  \]
   consisting of an elliptic curve $E\subset P_{3}(\mathbb C)$, an
   automorphism $\sigma$ of $E$ and an invertible $\calo_{E}$-module
   $\call$ (cf.
   \cite{art:1}, \cite{art:2}, \cite{Od:fe}, \cite{sm:st}). This
   geometric data is invariantly defined for any graded algebra and in
   the above case of $S(J_{k\ell})$, it degenerates when one of the
   parameters $J_{k\ell}$ vanishes (or in the case $J_{k\ell}=1,
   J_{\ell r}=-1$, cf. \cite{sm:st} for a careful discussion).\\
\noindent It is straightforward to extend the computations of  \cite{sm:st}
to the present situation in order to cover all cases. Up to the action of the
group $W$ the critical set $C$ is the union of the point 
$P=(\pi/2,\pi/2,\pi/2)$
with the two circles,
  \[
C_+=\{\mathbf{u} \,;\varphi_{1}=\varphi_{2}, \varphi_{3}=0 \} \,, \, \, \,
C_-=\{\mathbf{u} \,; 
\varphi_1=\frac{\pi}{2}+\varphi_3,\varphi_2=\frac{\pi}{2} \}
\]
\noindent For $\mathbf{u}=P$, the geometric data is very degenerate,
$E= P_{3}(\mathbb C)$, while $\sigma$ is a symmetry of determinant $-1$.
In fact there are two other $W$-orbits, those of $P'=(\pi/2,\pi/2,0)$ and
of $O=(0,0,0)$ for which $E= P_{3}(\mathbb C)$. For $P'$, the correspondence 
$\sigma$ is a symmetry of determinant $1$, while for $O$ it is the identity. 
 
\noindent For $\mathbf{u}\in C_+$, $\mathbf{u}\neq O$, $\mathbf{u}\notin W(P')$ 
the geometric data degenerates to the union of 6
    projective lines $P_{1}(\mathbb C)$, with $\sigma$ given by
    multiplication by $1$ for two of them, by $e^{2i\varphi_{1}}$ for two
others and $e^{-2i\varphi_{1}}$ for the last two.
The case $\mathbf{u}\in C_-$ is similar, but not identical.
$E$ is the union of six lines but $\sigma$ is given by
    multiplication by $-1$ for two of them, it exchanges two of the 
remaining lines
with $\sigma^2$ given by multiplication by $e^{4i\varphi_{1}}$
and exchanges the last two
with $\sigma^2$ given by multiplication by $e^{-4i\varphi_{1}}$.
  \\

\noindent For $\mathbf{u}\notin C$, we can assume as above that
$\mathbf{u}\notin D$. Then, provided that $\cos(\varphi_{j}) \neq 0 
\, \,\forall j$
we can reduce as above to Sklyanin algebras. In that case (\cite{sm:st})
the geometric data  $E\subset P_{3}(\mathbb C)$ is the union of 4 points
with a non-singular elliptic curve, except (up to signed 
permutations) for the following degenerate case:
  
  \[
  F_1=\{\mathbf{u};  \, J_{23}=-a, J_{31}=a, J_{12}=0 \}
\]
 In that case, $E$ is the union of 2 points, one line and 2 
circles, the correspondence
$\sigma$ fixes the 2 points and the line pointwise. It restricts to 
both circles
$\Gamma_j \sim P_{1}(\mathbb C)$ and is given in terms of a rational parameter as 
the multiplication
by\linebreak[4] $(i + a^{1/2})/(i - a^{1/2})$ where each circle 
corresponds to a different choice
of the square root $a^{1/2}$.

\noindent In the case

\[
   F_2=\{\mathbf{u};  \, \varphi_{1}= \frac{\pi}{2}   \, \, , \, \, 
\varphi_{2} \neq  \varphi_{3},\> \varphi_2\not =\frac{\pi}{2},\> 
\varphi_3\not =\frac{\pi}{2}  \}
\]
  where $\mathbf{u}\notin D$ but  $\cos(\varphi_{j})= 0$ for some $j$
say j=1, the above change of variables breaks down, but the direct 
computation shows that
as for $\mathbf{u} \in F_1$, $E$ is the union of 2 points, one line and 2 circles. 
However the correspondence
$\sigma$ is different from that case. It fixes the 2 points and is 
multiplication by $-1$
on the line. It exchanges the two circles
$\Gamma_j \sim P_{1}(\mathbb C)$ and its square $\sigma^2 $ is given in terms of a rational
  parameter as the multiplication
by the square of $(i + b^{1/2})/(i - b^{1/2})$, $b=- J_{31}$,
  where each circle corresponds to a different choice
of the square root $b ^{1/2}$.\\

\noindent On the circle,

\[
   L =\{\mathbf{u};  \, \varphi_{1}= \frac{\pi}{2}   \, \, , \, \, 
\varphi_{2} =  \varphi_{3}\}
\]
the first of the six relations (\ref{aah})  becomes trivial and
the quadratic algebra is independent of the
value of $\varphi_{2} =  \varphi_{3}$ except for the isolated values 
$0$ and $\pi/2$, which
correspond to the orbit $ W( P)  $ of the point 
$P=(\pi/2,\pi/2,\pi/2)$ discussed above and for which 
$\tilde\ch_{\frac{3}{2}}(U_\us)$ vanishes as explained in the last 
section.
  For points of the circle $L$ not on this orbit, $E$ is the union of
six lines. The correspondence $\sigma$ is $1$ on one line, $-1$
on another line, and permutes cyclically the
remaining 4 lines, inducing twice an isomorphism and twice the coarse 
correspondence.

\noindent Finally on the circle,

\[
   L' =\{\mathbf{u};  \, \varphi_{1}= \frac{\pi}{2}   \, \, , \, \, 
\varphi_{2} =  \frac{\pi}{2}   \}
\]
except for the special cases treated above $E$ is the union of
a point with $P_{2}(\mathbb C)$ and the correspondence $\sigma$ is 
a symmetry of determinant $-1$. \\

\noindent Let us now end the proof of theorem 3. We work modulo $W$.
  For $\mathbf{u}\in C$  the geometric data
allows to distinguish it from any $\mathbf{v}$ in a
neighborhood (one checks this for $\mathbf{u}=P$ and $\mathbf{u}\in C_{+,-}$).
For $\mathbf{u}\notin C$ the flow line through $\mathbf{u}$ is non-trivial.
For $\mathbf{u}\in L$ or $ L'$ the nearby points having the same
geometric data are necessarily on $L$ or $L'$ and the scaling flow is
locally transitive on both, so the answer follows.
Each of the faces $F_j$ is globally invariant under the flow $Z$.
For $\mathbf{u}\in F_j$ the nearby points having the same
geometric data are necessarily on $F_j$ and the correspondence $\sigma$
gives the required information to conclude that scale-equivalent nearby points
are on the same flow line.
Finally for points not $W$-equivalent to those treated so far, the 
geometric data
is a non-degenerate elliptic curve E whose j-invariant is given by
\[
\begin{array}{lll}
   j=256 (\lambda^2 - \lambda +1)^3/(\lambda^2 (1-\lambda)^2) \\
\\ \lambda=
\sin(2 \varphi_{1})\sin(2( \varphi_{2}- \varphi_{3}))/\sin(2 
\varphi_{2})\sin(2( \varphi_{1}- \varphi_{3}))
  \end{array}
\]
and a translation $\sigma$ which together allow for the local determination
of the parameters $J_{k\ell}$ and hence of the flow line of $\mathbf{u}$.\\

  \begin{corol}\label{corol1}
	  The critical points of the scaling foliation are given by the
union of the $W$-orbits of
$P$, of $C_+$ and of $C_-$.
	  \end{corol}

\noindent We shall now  analyse the noncommutative 3-spheres
associated to the critical points in $C_+$.
The easiest way to understand them is as special cases of the general
procedure of $\theta$-deformation (applied here to the usual 3-sphere
and also to $\mathbb R^4$) which lends itself to easy
    higher dimensional generalization. (The case of $C_-$
can be reduced to $C_+$ thanks to an easy "involutive twist"
which will be described in general in part II).\\

\section{The $\theta$-deformed $2n$-plane $\mathbb
R^{2n}_{\theta}$ and its Clifford algebra} \label{sec2}
\setcounter{equation}{0}

In the previous sections, we have obtained a multiparameter
noncommutative deformation $\calg(\mathbb R^4_{\us})$ of the graded
algebra $\calg(\mathbb R^4)$ of polynomial functions on $\mathbb R^4$
which induces a corresponding deformation $\calg(S^3_{\us})$ of the
algebra of polynomial functions on $S^3$ in such a way that all
dimensions are preserved as will be shown in Part II. Moreover this
is the generic deformation
under the above conditions. We also extracted from
this multiparameter deformation of $\calg(\mathbb R^4)$ a
one-parameter deformation $\calg(\mathbb R^4_{\theta})$ of
$\calg(\mathbb R^4)$ which is also a one-parameter deformation
$\calg(\mathbb C^2_{\theta})$ of $\calg (\mathbb C^2)$ whence $\mathbb
C^2$ is identified with $\mathbb R^4$ through (for instance)
$z^1=x^0+ix^3$, $z^2=x^1+ix^2$. The parameter $\theta$ corresponds
to the curve $\theta \mapsto \us(\theta)$ defined by
$\us_{1}=\us_{2}=e^{-\frac{i}{2}\theta}$ and $\us_{3}=1$, i.e. to
$\varphi_{1}=\varphi_{2}=-\frac{1}{2}\theta$ and $\varphi_{3}=0$ in
terms of the
previous parameters. Indeed for these values of $\us$, the
relations  (\ref{aag}), (\ref{aah}) for $x^0,x^1,x^2,x^3$ read in terms
of $z^1=x^0+ix^3,\bar z^1=x^0-ix^3$, $z^2=x^1+ix^2$, $\bar
z^2=x^1-ix^2$, (one has $z^{1\ast}=\bar z^1$ and $z^{2\ast}=\bar
z^2$) $z^1z^2=\lambda z^2 z^1$, $\bar z^1\bar z^2=\lambda \bar z^2\bar
z^1$, $z^1 \bar z^1=\bar z^1 z^1$, $z^2\bar z^2=\bar z^2 z^2$,
$\bar z^1 z^2=\lambda^{-1}z^2\bar z^1$, $z^1\bar
z^2=\lambda^{-1}\bar z^2z^1$ where we have set
$\lambda=e^{i\theta}$. This one-parameter deformation is well suited
for simple higher-dimensional generalizations (i.e. $\mathbb C^2$ is replaced
by $\mathbb C^n$ and $\mathbb R^4$ by $\mathbb R^{2n}$, $n\geq 2$). In
the following we shall describe and analyze them in details.
For this we shall generalize
$\theta$ as explained at the end of the introduction as an
antisymmetric matrix $\theta\in M_{n}(\mathbb R)$, the previous one
being identified as $\left( \begin{array}{cc}
0 & \theta\\
-\theta & 0
\end{array}
\right) \in M_{2}(\mathbb R)$, and we shall use the notations
explained at the end of the introduction.

\noindent Let $\pol(\mathbb R^{2n}_{\theta})$ be the complex unital associative
algebra generated by $2n$ elements $z^\mu,\bar z^\mu$
($\mu,\nu=1,\dots,n$) with relations
\begin{equation}
z^\mu z^\nu  =  \lambda^{\mu\nu} z^\nu z^\mu,   \>
\bar z^\mu \bar z^\nu = \lambda^{\mu\nu} \bar z^\nu \bar z^\mu,
\>
\bar z^\mu z^\nu  = \lambda^{\nu\mu} z^\nu\bar z^\mu
\label{a}
\end{equation}
for $\mu,\nu=1,\dots,n$ ($\lambda^{\mu\nu}=e^{i\theta_{\mu\nu}}$,
$\theta_{\mu\nu}=-\theta_{\nu\mu}\in \mathbb R$).  Notice that one
has $\lambda^{\nu\mu} =1/\lambda^{\mu\nu}=\overline{\lambda^{\mu\nu}}$ and that
$\lambda^{\mu\mu}=1$.  We endow $\pol(\mathbb R^{2n}_\theta)$ with
the unique $\mathbb C$-algebra involution $x\mapsto x^\ast$ such that
$z^{\mu\ast}=\bar z^\mu$.  Clearly the $\ast$-algebra $\pol(\mathbb
R^{2n}_{\theta})$ is a deformation of the commutative $\ast$-algebra
$\pol(\mathbb R^{2n})$ of complex polynomial functions on $\mathbb
R^{2n}$, (it reduces to the latter for $\theta=0$).  The algebra
$\pol(\rnt)$ will be refered to as the {\sl algebra of complex
polynomials on the noncommutative $2n$-plane $\mathbb
R^{2n}_{\theta}$}. \\

\noindent In fact the relations (\ref{a}) define  a deformation 
$\mathbb C^n_{\theta}$ of
   $\mathbb C^n$ and we can identify $\mathbb C^n_{\theta}$ and $\rnt$
   by writing $\calg(\rnt)=\calg(\mathbb C^n_{\theta})$.
Correspondingly, the unital subalgebra
   $\Halg(\mathbb C^n_{\theta})$ generated by the $z^\mu$ is a
   deformation of the algebra of holomorphic polynomial functions on
   $\mathbb C^n$. \\

\noindent There is a unique group-homomorphism  $s\mapsto
   \sigma_{s}$ of the abelian group $T^n$ into the group $\aut(\calg(\rnt))$
   of unital $\ast$-automorphisms of $\calg(\rnt)$ which is such that
   $\sigma_{s}(z^\nu) = e^{2\pi i s_{\nu}}z^\nu$,  ($
   \sigma_{s}(\bar z^\nu) = e^{-2\pi i s_{\nu}}\bar z^\nu$). This
   definition is independent of $\theta$, in particular $s\mapsto
   \sigma_{s}$ is well defined as a group-homomorphism of $T^n$ into
   $\aut(\calg(\mathbb R^{2n}))$ where it is induced by a smooth action
   of $T^n$ on the manifold $\mathbb R^{2n}$.  It follows from the
relations (\ref{a}) that the $z^\mu z^{\mu\ast}=z^{\mu\ast} z^\mu$ ($1\leq
   \mu\leq n$) are in the center of $\pol(\mathbb R^{2n}_{\theta})$.
   Furthermore these hermitian elements generate the center as unital
   subalgebra of $\pol(\rnt)$ whenever $\theta$ is {\sl generic}, i.e.
   for $\theta_{\mu\nu}$ irrational $\forall \mu,\nu$ with $1\leq
   \mu<\nu\leq n$. On the other hand these elements $\bar z^\mu z^\mu$
   generate the subalgebra $\calg(\rnt)^\sigma$ of elements of
   $\calg(\rnt)$ which are invariant by the action $\sigma$ of $T^n$.
   Thus $\calg(\rnt)^\sigma$ is contained in the center of
   $\calg(\rnt)$. This is not an accident, moreover the subalgebra of
   invariant elements of $\calg(\rnt)$ is  not deformed (i.e. does not
   depend on $\theta$) and is canonically isomorphic to $\calg(\mathbb
   R^{2n})^\sigma$.\\

\noindent Let $\clif(\mathbb R^{2n}_{\theta})$ be the unital associative
   $\mathbb C$-algebra generated by $2n$ elements $\Gamma^\mu$,
   $\Gamma^{\nu\ast}$ ($\mu, \nu=1,\dots,n$) with relations
   \begin{eqnarray}
       \Gamma^\mu\Gamma^\nu+\lambda^{\nu\mu} \Gamma^\nu\Gamma^\mu & = &
       0 \label{d}\\
      \Gamma^{\mu\ast}\Gamma^{\nu\ast}+\lambda^{\nu\mu}
      \Gamma^{\nu\ast}\Gamma^{\mu \ast} & = & 0 \label{e}\\
      \Gamma^{\mu\ast}\Gamma^\nu+\lambda^{\mu\nu} \Gamma^\nu\Gamma^{\mu
      \ast} & = & \delta^{\mu\nu}\bbbone\label{f}
      \end{eqnarray}
where $\bbbone$ denotes the unit of the algebra.  For $\theta=0$ one
recovers the usual Clifford algebra of $\mathbb R^{2n}$; the familiar
generators $\gamma^a\ (a=1,2,\dots,2n)$ associated to the canonical
basis of $\mathbb R^{2n}$ being then given by
$\gamma^\mu=\Gamma^\mu+\Gamma^{\mu\ast}$ and
$\gamma^{\mu+n}=-i(\Gamma^\mu-\Gamma^{\mu\ast})$.  There is a unique
involution $\Lambda \mapsto \Lambda^\ast$ such that
$(\Gamma^\mu)^\ast=\Gamma^{\mu\ast}$ for which $\clif(\mathbb
R^{2n}_{\theta})$ is a unital complex $\ast$-algebra.  One also endows
$\clif(\mathbb R^{2n}_{\theta})$ with a $\mathbb Z_{2}$-grading of
algebra by giving odd degree to the $\Gamma^\mu,\Gamma^{\nu\ast}$.
  The relations (\ref{d}), (\ref{e}) and (\ref{f}) imply
that the hermitian element
$[\Gamma^{\mu\ast},\Gamma^\mu]=\Gamma^{\mu\ast}\Gamma^\mu -\Gamma^\mu
\Gamma^{\mu\ast}$ anticommutes with $\Gamma^\mu$ and
$\Gamma^{\mu\ast}$ whereas it commutes with $\Gamma^\nu$ and
$\Gamma^{\nu\ast}$ for $\nu\not= \mu$ and that furthermore one has
$([\Gamma^{\mu\ast},\Gamma^\mu])^2=\bbbone$.  It follows that
$\gamma\in\clif(\mathbb R^{2n}_{\theta})$ defined by

\begin{equation}
	\gamma=[\Gamma^{1\ast},\Gamma^1] \dots [\Gamma^{n\ast},
	\Gamma^n] = \prod^n_{\mu=1}
	[\Gamma^{\mu\ast},\Gamma^\mu]\label{g}
        \end{equation}
is hermitian $(\gamma=\gamma^\ast)$ and satisfies

\begin{equation}
      \gamma^2=1, \ \ \gamma\Gamma^\mu+ \Gamma^\mu\gamma=0, \>
      \gamma\Gamma^{\mu\ast} + \Gamma^{\mu\ast}\gamma=0 \label{h}
      \end{equation}
      in fact $\Lambda\mapsto \gamma\Lambda\gamma$ is the $\mathbb
      Z_2$-grading.  The very reason why we have imposed the relations
      (\ref{d}), (\ref{e}) and (\ref{f}) is the following easy lemma.

   \begin{lemma}\label{lem1}
   In the algebra $\clif(\mathbb R^{2n}_{\theta})\otimes \pol(\mathbb
   R^{2n}_{\theta})$, the elements $\Gamma^{\mu\ast}z^\mu$ and  $\Gamma^\rho
   \bar z^\rho$ $=\Gamma^\rho z^{\rho\ast}$, $\mu,\rho=1,\dots,n$,
   satisfy the following anticommutation relations
   $\Gamma^{\mu\ast}z^\mu \Gamma^{\rho\ast}$ $z^\rho + \Gamma^{\rho
   \ast}z^\rho \Gamma^{\mu\ast} z^\mu=0 $ $(\Gamma^\mu \bar z^\mu
   \Gamma^\rho \bar z ^\rho + \Gamma^\rho \bar z^\rho\Gamma^\mu \bar
   z^\mu =0)$ and $\Gamma^{\mu\ast}z^\mu\Gamma^\rho\bar z^\rho +
   \Gamma^\rho \bar z^\rho \Gamma^{\mu\ast} z^\mu = \delta^{\mu\rho}
   z^\mu \bar z^\mu$ which do not depend on $\theta$.
     \end{lemma}
     This straightforward result is a key to reduce lots of computations
     to the classical case $\theta=0$, (see below).  The
     next result shows that $\clif(\rnt)$ is isomorphic to the usual
     $\clif(\mathbb R^{2n})$ as $\ast$-algebra and as $\mathbb
     Z_2$-graded algebra.

     \begin{proposition}\label{prop1}
         The following equality gives a faithful $\ast$-representation $\pi$ of
         $\clif(\rnt)$ in the Hilbert space $\otimes^n\mathbb C^2$,
         \[
         \begin{array}{l}
        \!\!\!\pi(\Gamma^{\mu\ast})=\left(
         \begin{array}{cc}
	   -\lambda^{1\mu} & 0\\
	   0 & 1 \end{array}\right)
	   \otimes \dots \otimes
	   \left(\begin{array}{cc} -\lambda^{\mu-1\mu} & 0\\
	   0 & 1
	   \end{array}\right)
	   \otimes \left( \begin{array}{cc}
	   0 & 1\\
	   0 & 0\end{array}\right)
	   \otimes \bbbone_2\otimes \dots
	   \otimes \bbbone_2\\
	 \>\>\>\>\>\>\>\>\>\>\>\>=\pi(\Gamma^\mu)^\ast
	   \end{array}
	   \]
     and $\pi$ is the unique irreducible $\ast$-representation of
     $\clif(\rnt)$ up to a unitary equivalence.
     \end{proposition}
     The proof is straightforward.  Note that $\otimes^n\mathbb C^2$
     , viewed as the graded tensor product of $\mathbb C^2$ graded by
$\left (
     \begin{array}{cc} 1 & 0\\
     0 & -1\end{array}\right)$ is a $\mathbb Z_{2}$-graded $\clif(\rnt)$-module.
   One has $\pi(\gamma)=\otimes^n\left(\begin{array}{cc} 1
     & 0\\ 0 & -1 \end{array}\right)$.  In the following we will use the
     above representation to identify $\clif(\rnt)$ with
     $M_{2^n}(\mathbb C)$.

     \section{Spherical property of $\theta$-deformed spheres}\label{sec3}
   \setcounter{equation}{0}
   Let  $\pol(\rntu)$, {\sl the algebra of polynomial functions on the
   noncommutative $(2n+1)$-plane $\rntu$},  be the unital complex
   $\ast$-algebra obtained by adding an hermitian
   generator $x$ to $\pol(\rnt)$ with relations $xz^\mu=z^\mu x$
   ($\mu=1,\dots,n$), i.e. $\pol(\rntu)\simeq \pol(\rnt)\otimes \mathbb
   C [x]\simeq \pol(\rnt)\otimes \pol(\mathbb R)$.  One knows that the
   $z^\mu\bar z^\mu=\bar z^\mu z^\mu$ and $x$ are in the center so
   $\sum^n_{\mu=1} z^\mu\bar z^\mu+x^2$ is also in the center
    $\pol(\rntu)$.  We let  $\pol(\snt)$  be
   the $\ast$-algebra quotient of $\pol(\rntu)$ by the ideal generated by
   $\sum^n_{\mu=1} z^\mu \bar z^\mu+x^2-\bbbone$.  In the following, we
   shall denote by $u^\mu, \bar u^\nu=u^{\nu\ast}$, $u$ the canonical
   images of $z^\mu, \bar z^\nu, x$ in $\pol(\snt)$.  On the unital
   complex $\ast$-algebra $\pol(\snt)$ there is a greatest
   $C^\ast$-seminorm which is a norm; the $C^\ast$-algebra $C(\snt)$
   obtained by completion will be refered to as {\sl the algebra of
   continuous functions on the noncommutative $2n$-sphere $\snt$}.\\

\noindent It is worth noticing that the noncommutative $2n$-sphere $\snt$ can
    be viewed as ``one-point compactification" of the noncommutative
    $2n$-plane $\rnt$.  To explain this, let us slightly enlarge the
    $\ast$-algebra $\pol(\rnt)$ by adjoining a hermitian central
    generator $( 1+\sum^n_{\mu=1} \bar z^\mu z^\mu )^{-1}=
    (1+\vert z\vert^2)^{-1}$ with relation\linebreak[4] 
$(1+\sum^n_{\mu=1}\bar z^\mu
    z^\mu)(1+\vert z\vert^2)^{-1} = (1+\vert z\vert^2)^{-1}(
    1+\sum^n_{\mu=1}\bar z^\mu z^\mu)=1$.  As will become
    clear  $(1+\vert z\vert^2)^{-1}$ is  smooth so that in fact we are
    staying in the algebra $C^\infty(\rnt)$ of smooth functions on $\rnt$.
    By setting
    \[
    \tilde u^\mu = 2z^\mu(1+\vert z\vert^2)^{-1},\ \tilde
    u^{\nu\ast}=2\bar z^\mu(1+\vert z\vert^2)^{-1},\ \tilde
    u=(1-\sum^n_{\mu=1}\bar z^\mu z^\mu)(1+\vert z\vert^2)^{-1},
    \]
    one sees that the $\tilde u^\mu$, $\tilde u^{\nu\ast}, \tilde u$
    satisfy the same relations as the $u^\mu$, $u^{\nu\ast}, u$.  The
    ``only difference" is that the classical point $u^\mu=0, \bar
    u^\mu=0, u=-1$ of $\snt$ does not belong to the spectrum of $\tilde
    u^\mu, \tilde u^{\nu\ast},\tilde u$.  In the same spirit, one can
    cover $\snt$ by two ``charts" with domain $\rnt$ with transition on
    $\rnt\backslash \{0\}$, ($z^\mu=0, \bar z^\nu=0$ being a classical
    point of $\rnt$).\\

\noindent Let  $\calg(S^{2n-1}_{\theta})$  be
    the quotient of the $\ast$-algebra $\calg(\rnt)$ by the two-sided
ideal generated by
    the element $\sum^n_{\mu=1} z^\mu\bar z^\mu-\bbbone$ of the center
    of $\calg(\rnt)$. This defines by duality the noncommutative
    $(2n-1)$-sphere $S^{2n-1}_{\theta}$. In the following, we shall
    denote by $v^\mu$, $\bar v^\nu$ the canonical images of $z^\mu$,
    $\bar z^\nu$ in $\calg(S^{2n-1}_{\theta})$. Again there is a greatest
    $C^\ast$-seminorm which is a norm on $\calg(S^{2n-1}_{\theta})$;
    the $C^\ast$-algebra obtained by completion will be refered to as
    {\sl the algebra of continuous functions on the noncommutative
    $(2n-1)$-sphere $S^{2n-1}_{\theta}$}. It is clear that, in an
    obvious sense, $S^{2n}_{\theta}$ is the suspension of
    $S^{2n-1}_{\theta}$.\\

\noindent As for the case of $\rnt$, one has an action $\sigma$ of $T^n$ on
    $\rntu$, $\snt$ and $S^{2n-1}_{\theta}$ which is induced by an
    action on the corresponding classical spaces. More precisely the
    group-homomorphism $s\mapsto\sigma_{s}$ of $T^n$ into
    $\aut(\calg(\rnt))$
    extends as a group-homomorphism $s\mapsto\sigma_{s}$ of $T^n$
    into $\aut(\calg(\rntu))$ and these group-homomorphisms induce group
    homomorphisms $s\mapsto\sigma_{s}$ of $T^n$ into
    $\aut(\calg(S^{2n-1}_{\theta}))$
    and of $T^n$ into $\aut(\calg(\snt))$. As for $\rnt$, one
    checks that the subalgebras of $\sigma$-invariant elements are in
    the respective centers, are not deformed, and are
    isomorphic to the subalgebras of $\sigma$-invariant elements of
    $\calg(\mathbb R^{2n+1})$, $\calg(S^{2n})$  and $\calg(S^{2n-1})$
    respectively.\\

\noindent In order to formulate the last part of the next theorem, let us
    notice that, in view of (\ref{g}) and  (\ref{h}), there is an
    injective representation of $\clif(\rnt)$ for which
    $\gamma=\left( \begin{array}{cc}
    \bbbone & 0\\
    0 & -\bbbone
    \end{array}
    \right)$ where $\bbbone$ denotes the unit of $M_{2^{n-1}}(\mathbb C)$. In
    such a representation one has in view of (\ref{h})
    \[
    \Gamma^\mu=\left(
    \begin{array}{cc}
        0 & \sigma^\mu\\
        \bar \sigma^{\mu \ast} & 0
        \end{array}
        \right), \> \>
     \Gamma^{\mu\ast}=\left(
    \begin{array}{cc}
        0 & \bar\sigma^\mu\\
        \sigma^{\mu \ast} & 0
        \end{array}
        \right)\]
        where $\sigma^\mu$ and $\bar\sigma^\mu$ are in
        $M_{2^{n-1}}(\mathbb C)$.

   \begin{theo}\label{theo1}
       $(i)$ One obtains a hermitian projection $e\in M_{2^n}(\calg(\snt))$ by
       setting $e=\frac{1}{2}(\bbbone +
       \sum^n_{\mu=1}(\Gamma^{\mu\ast}u^\mu + \Gamma^\mu
       u^{\mu\ast})+\gamma u)$.  Furthermore one has $\ch_{m}(e)=0$ for
       $0\leq m\leq n-1$.\\
       \phantom {THEOREM 3 } $(ii)$ One obtains a unitary $U\in
       M_{2^{n-1}}(\calg(S^{2n-1}_{\theta}))$  by setting
       $U=\sum^n_{\mu=1}(\bar \sigma^\mu v^\mu+\sigma^\mu \bar v^\mu)$ ,
       where $\sigma^\mu$ and $\bar\sigma^\mu$ are as above. Furthermore
       one has $\ch_{m-\frac{1}{2}}(U)=0 $ for $1\leq m\leq n-1$.
       \end{theo}

       \noindent The relation $e=e^\ast$ is obvious.
       It follows from Lemma~\ref{lem1} that
       \[
     \left(\sum^n_{\mu=1}(\Gamma^{\mu\ast} z^\mu+\Gamma^\mu
       z^{\mu\ast})\right )^2=\sum^n_{\mu=1} z^\mu \bar z^\mu,
       \]
      which in terms of the $\sigma^\mu$ reads
      \[
      (\bar \sigma^\mu z^\mu+\sigma^\mu\bar z^\mu)  (\bar \sigma^\mu
      z^\mu+\sigma^\mu\bar z^\mu)^\ast=  (\bar \sigma^\mu
      z^\mu+\sigma^\mu\bar z^\mu)^\ast  (\bar \sigma^\mu
      z^\mu+\sigma^\mu\bar z^\mu)= \sum^n_{\mu=1}z^\mu\bar z^\mu.
      \]

\noindent On the
       other hand relations (\ref{h}) imply then
       \[
       \left(\sum^n_{\mu=1}
       (\Gamma^{\mu\ast} z^\mu+\Gamma^\mu z^{\mu\ast} )+\gamma x\right
       )^2= \sum^n_{\mu=1} z^\mu \bar z^\mu +x^2
       \]
       which reduces to
       $\bbbone\in M_{2^n} (\calg(\snt))$.  This is equivalent to $e^2=e$.
       Using again Lemma~\ref{lem1}, $\ch_{m}(e)=0$ for $m<n$ follows from the
       vanishing of the corresponding traces of products of the
       $\Gamma^\mu, \Gamma^{\mu\ast},\gamma$ in the representation of
       Proposition~\ref{prop1}.  The unitarity of $U\in
       M_{2^{n-1}}(\calg(S^{2n-1}_{\theta}))$ is clear whereas one has
       \begin{equation}
        \ch_{m-\frac{1}{2}}(U)=\tr\left( (U\circledcirc
        U^\ast)^{\circledcirc m}-(U^\ast\circledcirc U)^{\circledcirc
        m}\right)
        \label{ch}
        \end{equation}
        which implies
        \begin{equation}
         \ch_{m-\frac{1}{2}}(U)=\tr\left(\frac{1+\gamma}{2}
         \Gamma^{\circledcirc 2m}-\frac{1-\gamma}{2}
         \Gamma^{\circledcirc 2m}\right) =\tr(\gamma\Gamma^{\circledcirc 2m})
         \label{ch'}
         \end{equation}
         where $\Gamma=\sum_{\mu}(\Gamma^{\mu\ast} v^\mu+\Gamma^\mu
         \bar v^\mu)\in M_{2^n}(\calg(S^{2n-1}_{\theta} ))$ and where
         in (\ref{ch'}) $\tr$ and $\circledcirc$ are taken for
         $M_{2^n}$ instead of $M_{2^{n-1}}$ as in (\ref{ch}), (see the
         definitions at the end of the introduction). It follows from
         (\ref{ch'}) that one has $\ch_{m-\frac{1}{2}}(U)=0$ for $1\leq
         m \leq n-1$ for the same reasons as $\ch_{m}(e)=0$ for $m\leq
         n-1$.\\

\noindent This theorem combined with the last theorem of Section
         \ref{sec10} and the last theorem of Section
         \ref{sec11} implies that $S^m_{\theta}$   is an
         $m$-dimensional noncommutative spherical manifold.\\

\noindent It follows from  $\ch_{m}(e)=0$ for $0\leq m\leq n-1$ that
    $\ch_{n}(e)$ is a Hochschild cycle which corresponds to the volume
    form on $\snt$. In fact it is obvious that the whole analysis of
    Section III and IV of \cite{connes:08} generalizes from
    $S^4_{\theta}$ to $\snt$. This is in particular the case of Theorem
    3 of \cite{connes:08}  (with the appropriate changes e.g. $4
    \mapsto 2n$ and $M_{4}(\mathbb C)\mapsto M_{2^n}(\mathbb C)$). The
    odd case is obviously similar. This will be discussed in more
    details in Section~\ref{sec11}.\\

\noindent The projection $e$  is a
    noncommutative version of the projection-valued field $P_{+}$ on the
    sphere $S^{2n}$ described in Section 2.7 of \cite{mdv2}; one has
    $P_{+}=e\vert_{\theta=0}$. As was shown there, $P_{+}$ satisfies
    the following self-duality equation
    \begin{equation}
        \ast P_{+}(dP_{+})^{n}=i^n P_{+}(dP_{+})^{n}
        \label{sd}
        \end{equation}
where $\ast$  is the usual Hodge duality of forms on $S^{2n}$.
Since $\ast$ is conformally invariant on forms of degree $n$, this
equation is conformally invariant. The above equation
generalizes to $e$ (i.e. on $\snt$) once the appropriate differential
calculus and metric are defined,  (see Theorem 6 of section 12 below).
  For $n$ even, Equation
   (\ref{sd}) describes an intanton  (the ``round" one) for a
   conformally invariant generalization of the classical Yang-Mills
   action
   on $S^{2n}$ (which reduces to the Yang-Mills action on $S^4$),
   \cite{mdv2}. The fact, which was pointed out and used
   in \cite{mdvg}, that classical gauge theory can be formulated in
   terms of projection-valued fields is a direct consequence of the
   theorem of Narasimhan and Ramanan on the existence of universal
   connections \cite{nara}, \cite{nara2}, (see also in \cite{mdv1} for a
   short economical proof of this theorem).\\
   It is clear that by changing $(u^\mu,u)$ into $(-u^\mu,-u)$ one
    also obtains a hermitian projection $e_{-}\in M_{2^n}(\calg(\snt))$
   satisfying $\mbox{ch}_{m}(e_{-})=0$ for $0\leq m\leq n-1$. For
   $\theta=0$, $e_{-}$ coincides with the projection-valued field $P_{-}$
   on $S^{2n}$ of \cite{mdv2} which satisfies  $\ast
   P_{-}(dP_{-})^{n}=-i^n P_{-}(dP_{-})^{n}$. What
   replaces $e\mapsto e_{-}$ for the odd-dimensional case is $U\mapsto
   U^\ast$.

   \section{The graded differential algebras
   $\Omega_{\mathrm{alg}}(\mathbb R^{m}_{\theta})$ and
   $\Omega_{\mathrm{alg}}(S^{m}_{\theta})$} \label{sec4}
   \setcounter{equation}{0}
   There are canonical differential calculi,
   $\Omega_{\mathrm{alg}}(\rnt)$ and $\Omega_{\mathrm{alg}}(\rntu)$, on
   the noncommutative planes $\rnt$ and $\rntu$, which are
deformations of the differential
   algebras of polynomial differential forms on $\mathbb R^{2n}$ and
   $\mathbb R^{2n+1}$ and which are such that the $z^\mu\bar z^\mu =
   \bar z^\mu z^\mu$ are in the center of $\Omega_{\mathrm{alg}}(\rnt)$
   and $\Omega_{\mathrm{alg}}(\rntu)$ as well as $x$ in the case
   $\Omega_{\mathrm{alg}}(\rntu)$.  Let us first give a detailed
   description of the graded differential algebra
   $\Omega_{\mathrm{alg}}(\rnt)$.\\

\noindent As a complex unital associative graded algebra
   $\Omega_{\mathrm{alg}}(\rnt) = \oplus_{p\in\mathbb
   N}\Omega_{\mathrm{alg}}^p(\rnt)$ is generated by $2n$ elements
   $z^\mu,\bar z^\nu$ of degree 0 with relations (\ref{a}) and by $2n$
elements $dz^\mu,d\bar z^\nu$ of degree 1 with
   relations
   \begin{equation}
       dz^\mu dz^\nu  + \lambda^{\mu\nu} dz^\nu dz^\mu = 0,  \>
       d\bar z^\mu d\bar z^\nu  + \lambda^{\mu\nu} d\bar z^\nu d\bar
       z^\mu = 0,\>
       d\bar z^\mu dz^\nu  +  \lambda^{\nu\mu} dz^\nu d\bar z^\mu =
       0  \label{i}
        \end{equation}
   \begin{equation}
       z^\mu dz^\nu  = \lambda^{\mu\nu} dz^\nu z^\mu, \>
       \bar z^\mu d\bar z^\nu =  \lambda^{\mu\nu} d\bar z^\nu \bar
       z^\mu, \>
       \bar z^\mu dz^\nu  =  \lambda^{\nu\mu} dz^\nu \bar
       z^\mu, \>
       z^\mu d\bar z^\nu =  \lambda^{\nu\mu} d\bar z^\nu
       z^\mu\label{j}
       \end{equation}
for any $\mu,\nu\in\{1,\dots,n\}$.  There is a unique
       differential $d$ of $\Omega_{\mathrm{alg}}(\rnt)$, (i.e. a unique
       antiderivation $d$ satisfying $d^2=0$), which extends the mapping
       $z^\mu\mapsto dz^\mu$, $\bar z^\nu\mapsto d\bar z^\nu$.  One
       extends $z^\mu\mapsto \bar z^\mu$, $dz^\nu\mapsto d\bar
       z^\nu=\overline{(dz^\nu)}$ as an antilinear involution
       $\omega\mapsto\bar \omega$ of $\Omega_{\mathrm{alg}}(\rnt)$ such
       that $\overline{\omega\omega'}=(-1)^{pp'}\bar\omega'\bar \omega$
       for $\omega\in \Omega_{\mathrm{alg}}^p(\rnt)$ and $\omega'\in
       \Omega_{\mathrm{alg}}^{p'}(\rnt)$.  One has $d\bar
       \omega=\overline{d\omega}$, $\forall
       \omega\in\Omega_{\mathrm{alg}}(\rnt)$.  Elements
       $\omega\in\Omega_{\mathrm{alg}}(\rnt)$ satisfying $\omega=\bar
       \omega$ will be refered to as {\sl real elements}.  Notice that
       the $\bar z^\mu z^\mu$, $\bar z^\mu dz^\mu$, $ z^\mu d\bar
       z^\mu$, $d\bar z^\mu dz^\mu$ for $\mu\in\{1,\dots,n\}$ generate
  a graded differential subalgebra of the
       graded center of $\Omega_{\mathrm{alg}}(\rnt)$ which coincides
       with this graded center whenever $\theta$ is generic. Notice also
       that these elements are invariant by the canonical extension
       to $\Oalg(\rnt)$ of the action $\sigma$ of $T^n$ on
       $\calg(\rnt)=\Omega_{\mathrm{alg}}^0(\rnt)$ (see the end of
this section).

\noindent There is another useful way to construct
       $\Omega_{\mathrm{alg}}(\rnt)$ which we now describe.  Consider
       the graded algebra $\pol(\rnt)\otimes_{\mathbb R}\wedge \mathbb
       R^{2n}= \pol(\rnt)\otimes\wedge_{c} \mathbb R^{2n}$ where
       $\wedge_{c}\mathbb R^{2n}$ is the complexified
       exterior algebra of $ \mathbb R^{2n}$.  The graded algebra
       $\pol(\rnt)\otimes\wedge_{c} \mathbb R^{2n}$ is the unital
       complex graded algebra generated by $2n$ elements of degree zero,
       $z^\mu, \bar z^\nu$ ($\mu, \nu= 1,\dots, n$) satisfying
       relations (\ref{a}) and by $2n$ elements of
       degree one, $\xi^\mu$, $\bar\xi^\nu$ ($\mu,\nu=1,\dots, n$) with
       relations
       \begin{eqnarray}
	 \xi^\mu \xi^\nu + \xi^\nu \xi^\mu = 0, \bar\xi^\mu
	 \bar\xi^\nu + \bar\xi^\nu \bar\xi^\mu = 0, \bar\xi^\mu
	 \xi^\nu + \xi^\nu \bar\xi^\mu = 0 \label{p}\\
	 z^\mu \xi^\nu = \xi^\nu z^\mu,\ \bar z^\mu \xi^\nu = \xi^\nu
	 \bar z^\mu, z^\mu \bar\xi^\nu = \bar\xi^\nu z^\mu,\ \bar
	 z^\mu \bar\xi^\nu = \bar\xi^\nu \bar z^\mu\label{q}
	 \end{eqnarray}
for $\mu, \nu\in \{1,\dots,n\}$.  The $2n$ elements $\xi^\mu,\bar
\xi^\nu$ satisfying (\ref{p}) generate the complexified exterior
algebra $\wedge_{c}\mathbb R^{2n}$.  An involution
$\omega\mapsto\bar \omega$ of graded algebra on $\pol(\rnt)\otimes
\wedge_{c}\mathbb R^{2n}$ is obtained by setting $\overline{z^\mu}=\bar
z^\mu,\overline{\bar z^\mu}=z^\mu$ as before and by setting
$\overline{\xi^\mu}=\bar\xi^\mu$, $\overline{\bar\xi^\mu}=\xi^\mu$.
There is a unique differential $d$ on the graded differential algebra
$\pol(\rnt)\otimes \wedge _{c}\mathbb R^{2n}$ such that
\begin{eqnarray}
      d\xi^\mu  =  0 &,& \>
      d\bar\xi^\mu  =  0\label{r}\\
      dz^\mu  =  z^\mu \xi^\mu &,& \>
      d\bar z^\mu  =   \bar z^\mu \bar\xi^\mu \label{s}
      \end{eqnarray}
      for $\mu=1,\dots,n$.  One then has
      $d\bar\omega=\overline{d\omega}$ for any $\omega\in
      \pol(\rnt)\otimes\wedge_{c} \mathbb R^{2n}$.  It is readily
      verified that the $dz^\mu,d\bar z^\nu$ defined by (\ref{s})
satisfy relations
      (\ref{i}) to (\ref{j}).  In other words
      $\Omega_{\mathrm{alg}}(\rnt)$ is the differential subalgebra of
      $\pol(\rnt)\otimes\wedge_{c}\mathbb R^{2n}$ generated by the
      $z^\mu,\bar z^\nu$ ($\mu,\nu=1,\dots, n$).  Furthermore the
      involution $\omega\mapsto \bar\omega$ of
      $\pol(\rnt)\otimes\wedge_{c} \mathbb R^{2n}$ induces on
      $\Omega_{\mathrm{alg}}(\rnt)$ the previously defined involution.
        As $\pol(\rnt)$-bimodule, one has
      $\Omega_{\mathrm{alg}}^p(\rnt)\subset
      \pol(\rnt)\otimes\wedge^p_{c}\mathbb R^{2n}$ so that
      $\Omega_{\mathrm{alg}}^p(\rnt)$ is a sub-bimodule of the
diagonal bimodule $(\pol(\rnt))^{C^p_{2n}}$,
      thus the $\Omega_{\mathrm{alg}}^p(\rnt)$ are {\sl diagonal
      bimodules} over $\pol(\rnt)$ \cite{mdv:pm1}.  This implies in
particular that
      $\Omega_{\mathrm{alg}}(\rnt)$ is a quotient of the graded
      differential algebra $\Omega_{\diagth}(\pol(\rnt))$ \cite{mdv3}.\\

\noindent The differential algebra $\pol(\rnt)\otimes\wedge_{c}\mathbb
     R^{2n}$ has the following interpretation.  Let us ``suppress" the
     classical points $z^\mu=0$ ($\mu=1,\dots,n$) of $\rnt$ by adjoining
     $n$ real (hermitian) central generators of degree zero $\vert
     z^\mu\vert^{-2}$ to $\Omega_{\mathrm{alg}}(\rnt)$ with relations
     \[
     \bar z^\mu z^\mu\vert z^\mu \vert^{-2}= \vert z^\mu\vert^{-2} \bar
     z^\mu z^\mu=\bbbone
     \]
     for $\mu=1,\dots,n$.  This becomes a graded differential algebra
     $\tilde\Omega_{\mathrm{alg}}(\rnt)$ if one sets $d\vert
     z^\mu\vert^{-2}=-(\vert z^\mu\vert^{-2})^2 d(\bar z^\mu z^\mu)$ for
     $\mu=1,\dots,n$.

\noindent Then the algebra $\pol(\rnt)\otimes
\wedge_{c}\mathbb R^{2n}$
     is the subalgebra generated by the $z^\mu,\bar z^\nu$ and the
     $\xi^\mu=\vert z^\mu\vert^{-2} \bar z^\mu d z^\mu$,
     $\bar\xi^\nu=\vert z^\nu\vert^{-2} z^\nu d \bar z^\nu$ and it is a
     graded differential subalgebra of
     $\tilde\Omega_{\mathrm{alg}}(\rnt)$.  The algebra
     $\tilde\Omega_{\mathrm{alg}}(\rnt)$ is the $\theta$-deformation of
     the algebra of complex polynomial differential forms on $(\mathbb
     C\backslash \{0\})^n \subset \mathbb R^{2n}$.

\noindent The complex unital associative graded algebra
     $\Omega_{\mathrm{alg}}(\rntu)$ is defined as the graded tensor
product
$\Omega_{\mathrm{alg}}(\rnt)\otimes_{\gr}\Omega_{\mathrm{alg}}(\mathbb
R)$.
More concretely one adjoins to
     $\Omega_{\mathrm{alg}}(\rnt)$ one generator $x$ of degree zero and
     one generator $dx$ of degree one with relations
     \begin{equation}
         xdx  =  dxx, \>
         x\omega  =  \omega x,\>
         dx \omega  = (-1)^p \omega dx \label{x}
         \end{equation}
         for $\omega\in \Omega_{\mathrm{alg}}^p(\rnt)$.  One extends the
         differential $d$ of $\Omega_{\mathrm{alg}}(\rnt)$ as the unique
         differential $d$ of $\Omega_{\mathrm{alg}}(\rntu)$ mapping $x$
         on $dx$.  The graded involution of
         $\Omega_{\mathrm{alg}}(\rnt)$ is extended into a graded
         involution $\omega\mapsto \bar \omega$ of
         $\Omega_{\mathrm{alg}}(\rntu)$ by setting $\bar x=x$ and
         $\overline{dx}=dx$.  One has again $d\bar
         \omega=\overline{d\omega}$ for
         $\omega\in\Omega_{\mathrm{alg}}(\rntu)$.\\

\noindent Again $\Omega_{\mathrm{alg}}(\rntu)$ is the differential
         subalgebra of $\pol(\rntu)\otimes\wedge_{c} \mathbb R^{2n+1}$
         generated by the $z^\mu, \bar z^\nu,x$ where the $(2n+1)$-th
         basis element of $\mathbb R^{2n+1}$ is identified with $dx$
         i.e. $\pol(\rntu)\otimes\wedge_{c} \mathbb R^{2n+1} \simeq
         (\pol(\rnt)\otimes \wedge_{c}\mathbb R^{2n})\otimes \wedge
         (x,dx)$.  Thus again the $\Omega_{\mathrm{alg}}^p(\rntu)$ are
         diagonal bimodules over $\pol(\rntu)$ which implies that
         $\Omega_{\mathrm{alg}}(\rntu)$ is a quotient of
         $\Omega_{\diagth}(\pol(\rntu))$.  Notice that these
         identifications are compatible with the involutions of the
         corresponding graded differential algebras.\\

\noindent Let now $\Oalg(S^{2n-1}_{\theta})$ be the graded differential algebra
  quotient of $\Oalg(\rnt)$ by the differential two-sided ideal
generated by $\sum^n_{\mu=1}z^\mu \bar z^\mu-\bbbone$ and similarly
  $\Oalg(\snt)$  be the quotient of
$\Oalg(\rntu)$ by the differential two-sided ideal generated by
$\sum^n_{\mu=1} z^\mu \bar z^\mu + x^2-\bbbone$. These are again
graded-involutive algebras with real differentials. Furthermore, it
will be shown  using the splitting homomorphism that they are
diagonal bimodules over $\calg (S^{2n-1}_{\theta})$  and over
$\calg(\snt)$ repectively from which it follows that they are quotient
of $\Omega_{\diagth}(\calg(S^{2n-1}_\theta))$ and of
$\Omega_{\diagth}(\calg(\snt))$ respectively.\\

\noindent Let $m=2n$ or $2n+1$. The actions $s\mapsto \sigma_{s}$ of $T^n$ on
  $\calg(\mathbb R^{m}_{\theta})$ and
$\calg(S^{m-1}_{\theta})$  extend
canonically to actions  of $T^n$ as automorphisms of graded-involutive
differential algebras, $s\mapsto
\sigma_{s} \in \aut(\Oalg(\mathbb R^{m}_{\theta}))$,
and $s\mapsto
\sigma_{s} \in \aut(\Oalg (S^{m-1}_{\theta}))$.
   The differential
subalgebras $\Oalg(\mathbb R^{m}_{\theta})^\sigma$ and
$\Oalg(S^{m-1}_{\theta})^\sigma$
of
$\sigma$-invariant elements are in the graded centers of
$\Oalg(\mathbb R^{m}_{\theta})$ and $\Oalg(S^{m-1}_{\theta})$
  and they are undeformed, i.e.
isomorphic to the corresponding subalgebras $\Oalg(\mathbb
R^{m})^\sigma$ and
$\Oalg(S^{m-1})^\sigma$of $\Oalg(\mathbb
R^{m})$ and
$\Oalg(S^{m-1})$.

   \section{The quantum groups $GL_{\theta}(m,\mathbb
   R)$, $SL_{\theta}(m,\mathbb R)$  and $GL_{\theta}(n,\mathbb C)$}
  \label{sec5}
   \setcounter{equation}{0}
   In this section we shall give a concrete explicit description of 
the various quantum
groups of symmetries of the noncommutative spaces
   $\mathbb R^m_{\theta}$ and $\mathbb C^n_{\theta}$ for $m\geq 4$ and
   $n\geq 2$. There are
other approaches to quantum groups of symmetries of
$S^4_{\theta}$ and $\mathbb R^4_{\theta}$ and some generalizations \cite{sit},
\cite{var}, \cite{asch}. In
\cite{sit} the dual point of view is adopted and what is produced is
the deformation of the universal enveloping algebra whereas in
\cite{var} the deformation is on the same side of the duality as
developed here; both points of view are of course useful. However it
must be stressed that, beside the fact that our approach is
closely related to the differential calculus, the important point here
is the observation that the quantum groups we introduce arise
with their expected Hochschild dimensions which equals the dimensions of the
corresponding classical groups.
They are
   deformations  (called  {\sl $\theta$-deformations}) of the  classical groups
   $GL(m,\mathbb R)$, $SL(m,\mathbb R)$,  $GL(n,\mathbb C)$ and as will
   be shown in Section 12, the Hochschild dimension is an invariant of
   these deformations.
    It is worth noticing here  that
   there is no corresponding $\theta$-deformation of $SL(n,\mathbb 
C)$; the reason
   being that $dz^1\cdots dz^n$ is not central and not
   $\sigma$-invariant in $\Oalg(\mathbb C^n_{\theta})=\Oalg(\rnt)$.
   However, there is a $\theta$-deformation of the subgroup of
   $GL(n,\mathbb C)$ consisting of matrices with determinants of modulus one
    because $dz^1\cdots dz^n d\bar z^1\cdots d\bar z^n$ is
   $\sigma$-invariant and (consequently) central.\\

\noindent Let  $M_{\theta}(2n,\mathbb R)$  be the unital
         associative $\mathbb C$-algebra generated by $4n^2$ element
         $a^\mu_{\nu}, b^\mu_{\nu},\bar a^\mu_{\nu},\bar b^\mu_{\nu}$
         ($\mu, \nu=1,\dots,n$) with relations such that the elements
         $y^\mu$, $\bar y^\mu$, $\zeta^\mu$, $\bar\zeta^\mu$ of
         $M_{\theta}(2n,\mathbb R)\otimes \Omega_{\mathrm{alg}}(\rnt)$
         defined by
         \[
         y^\mu=a^\mu_{\nu}\otimes z^\nu +b^\mu_{\nu}\otimes \bar z^\nu,
         \> \> \bar{y}^\mu=\bar a^\mu_{\nu}\otimes \bar z^\nu+\bar
         b^\mu_{\nu} \otimes z^\nu,
         \]
         \[
         \zeta^\mu=a^\mu_{\nu}\otimes dz^\nu+b^\mu_{\nu}\otimes d\bar z,
         \> \> \bar{\zeta}^\mu = \bar a^\mu_{\nu}\otimes d\bar z^\nu
         +\bar b^\mu_{\nu} \otimes dz^\nu
         \]
         satisfy the relation
         \[
         y^\mu y^\nu=\lambda^{\mu\nu} y^\nu y^\mu, \>
         \bar{y}^\mu\bar{y}^\nu=\lambda^{\mu\nu} \bar{y}^\nu
         \bar{y}^\mu, \> \bar{y}^\mu y^\nu = \lambda^{\nu\mu}y^\nu
         \bar{y}^\mu,
         \]
         \[
         \zeta^\mu \zeta^\nu + \lambda^{\mu\nu} \zeta
         ^\nu\zeta^\mu=0, \> \bar \zeta^\mu\bar
         \zeta^\nu+\lambda^{\mu\nu}\bar{\zeta}^\nu \bar \zeta^\mu=0,\>
         \bar\zeta^\mu\zeta^\nu+\lambda^{\nu\mu}\zeta^\nu\bar{\zeta}^\mu=0.
         \]
         There is a unique $\ast$-algebra involution $a\mapsto a^\ast$
         on $M_{\theta}(2n,\mathbb R)$ such that
         $(a^\mu_{\nu})^\ast=\bar a^\mu_{\nu}$, $(b^\mu_{\nu})^\ast=\bar
         b^\mu_{\nu}$.  The relations between the generators are easy to
         write explicitely, they read
   \begin{eqnarray}
	   a^\mu_\nu a^\tau_{\rho}  =
	   \lambda^{\mu\tau}\lambda_{\rho\nu} a^\tau_{\rho}
	   a^\mu_{\nu} &,&
	   a^\mu_\nu \bar a^\tau_{\rho}  =
	   \lambda^{\tau\mu}\lambda_{\nu\rho} \bar a^\tau_{\rho}
	   a^\mu_{\nu}\label{y}\\
	   a^\mu_\nu b^\tau_{\rho}  =
	   \lambda^{\mu\tau}\lambda_{\nu\rho} b^\tau_{\rho}
	   a^\mu_{\nu}& , &
	   a^\mu_\nu \bar b^\tau_{\rho}  =
	   \lambda^{\tau\mu}\lambda_{\rho\nu} \bar b^\tau_{\rho}
	   a^\mu_{\nu}\label{z}\\
	   b^\mu_\nu b^\tau_{\rho}  =
	   \lambda^{\mu\tau}\lambda_{\rho\nu} b^\tau_{\rho}
	   b^\mu_{\nu} &,&
	   b^\mu_\nu \bar b^\tau_{\rho} =
	   \lambda^{\tau\mu}\lambda_{\nu\rho} \bar b^\tau_{\rho}
	   b^\mu_{\nu}\label{aa}
         \end{eqnarray}
       plus the relations obtained by hermitian conjugation, where we
       have also used the notation $\lambda_{\nu\rho}$ for
       $\lambda^{\nu\rho}$ to indicate that there is no summation in the
       above formulas.  This $\ast$-algebra becomes a $\ast$-bialgebra
       with coproduct $\Delta$ and counit $\varepsilon$ if we endow it
       with the unique algebra-homomorphism
         \[
         \Delta:M_{\theta}(2n,\mathbb R)\rightarrow
         M_{\theta}(2n,\mathbb R)\otimes M_{\theta}(2n,\mathbb R)
         \]
         and the unique character $\varepsilon:M_{\theta}(2n,\mathbb
         R)\rightarrow \mathbb C$ such that
         \begin{eqnarray}
	   \Delta a^\mu_{\nu} & = & a^\mu_{\lambda} \otimes
	   a^\lambda_{\nu} + b^\mu_{\lambda} \otimes \bar
	   b^\lambda_{\nu} ,\> \> \varepsilon(a^\mu_{\nu}) =
	   \delta^\mu_{\nu}\label{ae}\\
	   \Delta \bar a^\mu_{\nu} & = & \bar a^\mu_{\lambda} \otimes
	   \bar a^\lambda_{\nu} + \bar b^\mu_{\lambda} \otimes
	   b^\lambda_{\nu} ,\> \> \varepsilon(\bar a^\mu_{\nu}) =
	   \delta^\mu_{\nu}\label{af}\\
	   \Delta b^\mu_{\nu} & = & a^\mu_{\lambda} \otimes
	   b^\lambda_{\nu} + b^\mu_{\lambda} \otimes \bar
	   a^\lambda_{\nu} ,\> \> \varepsilon(b^\mu_{\nu} ) =
	   0\label{ag}\\
	   \Delta \bar b^\mu_{\nu} & = & \bar a^\mu_{\lambda} \otimes
	   \bar b^\lambda_{\nu} + \bar b^\mu_{\lambda} \otimes a
	   ^\lambda_{\nu} ,\> \> \varepsilon(\bar b^\mu_{\nu} ) =
	   0\label{ah}
	   \end{eqnarray}
	   for any $\mu,\nu\in\{1,\dots,n\}$.  It is easy to verify
	   that there is a unique algebra-homomorphism $\delta
	   :\Omega_{\mathrm{alg}}(\rnt)\rightarrow
	   M_{\theta}(2n,\mathbb R)\otimes\Omega_{\mathrm{alg}}(\rnt)$
	   such that $\delta z^\mu=y^\mu$, $\delta\bar z^\mu=\bar
	   y^\mu$, $\delta dz^\mu=\zeta^\mu$, $\delta d\bar
	   z^\mu=\bar\zeta^\mu$ and that this is furthermore a
	   graded-involutive algebra-homomorphism.  In fact, this is
	   another way to obtain $\Omega_{\mathrm{alg}}(\rnt)$
	   starting from $\pol(\rnt)$ and from the $\theta$-twisted
	   complexified exterior algebra $\wedge_{c}\rnt$ generated by
	   the $dz^\mu, d\bar z^\nu$ satisfying (\ref{i}).
          One has
	   \begin{equation}
	       (\Delta \otimes I) \circ \delta = (I\otimes
	       \delta)\circ \delta,\> \> \>(\varepsilon\otimes
	       I)\circ\delta = I \label{ai}
	       \end{equation}
	       and $\delta\Omega_{\mathrm{alg}}^p(\rnt)\subset
	       M_{\theta}(2n,\mathbb R)\otimes
	       \Omega_{\mathrm{alg}}^p(\rnt),\> \> \forall p\in
	       \mathbb N$.

  \noindent One has of course $\delta\pol(\rnt)$ $\subset
	       M_{\theta}(2n,\mathbb R)\otimes\pol(\rnt)$, (this is
	       the previous result for $p=0$ since
	       $\pol(\rnt)=\Omega_{\mathrm{alg}}^0(\rnt)$), and 
$\delta\wedge_{c}\mathbb R
	       ^{2n}_{\theta}\subset M_{\theta}(2n,\mathbb R)\otimes
	       \wedge_{c}\mathbb R^{2n}_{\theta}$ with
	       $\delta\wedge^p_{c}\mathbb R^{2n}_{\theta}\subset
	       M_{\theta}(2n,\mathbb R)\otimes \wedge^p_{c}\mathbb
	       R^{2n}_{\theta}$ for any $p\in\mathbb N$.  Since
	       $\wedge^{2n}_{ c}\mathbb R^{2n}_{\theta}$ is of
	       dimension 1 and spanned by $d\bar z^1dz^1\dots d\bar
	       z^ndz^n=\prod^n_{\mu=1} d\bar z^\mu dz^\mu$, it follows
	       that one defines an element $\det_{\theta}\in
	       M_{\theta}(2n,\mathbb R)$ by setting
	       \begin{equation} \delta \prod^n_{\mu=1}d\bar z^\mu
	       dz^\mu=\detf_{\theta}\otimes \prod^n_{\mu=1}d\bar z^\mu
	       dz^\mu\label{aj} \end{equation} which satisfies \begin{eqnarray}
	       \Delta \detf_{\theta} & = &
	       \detf_{\theta}\otimes\detf_{\theta}\label{ak}\\
	       \varepsilon(\detf_{\theta}) & = & 1\label{al}
\end{eqnarray} and from
	       the fact that $\prod^n_{\mu=1}d\bar z^\mu dz^\mu$ is central in
	       $\Omega_{\mathrm{alg}}(\rnt)$ and from the very definition of
	       $M_{\theta}(2n,\mathbb R)$ it also follows that
$\det_{\theta}$ {\sl
	       belongs to the center of} $M_{\theta}(2n,\mathbb R)$.
The element
	       $\det_{\theta}$ of $M_{\theta}(2n,\mathbb R)$ is
clearly hermitian,
	       $(\det_{\theta})^\ast=\det_{\theta}$.\\

	       \noindent \underbar{Remark.} It is worth noticing that
Relations (\ref{y}),  (\ref{z}), (\ref{aa}) and their hermitian
conjugate are the quadratic
relations associated with a $R$-matrix $\hat R$ satisfying the  braid equation
(Yang-Baxter) and which is of square equal to 1, (i.e. $\hat R$ represents an
elementary transposition). In other words, the bialgebra
$M_{\theta}(2n,\mathbb R)$ is the
	      bialgebra of the $R$-matrix $\hat R$. \\

\noindent Let
	       $\pol(GL_{\theta}(2n,\mathbb R))$ be the
$\ast$-bialgebra obtained by
	       adding to $M_{\theta}(2n,\mathbb R)$ a hermitian central element
	       $\det^{-1}_{\theta}$ with relation
	       $\det_{\theta}.\det^{-1}_{\theta}=\bbbone
	       =\det^{-1}_{\theta}.\det_{\theta}$ and by setting
	       $\Delta\det^{-1}_{\theta}=\det^{-1}_{\theta} \otimes
	       \det^{-1}_{\theta}$ and
$\varepsilon(\det^{-1}_{\theta})=1$.  It is
	       not hard (but cumbersome) to see that the introduction of
	       $\det^{-1}_{\theta}$ allows to invert the ($2n,2n$) matrix
	       $L= \left (
	       \begin{array}{cc} A & B\\ \bar B & \bar A\end{array}\right )$ in
	       $M_{2n}(\pol(GL_{\theta}(2n,\mathbb R))$ and to obtain
an antipode $S$ on $\pol(GL_{\theta}(2n,\mathbb
	       R))$ which of course satisfies
$S(\det_{\theta})=\det^{-1}_{\theta}$
	       and $S(\det^{-1}_{\theta})=\det_{\theta}$.  Thus
	       $\pol(GL_{\theta}(2n,\mathbb R))$ {\sl is a
$\ast$-Hopf algebra} and
	       the quantum group $GL_{\theta}(2n,\mathbb R)$ is
defined to be the
	       dual object.\\ The quotient $\pol(SL_{\theta}(2n,\mathbb R))$ of
	       $M_{\theta}(2n,\mathbb R)$ by the relation
$\det_{\theta}=\bbbone$ is
	       also the quotient of $\pol(GL_{\theta}(2n,\mathbb R)$
by the two-sided
	       ideal generated by $\det_{\theta}-\bbbone$ and
	       $\det_{\theta}^{-1}-\bbbone$ which is a $\ast$-Hopf ideal.  So
	       $\pol(SL_{\theta}(2n,\mathbb R))$ is again a
$\ast$-Hopf algebra which
	       defines the quantum group $SL_{\theta}(2n,\mathbb R)$
by duality.\\
	       Replacing $\Omega_{\mathrm{alg}}(\rnt)$ by
	       $\Omega_{\mathrm{alg}}(\rntu)$ one defines in a similar way the
	       $\ast$-bialgebra $M_{\theta}(2n$ $+1,\mathbb R)$, the 
$\ast$-Hopf
	       algebras $\pol(GL_{\theta}(2n+1,\mathbb R))$,
	       $\pol(SL_{\theta}(2n+1,\mathbb R))$ and therefore the
quantum groups
	       $GL_{\theta}(2n+1,\mathbb R)$ and $SL_{\theta}(2n+1,\mathbb
	       R)$. \\

\noindent Finally, we let $\calg(GL_{\theta}(n,\mathbb C))$ be
	       the quotient of $\calg(GL_{\theta}(2n,\mathbb R))$ by the
	       ideal generated by the $b^\mu_{\nu}$ and the $\bar
	       b^\mu_{\nu}$ which is a $\ast$-Hopf ideal. The coaction of the
	       corresponding Hopf algebra on $\Oalg(\mathbb C^n_{\theta})$ is
	       straightforwardly obtained.  This defines the quantum group
	       $GL_{\theta}(n,\mathbb C)$ and its action on $\mathbb
	       C^n_{\theta}$. The ideal generated by the image of
	       $\det_{\theta}-\bbbone$ in $\calg(GL_{\theta}(n,\mathbb C))$
	       is a $\ast$-Hopf ideal and the corresponding quotient Hopf
	       algebra defines by duality a quantum group which is a
	       deformation ($\theta$-deformation) of the subgroup of
	       $GL(n,\mathbb C)$ which consists of matrices with
	       determinants of modulus one.

	       \section{The quantum groups $O_{\theta}(m)$,
	       $SO_{\theta}(m)$ and $U_{\theta}(n)$}\label{sec6}
	       \setcounter{equation}{0}
	       Let $\pol(O_{\theta}(2n))$ be the quotient of
	       $M_{\theta}(2n,\mathbb R)$ by the two-sided ideal
generated by
	       \[ \sum^n_{\mu=1}(\bar a^\mu_{\alpha}
a^\mu_{\beta} +
	       b^\mu_{\alpha}\bar b^\mu_{\beta})
-\delta_{\alpha\beta}\bbbone,\>
	       \sum^n_{\mu=1}(\bar a^\mu_{\alpha} b^\mu_{\beta} +
b^\mu_{\alpha} \bar
	       a^\mu_{\beta}),\>  \sum^n_{\mu=1}(\bar b^\mu_{\alpha}
	       a^\mu_{\beta}+a^\mu_{\alpha} \bar b^\mu_{\beta})
\] for
	       $\alpha,\beta =1,\dots, n$.  This ideal is
$\ast$-invariant and is
	       also a coideal.  It follows that
$\pol(O_{\theta}(2n))$ is again a
	       $\ast$-bialgebra.  Furthermore, one can show that
	       $(\det_{\theta})^2-\bbbone$ is in the above ideal (see below) so
	       $\pol(O_\theta(2n))$ is a $\ast$-Hopf algebra which is
a quotient of
	       $\pol(GL_{\theta}(2n,\mathbb R))$.  One verifies that 
the homomorphism
	       $\delta:\Omega_{\mathrm{alg}}(\mathbb
R^{2n}_{\theta})\rightarrow
	       M_{\theta}(2n,\mathbb R)\otimes \Omega_{\mathrm{alg}}(\mathbb
	       R^{2n}_{\theta})$ yields a homomorphism
	       \[
	       \delta_{R} : \Omega_{\mathrm{alg}}(\mathbb
R^{2n}_{\theta})\rightarrow
	       \pol(O_{\theta}(2n))\otimes \Omega_{\mathrm{alg}}(\mathbb
	       R^{2n}_{\theta})
	       \]
	       of graded-involutive algebras.  This yields the
	       quantum group $O_{\theta}(2n)$ which is a deformation
of the group of
	       rotations in dimension $2n$ and its action on $\mathbb
	       R^{2n}_{\theta}$ (cf. \cite{asch}).  Indeed one has \[
\delta_{R}(\sum^n_{\mu=1}\bar
	       z^\mu z^\mu)= \bbbone \otimes (\sum^n_{\mu=1}\bar
z^\mu z^\mu) \] by
	       the very definition of $\pol(O_{\theta}(2n))$.  One can notice
	       here that $\pol(O_{\theta}(2n))$ is a quotient of the
Hopf algebra of
	       the quantum group of the non-degenerate bilinear form
$B$ on $\mathbb
	       C^{2n}$ with matrix $\left ( \begin{array}{cc} 0 & \bbbone_{n}\\
	       \bbbone_{n} & 0 \end{array}\right)$ defined in
\cite{mdvl}, the later bilinear form is
	       equivalent to the metric of $\mathbb R^{2n}$, (the
involution being
	       defined accordingly).  The coaction $\delta_{R}$ passes to the
	       quotient to give the coaction
	       \[
	       \delta_{R}:\Oalg(S^{2n-1}_{\theta})\rightarrow
	       \calg(O_{\theta}(2n))\otimes \Oalg(S^{2n-1}_{\theta})
	       \]
	       which is also a homomorphism of graded-involutive algebras.
	       By taking a further quotient by the relation
	       $\det_{\theta}=\bbbone$, one obtains the $\ast$-Hopf algebra
	       $\pol(SO_{\theta}(2n))$ defining the quantum group
$SO_{\theta}(2n)$.
	       Let $\rho:M_{\theta}(2n,\mathbb R)\rightarrow
\pol(O_{\theta}(2n))$ be
	       the canonical projection.  The algebra
$\pol(O_{\theta}(2n))$ is the
	       unital $\ast$-algebra generated by the $4n^2$ elements
	       $\rho(a^\mu_{\nu})$, $\rho(b^\mu_{\nu})$, $\rho(\bar
a^\mu_{\nu})$,
	       $\rho(\bar b^\mu_{\nu})$ with relations induced by (\ref{y}),
	       (\ref{z}), (\ref{aa}) and the
	       relations \[
	       \sum_{\mu} (\rho(\bar
  a^\mu_{\alpha})\rho(a^\mu_{\beta})+\rho(b^\mu_{\alpha}) \rho(\bar
	       b^\mu_{\beta}))  = \delta_{\alpha\beta}\bbbone,\>
	       \sum_{\mu} (\rho(\bar
  a^\mu_{\alpha})\rho(b^\mu_{\beta})+\rho(b^\mu_{\alpha}) \rho(\bar
	       a^\mu_{\beta})) =  0 \]
	       (for $\alpha,
	       \beta=1,\dots,n$), together with $\rho(\bar
	       a^\mu_{\nu})=\rho(a^\mu_{\nu})^\ast$ and $\rho(\bar
	       b^\mu_{\nu})=\rho(b^\mu_{\nu})^\ast$.  It follows
that, for any $C^\ast$-semi-norm $\nu$ on
$\pol(O_{\theta}(2n))$ one
	       has $\nu(a^\mu_{\nu})=\nu(\bar a^\mu_{\nu})\leq 1$ and
	       $\nu(b^\mu_{\nu})=\nu(\bar b^\mu_{\nu})\leq 1$ so that
there is a
	       greatest $C^\ast$-semi-norm on $\pol(O_{\theta}(2n))$
which is a norm
	       and the corresponding completion $C(O_{\theta}(2n))$ of
	       $\pol(O_{\theta}(2n))$ is a $C^\ast$-algebra.  This defines
	       $O_{\theta}(2n)$ as a {\sl compact matrix quantum group}
	       \cite{slw:0}.  The same
	       applies to $SO_{\theta}(2n)$ which is therefore also a
compact matrix
	       quantum group.\\

\noindent One proceeds similarily (with obvious modifications)
	       to obtain the quantum groups $O_{\theta}(2n+1)$ and
	       $SO_{\theta}(2n+1)$ which are again compact matrix
quantum groups.
	       One has also the coaction
	       \[ \delta_{R}:\Omega_{\mathrm{alg}}(\mathbb
	       R^{2n+1}_{\theta})\rightarrow \pol(O_{\theta}(2n+1))\otimes
	       \Omega_{\mathrm{alg}}(\mathbb R^{2n+1}_{\theta})
	       \]
	       which passes to the quotient to yield the coaction
	       \[ \delta_{R}:\Oalg(S^{2n}_{\theta})\rightarrow
	       \calg (O_{\theta}(2n+1))\otimes \Oalg (S^{2n}_{\theta})
	       \]
	       these coactions are homomorphisms of graded-involutive algebras.
	       This gives the action
	       of the quantum group $O_{\theta}(2n+1)$ on the noncommutative
	       $2n$-sphere $S^{2n}_{\theta}$.  One obtains similarily
the action of
	       $SO_{\theta}(2n)$ on $S^{2n-1}_{\theta}$ and of
	       $SO_{\theta}(2n+1)$ on $S^{2n}_{\theta}$.  \\

\noindent Finally one lets $\calg(U_{\theta}(n))$  be the quotient
	       of $\calg(O_{\theta}(2n))$ by the ideal generated by the
	       $\rho(b ^\mu_{\nu})$ and $\rho(\bar b^\mu_{\nu})$ which is
	       also a $\ast$-Hopf ideal. The coactions $\delta_{R}$ of
	       $\calg(O_{\theta}(2n))$ on $\Oalg(\rnt)=\Oalg(\mathbb
	       C^n_{\theta})$ and on $\calg(S^{2n-1}_\theta)$ pass
	       to quotient to give corresponding coactions of
	       $\calg(U_{\theta}(n))$ . \\
	       Again there is no corresponding $\theta$-deformation 
of $SU(n)$.\\

\noindent Let us denote by $z_{\mu},\bar z_{\nu}=z^\ast_{\nu}$ the
	       generators of $\calg(\mathbb R^{2n}_{-\theta}) =\calg(\mathbb
	       C^n_{-\theta})$ satisfying
	       $z_{\mu}z_{\nu}=\lambda_{\nu\mu}z_{\nu}z_{\mu}$ and
	      $\bar z_{\mu}z_{\nu}=\lambda_{\mu\nu}z_{\nu}\bar z_{\mu}$. One
	      verifies that one obtains a unique $\ast$-homomorphism
	      $\varphi$ of $M_{\theta}(2n,\mathbb R)$ into
	      $\calg(\rnt)\otimes \calg(\mathbb R^{2n}_{-\theta} )$ by
	      setting $\varphi(a^\mu_{\nu})=z^\mu\otimes z_{\nu}$ and
	      $\varphi(b^\mu_{\nu})=z^\mu \otimes z_{\nu}^\ast$. This
	      homomorphism is injective and its image is invariant by the
	      action $\sigma\otimes \sigma$ of $T^n\times T^n$ on
	      $\calg(\rnt)\otimes\calg(\mathbb R^{2n}_{-\theta})$. We shall
	      again denote by $\sigma\otimes \sigma$ the corresponding action
	      of $T^n\times T^n$ on $M_{\theta}(2n,\mathbb R)$, i.e. the
	      group-homomorphism of $T^n\times T^n$ into
	      $\aut(M_{\theta}(2n,\mathbb R))$, e.g. one writes
	      $\sigma_{s}\otimes \sigma _{t}(a^\mu_{\nu})= e^{2\pi i
	      (s_{\mu}+t_{\nu}))}a^\mu_{\nu}$,  $\sigma_{s}\otimes
\sigma _{t}(b^\mu_{\nu})= e^{2\pi i
	      (s_{\mu}-t_{\nu}))}b^\mu_{\nu}$, etc. . This induces a
	      group-homomorphism (also denoted by $\sigma\otimes \sigma$)
	      of $T^n\times T^n$ into the group of automorphisms of unital
	      $\ast$-algebras (not necessarily preserving the
coalgebra structure)
of the polynomial algebra $\calg$ on each of
	      the quantum groups defined in this section and in
	      Section~\ref{sec5}. In
	      each case, the subalgebra of $\sigma\otimes\sigma$-invariant
	      elements is in the center and is undeformed, that is
	      isomorphic to the corresponding subalgebra for $\theta=0$.

	\section{The graded differential algebras $\Omega_{\alg}(G_{\theta})$
	as \\ graded differential Hopf algebras}\label{sec7}
	\setcounter{equation}{0}
	The relations (\ref{y}) to  (\ref{aa}) define the $\ast$-algebra
	$M_{\theta}(2n,\mathbb R)$ as $\calg(\mathbb R^{2N}_{\Theta})$ with
	$N=2n^2$ and where $\Theta\in M_{N}(\mathbb R)$ is the appropriate
	antisymmetric matrix (which depends on $\theta\in M_{n}(\mathbb R)$).
	Let $\Oalg(\mathbb R^{2N}_{\Theta})$ be the corresponding
	graded-involutive differential algebra as in Section~\ref{sec4}.

	\begin{proposition}\label{prop2}
	 The coproduct $\Delta$ of $M_{\theta}(2n,\mathbb R)$ has a unique
	 extension as homomorphism of graded differential algebras, again
	 denoted by $\Delta$, of $\Oalg(\mathbb R^{2N}_{\Theta})$ into
	 $\Oalg(\mathbb R^{2N}_{\Theta})\otimes_{\gr}\Omega_{\alg}(\mathbb
	 R^{2N}_{\Theta})$. The counit $\varepsilon$ of
	 $M_{\theta}(2n,\mathbb R)$ has a unique extension as
	 algebra-homomorhism, again denoted by $\varepsilon$, of
	 $\Oalg(\mathbb R^{2N}_{\Theta})$ into $\mathbb C$ with
	 $\varepsilon\circ d=0$.  The coaction ${\delta:\calg(\rnt)\rightarrow
	 M_{\theta}(2n,\mathbb R)\otimes \calg(\rnt)}$ has a unique extension
	 as homomorphism of graded differential algebras, again denoted by
	 $\delta$,  of $\Oalg(\rnt)$ into $\Oalg(\mathbb
	 R^{2N}_{\Theta})\otimes_{\gr} \Oalg(\rnt)$.
	 The extended $\Delta$ is coassociative and
	 the extended $\varepsilon$ is a counit for it and one has
	 $(\Delta\otimes I)\circ \delta=(I\otimes\delta)\circ \delta$,
	 $(\varepsilon\otimes I)\circ\delta=I$. These extended
	 homomorphisms are real.
	 \end{proposition}
	 In this proposition, $N=2n^2$ and $\Theta$ are as explained
above and one
	 endows $\Oalg(\mathbb R^{2N}_{\Theta})\otimes_{\gr}\Omega(\mathbb
	 R^{2N}_{\Theta})$ of the involution $\omega'\otimes \omega'\mapsto
	 \overline{\omega'\otimes\omega''}=\bar \omega' \otimes \bar
	 \omega''$. So equipped $\Oalg(\mathbb
	 R^{2N}_{\Theta})\otimes_{\gr}\Oalg  (\mathbb
	 R^{2N}_{\Theta})$ is a graded-involutive differential algebra and
	 the reality of $\Delta$ means
	 $\overline{\Delta(\omega)}=\Delta(\bar\omega)$. The uniqueness in
	 the proposition is obvious and the only thing to verify is the
	 compatibility of the extension with the relations
                     $da^\mu_{\nu}da^\tau_{\rho}+\lambda^{\mu\tau}
                      \lambda_{\rho\nu}
                      da^\tau_{\rho}da^\mu_{\nu}=0,\dots,
	a^\mu_{\nu}da^\tau_{\rho}=\lambda^{\mu\tau}
                      \lambda_{\rho\nu}
                      da^\tau_{\rho}a^\mu_{\nu},\dots,$ etc. which is
                      easy.  One proceeds similarily for $\delta$.
	 In short, $\Oalg(\mathbb
	 R^{2N}_{\Theta})$ is a {\sl graded-involutive differential
	 bialgebra}  and $\Oalg(\rnt)$ is a {\sl graded-involutive
	 differential comodule} over $\Oalg(\mathbb R^{2N}_{\Theta})$. Notice
	 that to say that $\Delta$ is a homomorphism of graded differential
	 algebras means that $\Delta$ is a homomorphism of graded algebras
	 and that one has {\sl the graded co-Leibniz rule} $\Delta\circ
	 d=(d\otimes I + (-I)^\gr\otimes d)\circ \Delta$.\\

\noindent By a {\sl graded differential Hopf algebra} we mean a
	 graded differential bialgebra which  admits an antipode; the
	 antipode $S$ is then necessarily unique and satisfies $S\circ
	 d=d\circ S$. The notion of graded-involutive differential Hopf
	 algebra is clear. By adding $\det^{-1}_{\theta}$ to
	 $M_{\theta}(2n,\mathbb R)=\Omega^0_{\alg}(\mathbb R^{2N}_{\Theta})$
	 as in Section~\ref{sec5} to obtain the Hopf algebra
	 $\calg(GL_{\theta}(2n,\mathbb R))$  and by setting
	 \[
	 \begin{array}  {lll}
	 [{\mathrm{det}}_{\theta}^{-1},\omega]  &= & 0 , \  \forall \omega \in
	 \Oalg (\mathbb R^{2N}_{\Theta})\\
	 \\
	 d({\mathrm{det}}^{-1}_{\theta}) &=& -({\mathrm{det}}^{-1}_{\theta})^2
	 d({\mathrm{det}}_{\theta})
	 \end{array}
	 \]
	 one defines the graded-involutive differential algebra
	 $\Oalg(GL_{\theta}(2n,\mathbb R))$ (writing
	 $\Omega^0_{\alg}(GL_{\theta}(2n,\mathbb
	 R))=\calg(GL_{\theta}(2n,\mathbb R))$, etc.)  which is naturally a
	 graded-involutive differential bialgebra and it is easy to show
	 that the antipode $S$ of $\calg(GL_{\theta}(2n,\mathbb R))$ extends
	 (uniquely) as an antipode, again denoted by $S$, of
	 $\Oalg(GL_{\theta}(2n,\mathbb R))$. One proceeds similarily to
	 define $\Oalg(GL_{\theta}(2n+1,\mathbb R))$. One  thus gets the
	 following result.

	 \begin{theo}\label{theo2}
	  Let $m$ be either $2n$ or $2n+1$. Then the differential algebra
	  $\Oalg(GL_{\theta}(m,\mathbb R))$ is a graded-involutive
	  differential Hopf algebra and $\Oalg(\mathbb R^m_{\theta})$
is canonically
	  a graded-involutive differential comodule over
	  $\Oalg(GL_{\theta}(m,\mathbb R))$.
	  \end{theo}

\noindent Let $G_{\theta}$ be any of the quantum groups defined in
	  Sections~\ref{sec5}
	  and \ref{sec6}. Then $\calg(G_{\theta})$ is a $\ast$-Hopf
algebra which is a
	 quotient of $\calg(GL_{\theta}(m,\mathbb R))$ by a real Hopf ideal
	 $I(G_{\theta})$ for $m=2n$ or $m=2n+1$. Let $[I(G_{\theta})]$ be the
	 closed graded two-sided ideal of $\Oalg(GL_{\theta}(m,\mathbb R))$
	 generated by $I(G_{\theta})$ and let $\Oalg(G_{\theta})$ be the
	 quotient of $\Oalg(GL_{\theta}(m,\mathbb R))$ by $[I(G_{\theta})]$.
	 The above result has the following corollary.

	 \begin{corol}\label{corol2}
	  The differential algebra $\Oalg(G_{\theta})$ is a
	  graded-involutive differential Hopf algebra and $\Oalg(\mathbb
	  R^m_{\theta})$ is a graded-involutive differential comodule over
	  $\Oalg(G_{\theta})$.
	  \end{corol}

\noindent Similarly the algebra $\Oalg(S^m_{\theta})$ is a graded-involutive
	  differential comodule over $\Oalg(SO_{\theta}(m+1))$ and a similar
result holds for
	  $GL_{\theta}(n,\mathbb C)$, $m=2n$ and
	  $\Oalg(\mathbb C^n_{\theta})$ $=\Oalg(\mathbb
	  R^{2n}_{\theta})$.\\

	       \section{The splitting homomorphisms} \label{sec8}
	       \setcounter{equation}{0}
	       We let $\pol(T^n_{\theta})$ be the
	       $\ast$-algebra of polynomials on the
	       noncommutative $n$-torus $T^n_{\theta}$ i. e. the unital
$\ast$-algebra
	       generated by $n$ unitary elements $U^\mu$ with relations
	       \begin{equation} U^\mu U^\nu = \lambda^{\mu\nu} U^\nu
U^\mu\label{ao}
	       \end{equation} for $\mu, \nu=1,\dots, n$.
	       We denote by
	       $s\mapsto \tau_{s} \in \aut(\calg(T^n_{\theta}))$ the
natural action of $T^n$ on $T^n_{\theta}$ (\cite{connes:01})
	       such that $\tau_{s}(U^\mu)=e^{2\pi i s_{\mu}}U^\mu$
$\forall s\in
	       T^n$ and $\mu \in \{1,\dots, n\}$.\\

\noindent We let as in  Section~\ref{sec2},
	        $s\mapsto \sigma_{s} \in \aut(\calg(\rnt)) $ be the
natural action of $T^n$ on $\calg(\rnt)$. It is defined for any $\theta$ (real
	       antisymmetric $(n,n)$-matrix) and in particular for $\theta=0$.
	       This yields two actions $\sigma$ and $\tau$ of $T^n$
  on $\mathbb R^{2n}\times T^n_{\theta}$ given by the
group-homomorphisms $s\mapsto \sigma_{s}\otimes
	       I$ and $s\mapsto I\otimes \tau_{s}$ of $T^n$ into
	       $\aut(\calg(\mathbb R^{2n})\otimes \calg(T^n_{\theta}))$ with
	       obvious notations. The noncommutative space $\mathbb
	       R^{2n}\times T^n_{\theta}$ is here defined by
	        $\calg(\mathbb R^{2n}\times
	       T^n_{\theta})=\calg(\mathbb R^{2n})\otimes \calg(T^n_{\theta})$.
	       We shall use the
	       actions $\sigma$ and {\sl the diagonal action} $\sigma\times
	       \tau^{-1}$ of $T^n$ on $\mathbb R^{2n}\times T^n_{\theta}$,
	       where $\sigma\times \tau^{-1}$ is defined by $s\mapsto
	       \sigma_{s}\otimes \tau_{-s}=(\sigma\times \tau^{-1})_s$ (as
	       group homomorphism of $T^n$ into $\aut(\calg(\mathbb
	       R^{2n}\times T^n_{\theta}))$).\\

\noindent In the following statement, ${z^\mu}_{(0)}$ denotes
	       the classical coordinates of $\mathbb C^n$ corresponding to
	       $z^\mu$ for $\theta=0$.

	       \begin{theo}\label{theo3}

	a) There is a unique homomorphism of unital $\ast$-algebra
	\[
	st :
	\calg(\rnt)\rightarrow \calg(\mathbb R^{2n})\otimes \calg
	(T^n_{\theta})
	\]
	such that $st(z^\mu)={z^\mu}_{(0)}\otimes
	U^\mu$ for $\mu=1,\dots,n$.\\
	b) The homomorphism $st$ induces an isomorphism of
$\calg(\rnt)$ onto the
	subalgebra $\calg(\mathbb R^{2n}\times T^n_{\theta})^{\sigma\times
	\tau^{-1}}$ of $\calg(\mathbb R^{2n}\times T^n_{\theta})$ of
  fixed points of the diagonal action of $T^n$.
	\end{theo}

\noindent One has $st(\bar z^\mu)=st(z^\mu)^\ast$ and, using (\ref{ao}),
	one checks that $st(z^\mu)$, $st(\bar z^\mu)$ fulfill the
	relations (\ref{a}). On the other hand,
it is obvious that the
	$st(z^\mu)$ are invariant by the diagonal action of $T^n$. Thus the
	only non-trivial parts of the statement, which are not difficult to
	show, are the injectivity of $st$ and the fact that $\calg(\mathbb
	R^{2n}\times T^n_{\theta})^{\sigma\times \tau^{-1}}$ is generated by
	the $z^\mu$ as unital $\ast$-algebra.

\noindent This extends trivially to
	\[
	st : \calg(\rntu)\rightarrow \calg(\mathbb R^{2n+1})\otimes \calg
	(T^n_{\theta})=\calg (\mathbb R^{2n+1}\times T^n_{\theta})
	\]
	with $st(x) ={x}_{(0)} \otimes \bbbone$ and
	$st(z^\mu)={z^\mu}_{(0)}\otimes U^\mu$. This is again an
	isomorphism of $\calg(\rntu)$ onto $\calg(\mathbb R^{2n+1}\times
	T^n_{\theta})^{\sigma\times \tau^{-1}}$.\\

\noindent The above homomorphisms $st$ pass to the quotient to define
	homomorphisms of unital $\ast$-algebras ($m=2n, 2n+1$)
        \[
	\begin{array}{l}
	 st:\calg(S^{m}_{\theta}))  \rightarrow  \calg(S^{m}) \otimes \calg(
	 T^n_{\theta})=\calg(S^{m}\times T^n_{\theta}) \end{array}
	  \]
	  which are isomorphisms of $\calg(S^{m}_{\theta})$
         with
	  $\calg(S^{m}\times T^n_{\theta})^{\sigma\times \tau^{-1}}$, the
          fixed points of the diagonal action $\sigma\times \tau^{-1}$
	  of $T^n$ (recall that $\sigma$ was previously defined for any
	  $\theta$, in particular for $\theta=0$).\\

\noindent We shall refer to the above homomorphisms $st$ as the {\sl splitting
	  homomorphisms}. They satisfy
	  $st \circ \sigma_{s}=(\sigma_{s}\otimes I)\circ st$
	  for any $s=(s_{1},\dots,s_{n})\in T^n$ and thus $st$
	  induce isomorphisms
	  \[
	  st:\calg(M_{\theta})^\sigma
	  \stackrel{\simeq}{\rightarrow} \calg (M)^\sigma\otimes
	  \bbbone (\subset \calg(M)\otimes \calg (T^n_{\theta}))
	  \]
	  for
	  $M=\mathbb R^{m}$ and $S^{m}$.\\

\noindent In a similar manner, with $M$ as above, $st$ extends
to isomorphisms of unital
	       graded-involutive differential algebras

	       \[
	       \begin{array}{l}
	       st:
	       \Oalg(M_{\theta})\rightarrow
	       (\Oalg(M)\otimes \calg(T^n_{\theta}))^{\sigma\times \tau^{-1}}
	       	       \end{array}
	         \]
	 by setting
	 \[
	 st(dz^\mu)=d{z^\mu}_{(0)}\otimes U^\mu\ \   \mbox{and} \ \
	st(dx)=d{x}_{(0)}\otimes \bbbone
	 \]
	 using the previously defined action $\sigma$ of $T^n$ on
	 $\Oalg(M_{\theta})$ for any $\theta$ (in particular $\theta=0$).

\noindent The compatibility with the differential and the action of $T^n$
is explicitly given by,

	\begin{eqnarray}
	  st \circ d  & =  & (d\otimes I)\circ st  \label{at}\\
	  st\circ \sigma_{s} & = & (\sigma_{s}\otimes I)  \circ st  \label{au}
	  \end{eqnarray}

  \noindent
	       \underbar{Remark.} We shall use the splitting
homomorphisms $st$ to reduce
computations involving $\theta$-deformations to the classical case
($\theta=0$).
  For instance we shall later define the Dirac
	       operator, $D_{\theta}$, on $M_{\theta}$ in such a way that is
	       satisfies with obvious notations $st\circ
ad(D_{\theta})=(ad(D)\otimes
	       I )\circ
	       st$ where on the right-hand side $D$ is the ordinary
Dirac operator on
	       the riemannian spin manifold $M$, ($M=\mathbb R^{2n}$, $\mathbb
	       R^{2n+1}$, $S ^{2n-1}$, $S^{2n}$); this will imply the
first order
	       condition, the reality condition and the identification of the
	       differential algebra $\Omega_{D}$ with
	       $\Omega_{\mathrm{alg}}(M_{\theta})$, (see
	       Section~\ref{sec11}).\\

\noindent A similar discussion applies to the various $\theta$-deformed groups
mentionned above.
To be specific,
we introduce the $n$ unitary elements
	       $U_{\mu}$ with relations \begin{equation} U_{\mu} U_{\nu }=
	       \lambda_{\nu\mu} U_{\nu} U_{\mu}\label{ap} \end{equation} for
	       $\mu,\nu=1,\dots,n$, (recall that $\lambda_{\mu\nu}=e^{
	       i\theta_{\mu\nu}}=\lambda^{\mu\nu}$, $\forall \mu,\nu$)
which generate $\calg(T^n_{-\theta})$,  the opposite algebra of
	       $\calg(T^n_{\theta})$.  \\

\noindent Let us consider for $m=2n$ or $m=2n+1$ the homomorphism
$ r_{23}\circ(st\otimes st)$:
	\begin{eqnarray}
	 \calg(\mathbb R^m_{\theta})\otimes
	  \calg(\mathbb R^m_{-\theta})\rightarrow \calg(\mathbb R^m)\otimes
	  \calg(\mathbb R^m)\otimes \calg(T^n_{\theta})
	  \otimes \calg(T^n_{-\theta}) \nonumber
	  \end{eqnarray}
	  where $r_{23}$ is the transposition of the
	  second and the third factors in the tensor product, (i.e.
	  $\calg(T^n_{\theta})\otimes\calg(\mathbb R^m)$ is replaced by
	  $\calg(\mathbb R^m)\otimes\calg(T^n_{\theta})$ there). This
	  $\ast$-homomorphism restricts to give a homomorphism, again
	  denoted by $st$
	  \[
	  st : M_{\theta}(m,\mathbb R)\rightarrow M(m,\mathbb R)\otimes
	  \calg(T^n_{\theta})\otimes\calg(T^n_{-\theta})
	  \]
	  which is again a homomorphism of unital $\ast$-algebras and will
	  be also refered to as splitting homomorphism. For instance, for
	  $m=2n$, it is the unique unital $\ast$-homomorphism such that
	  \begin{eqnarray}
	   st (a^\mu_{\nu})=\stackrel{(0)}{a^\mu_{\nu}}\otimes U^\mu\otimes
	   U_{\nu} \label{ar}\\
	   \nonumber\\
	   st (b^\mu_{\nu})=\stackrel{(0)}{b^\mu_{\nu}}\otimes U^\mu\otimes
	   U_{\nu}^\ast \label{as}
	   \end{eqnarray}
	   for $\mu,\nu=1,\dots, n$ where $\stackrel{(0)}{a^\mu_{\nu}}$ and
	   $\stackrel{(0)}{b^\mu_{\nu}}$ are the classical coordinates
	   corresponding to $a^\mu_{\nu}$ and $b^\mu_{\nu}$ for $\theta=0$.
	   The counterpart of $b)$ in Theorem~\ref{theo3} is that
$st$ induces here an
	   isomorphism of $M_{\theta}(m,\mathbb R)$ onto the subalgebra of
	   the elements $x$ of\linebreak[4] $M(m,\mathbb R)\otimes \calg
	   (T^n_{\theta})\otimes \calg(T^n_{-\theta})$ which are invariant by
	   the diagonal action $(\sigma\otimes\sigma)\times (\tau\otimes
	   \tau)^{-1}$ of $T^n\times T^n$ i.e. which satisfy
	   $(\sigma_{s}\otimes\sigma_{t})(\tau_{-s}\otimes \tau_{-t})(x)=x$,
	   $\forall (s,t)\in T^n\times T^n$ (with the notations of the end of
	   last section). One has
	   \[
	   st\circ (\sigma_{s}\otimes\sigma_{t})=((\sigma_{s}\otimes\sigma_{t})
	   \otimes I\otimes I)\circ st
	   \]
	   which then implies that $st$ induces an isomorphism of
$M_{\theta}(m,\mathbb
	   R)^{\sigma\otimes\sigma}$ onto \break $M(m,\mathbb 
R)^{\sigma\otimes\sigma}
	   \otimes \bbbone\otimes \bbbone$
	   where $M_{\theta}(m,\mathbb R)^{\sigma\otimes\sigma}$ denotes the
	   subalgebra of elements which are invariant by the action
	   of $T^n\times T^n$, (the same for $\theta=0$ on the
right-hand side).
	   This in particular implies that $st(\det_{\theta})$ is in
	   $M(2n,\mathbb R)^{\sigma\otimes\sigma}\otimes
	   \bbbone\otimes\bbbone$; in fact one has
	   $st(\det_{\theta})=\det\otimes \bbbone\otimes \bbbone$
	    where $\det=\det_{\theta=0}$ is the ordinary determinant.

\noindent The above homomorphism passes to the quotient to yield
	    homomorphisms
	    \[
	    st:\calg(G_{\theta})\rightarrow
	    \calg(G)\otimes\calg(T^n_{\theta})\otimes \calg(T^n_{-\theta})
	    \]
	    where $G$ is any of the classical
	    groups $GL(m,\mathbb R)$, $SL(m,\mathbb R)$, $O(m)$, $SO(m)$,
	    $GL(n,\mathbb C)$, $U(n)$ or the subgroup $GL_{(1)}(n,\mathbb
	    C)$ of $GL(n,\mathbb C)$ consisting of matrices with
	    determinants of modulus one, $m=2n$ or $m=2n+1$, and
	    where $G_{\theta}$ denote the corresponding quantum groups
	    defined in Section~\ref{sec5} and in Section~\ref{sec6}.
These homomorphisms $st$
	    which will still be refered to as the splitting homomorphisms,
	    have the property that they induce isomorphisms of
	    $\calg(G_{\theta})$ onto
	    $(\calg(G)\otimes\calg(T^n_{\theta})\otimes
	    C(T^n_{-\theta}))^{(\sigma\otimes\sigma)\times
	    (\tau\otimes\tau)^{-1}}$ for these groups $G$.\\

\noindent Thus, one sees that the situation is the same for the above
	    quantum groups as for the noncommutative spaces $M_{\theta}$
	    with $M=\mathbb R^{m}$,  $S^{m}$
	    excepted that the action of $T^n$ is replaced by an action of
	    $T^n\times T^n=T^{2n}$ and that the noncommutative $n$-torus
	    $T^n_{\theta}$ is replaced by the noncommutative $2n$-torus
	    $T^{2n}_{\theta\times(-\theta)}$ where $\theta\times(-\theta)$
	    is the real antisymmetric $(2n,2n)$-matrix $\left(
	    \begin{array}{cc} \theta & 0\\ 0 & -\theta
	     \end{array}\right)\in M_{2n}(\mathbb R)$; one has of course
	     $\calg(T^{2n}_{\theta\times(-\theta)}
	     )=\calg(T^n_{\theta})\otimes \calg(T^n_{-\theta})$.

\section{Smoothness}\label{sec9}
\setcounter{equation}{0}
Beside their usefulness for computations, the splitting homomorphisms
give straightforward unambiguous notions of smooth functions on
$\theta$-deformations.\\

\noindent The locally convex $\ast$-algebra $C^\infty(T^n_{\theta})$ of smooth
functions on the noncommutative torus $T^n_{\theta}$ was
defined in \cite{connes:01}. It is the
completion of $\pol(T^n_{\theta})$ endowed with the locally convex
topology generated by the seminorms
\[
\vert u\vert_{r}=\sup_{r_{1}+\dots+r_{n}\leq r}\parallel
X^{r_{1}}_{1}\dots X^{r_{n}}_{n}(u)\parallel
\]
where $\parallel\cdot \parallel$ is the $C^\ast$-norm (which is the
$\sup$ of the $C^\ast$-seminorms) and where the $X_{\mu}$ are the
infinitesimal generators of the action $s\mapsto \tau_{s}$ of $T^n$ on
$T^n_{\theta}$. They are the
unique derivations of $\pol(T^n_{\theta})$ satisfying
\begin{equation}
X_{\mu}(U^\nu)= 2\pi i\delta^\nu_{\mu} U^\nu
\label{ba}
\end{equation}
for $\mu, \nu=1,\dots, n$.  Notice that these derivations are real and
commute between themselves, i.e. $X_{\mu}(u^\ast)=(X_{\mu}(u))^\ast$
and $X_{\mu}X_{\nu}-X_{\nu}X_{\mu}=0$.
This locally convex
$\ast$-algebra is a nuclear Fr\'echet space and it
follows from the general theory of topological tensor products
that the $\pi$-topology and $\varepsilon$-topology
coincide \cite{gro} on any tensor product,
\cite{tre} i.e.
\[
E\otimes_{\pi }C^\infty
(T^n_{\theta})=E\otimes_{\varepsilon}
C^\infty(T^n_{\theta})
\]
so that on $E\otimes C^\infty(T^n_{\theta})$ there
is essentially one reasonable locally convex topology and we denote by
$E\widehat{\otimes}
C^\infty(T^n_{\theta})$ the
corresponding completion. \\

\noindent It is then straightforward to define the function spaces
$C^\infty(M_{\theta})$ (of smooth functions) and
$C_c^\infty(M_{\theta})$ (of smooth functions with compact support)
for any of the $\theta$-deformed spaces mentionned above,
as the fixed point algebra of the diagonal action of $T^n$
on the completed tensor product $C^\infty(M)\widehat{\otimes}
C^\infty(T^n_{\theta})$ (and on $C_c^\infty(M)\widehat{\otimes}
C^\infty(T^n_{\theta}))$.

\noindent Using the appropriate splitting homomorphisms, one defines 
in the same
way the locally convex $\ast$-algebras
$C^\infty(G_{\theta})$ and $C_c^\infty(G_{\theta})$
  of smooth functions on the different
quantum groups defined in Section~\ref{sec5} and in
Section~{\ref{sec6}. The same discussion applies to the
algebras  $\Omega(M_{\theta})$ and $\Omega_c(M_{\theta})$of smooth
differential forms.\\

\section{Differential forms, self-duality, Hochschild
cohomology for $\theta$-deformations}\label{sec10}
\setcounter{equation}{0}
Let $M$ be a smooth $m$-dimensional manifold
endowed with a smooth action $s\mapsto \sigma_{s}$ of the compact
abelian Lie group $T^n$,  (the $n$-torus). We also denote by
$s\mapsto \sigma_{s}$ the corresponding group-homomorphism of $T^n$
into the group $\aut(C^\infty(M))$ (resp $\aut(\Omega(M))$)
of automorphisms of the unital
$\ast$-algebra $C^\infty(M)$ of complex smooth functions on $M$
with its standard topology (resp of the
graded-involutive differential algebra $\Omega(M)$ of smooth
differential forms).

\noindent Let $C^\infty(M_{\theta})$ be the $\theta$-deformation of the
$\ast$-algebra $C^\infty(M)$ associated by \cite{rief:3} to the above data.
We shall find it convenient to give the following (trivially  equivalent)
direct description of
$C^\infty(M_{\theta})$ as a fixed point algebra.

\noindent The completed tensor product $C^\infty(M)\widehat\otimes
   C^\infty(T^n_{\theta})$ is
unambiguously defined by nuclearity and is a unital locally convex
   $\ast$-algebra which is a complete nuclear space. We define by
   duality the noncommutative smooth manifold $M\times T^n_{\theta}$ by
   setting $C^\infty(M\times T^n_{\theta})=C^\infty(M)\widehat\otimes
   C^\infty(T^n_{\theta})$ ; elements of $C^\infty(M\times
   T^n_{\theta})$ will be refered to as {\sl the smooth functions on
   $M\times T^n_{\theta}$}. Let $C^\infty(M\times
   T^n_{\theta})^{\sigma\times \tau^{-1}}$ be the subalgebra of
    the $f\in C^\infty(M\times T^n_{\theta})$ which
   are invariant by the diagonal action $\sigma\times\tau^{-1}$ of
   $T^n$, that is such that $\sigma_{s}\otimes \tau_{-s}(f)=f$ for any
   $s\in T^n$. One defines by duality the noncommutative manifold
   $M_{\theta}$ by setting $C^\infty(M_{\theta})=C^\infty(M\times
   T^n_{\theta})^{\sigma\times \tau^{-1}}$  and the elements of
   $C^\infty(M_{\theta})$ will be refered to as {\sl the smooth
   functions on $M_{\theta}$}. This definition clearly
   coincides with the one used before for the examples of the
   previous sections once identified using the splitting
   homomorphisms.\\

\noindent Let us now give a first construction of smooth differential forms on
   $M_{\theta}$ generalizing the one given before in the examples.
   Let $\Omega(M_{\theta})$ be the graded-involutive subalgebra
   $(\Omega(M)\widehat\otimes C^\infty(T^n_{\theta}))^{\sigma\times
   \tau^{-1}}$  of $\Omega(M)\widehat\otimes C^\infty (T^n_{\theta})$
   consisting of elements which are invariant by the diagonal
   action $\sigma\times \tau^{-1}$ of $T^n$. This subalgebra is stable
   by $d\otimes I$ so $\Omega(M_{\theta})$ is a locally convex
   graded-involutive differential algebra which is a deformation of
   $\Omega(M)$ with $\Omega^0(M_{\theta})=C^\infty(M_{\theta})$ and
   which will be refered to as the {\sl algebra of smooth differential
   forms on $M_{\theta}$}.  The action $s\mapsto\sigma_{s}$ of $T^n$ on
   $\Omega(M)$ induces $s\mapsto \sigma_{s}\otimes I$ on
   $\Omega(M)\widehat \otimes C^\infty (T^n_{\theta})$ which gives by
   restriction a group-homomorphism, again denoted $s\mapsto
   \sigma_{s}$, of $T^n$ into the group $\aut(\Omega(M_{\theta}))$ of
   automorphisms of the graded-involutive differential algebra
   $\Omega(M_{\theta})$.

   \begin{proposition}\label{prop3}
    The graded-involutive differential subalgebra
    $\Omega(M_{\theta})^\sigma$ of $\sigma$-invariant elements of
    $\Omega(M_{\theta})$ is in the graded center of $\Omega(M_{\theta})$
    and identifies canonically with the graded-involutive differential
    subalgebra $\Omega(M)^\sigma$ of $\sigma\mbox{-}$invariant elements of
    $\Omega(M)$.
    \end{proposition}
    In other words the subalgebra of $\sigma$-invariant elements of
    $\Omega(M_{\theta})$ is not deformed (i.e. independent of
    $\theta$). One has
    $\Omega(M_{\theta})^\sigma=\Omega(M)^\sigma\otimes \bbbone$
    ($\subset \Omega(M)\widehat\otimes C^\infty(T^n_{\theta})$).

\noindent The notations $M_{\theta}$,
   $C^\infty(M_{\theta})$ introduced here are coherent with the
   standard ones $T^n_{\theta}$, $C^\infty(T^n_{\theta})$ used for the
   noncommutative torus. Indeed it is true that one has
   $C^\infty(T^n_{\theta})=(C^\infty(T^n)\widehat\otimes
   C^\infty(T^n_{\theta}))^{\sigma\times \tau^{-1}}$  where
   $\sigma$ is the canonical action of $T^n$ on itself. Furthermore
   there is a natural definition of the graded differential algebra of
   smooth differential forms on the noncommutative $n$-torus
   $T^n_{\theta}$ \cite{connes:01} and it turns out that it coincides
   with the above one for $M=T^n$, that is with $\Omega(T^n_{\theta})$,
   as easily verified.\\

\noindent Although simple and useful, the previous definition of  smooth
   differential forms on $M_{\theta}$ is not the most natural one.
   Indeed the construction has  the
   following geometric interpretation. The noncommutative
   manifold $M_{\theta}$ is the quotient of the product $M\times
   T^n_{\theta}$ by the diagonal action of $T^n$, and one has a
   noncommutative fibre bundle
   \[
   M\times T^n_{\theta}\stackrel{T^n}{\longrightarrow} M_{\theta}
   \]
   with fibre $T^n$. In such a context it is natural to describe
   differential forms on $M_{\theta}$ as the basic forms on $M\times
   T^n_{\theta}$ for the operation of $\lie(T^n)$ corresponding to the
   infinitesimal diagonal action of $T^n$.
  More precisely, let
$Y_{\mu}$, $\mu\in\{1,\dots,n\}$ be the vector fields on $M$
corresponding to the infinitesimal action of $T^n$
\begin{equation}
   Y_{\mu}(x) =\frac{\partial}{\partial s_{\mu}}\sigma_{s}(x)\mid_{s=0}
   \label{bb}
\end{equation}
for $x\in M$. These vector fields are real and define $n$
derivations of $C^\infty(M)$, again denoted by $Y_{\mu}$, which are
real and commute between themselves. The inner anti-derivations
$Y_{\mu}\mapsto i_{Y_{\mu}}$ define {\sl an operation of the (abelian)
Lie algebra $\lie(T^n)$ in the graded differential algebra}
$\Omega(M)$ \cite{cart:01}, \cite{ghv} and the corresponding Lie
derivatives $L_{Y_{\mu}}=di_{Y_{\mu}}+i_{Y_{\mu}}d$ are derivations of
degree zero of $\Omega(M)$ which extend the $Y_{\mu}$ and correspond
to the infinitesimal action of $T^n$ on $\Omega(M)$.
  The natural
   graded differential algebra of smooth differential forms on $M\times
   T^n_{\theta}$ is $\Omega(M\times
   T^n_{\theta})=\Omega(M)\widehat\otimes_{\gr}  \Omega(T^n_{\theta})$, and the
   operation \cite{cart:01}, \cite{ghv} of $\lie(T^n)$ in
   $\Omega(M\times T^n_{\theta})$ corresponding to the diagonal action
   of $T^n$ is described by the antiderivations
   $i_{\mu}=i_{Y_{\mu}}\otimes I - (-I)^{\gr}\otimes i_{X_{\mu}}$ of
   $\Omega(M\times T^n_{\theta})$ where $i_{X_{\mu}}$ is the
   antiderivation of degree -1 of
   $\Omega(T^n_{\theta})=C^\infty(T^n_{\theta})\otimes_{\mathbb
   R}\wedge \mathbb R^n$ \cite{connes:01} such that
   $i_{X_{\mu}}(\omega^\nu)=\delta_{\mu}^\nu$ with
   $\omega^\mu=\frac{1}{2\pi i} U^{\mu\ast} dU^\mu$. The infinitesimal
   diagonal action of $T^n$ is described by the Lie derivatives
   $L_{\mu}=di_{\mu}+i_{\mu}d$ on $\Omega(M\times T^n_{\theta})$ and
   the differential subalgebra $\Omega_{B} (M\times T^n_{\theta})$ of
   the basic elements of $\Omega(M\times T^n_{\theta})$, that is of the
   elements $\alpha$ satisfying $i_{\mu}(\alpha)=0$ and
   $L_{\mu}(\alpha)=0$ for $\mu\in\{1 ,\cdots,n\}$, is a natural
   candidate to be the algebra of smooth differential forms on
   $M_{\theta}$. Fortunately, it is not hard to show that one has the
   following result which allows to use either point of view.

   \begin{proposition}\label{prop4}
    As graded-involutive differential algebra $\Omega_{B}(M\times
    T^n_{\theta})$ is isomorphic to $\Omega(M_{\theta})$.
    \end{proposition}

\noindent The  (first) construction of $\Omega(M_{\theta})$ admits the
    following generalization. Let $S$ be a smooth complex vector bundle
    of finite rank over $M$ and let $C^\infty(M,S)$ be the
    $C^\infty(M)$-module of its smooth sections,  endowed
    with its usual topology of complete nuclear space.
The vector bundle $S$ will be called
    $\sigma$-{\sl equivariant} if it is endowed with a group-homomorphism
    $s\mapsto V_{s}$ of $T^n$ into the group $\aut(S)$ of automorphisms
    of $S$ which covers the action $s\mapsto \sigma_{s}$ of $T^n$ on
    $M$. In terms of smooth sections this means that one has
    \begin{equation}
     V_{s}(f\psi)=\sigma_{s}(f) V_{s}(\psi)
     \label{bc}
     \end{equation}
     for $f\in C^\infty(M)$ and $\psi\in C^\infty(M,S)$ with an obvious
     abuse of notations. Let $C^\infty(M_{\theta},S)$ be the closed
     subspace of $C^\infty(M,S)\widehat\otimes C^\infty(T^n_{\theta})$
     consisting of elements $\Psi$ which are invariant by the
     diagonal action $V\times\tau^{-1}$ of $T^n$, i.e. which satisfy
     $V_{s}\otimes \tau_{-s}(\Psi)=\Psi$ for any $s\in T^n$. The locally
     convex space $C^\infty(M_{\theta},S)$
     is also canonically a  topological bimodule over
     $C^\infty(M_{\theta})$, or which is the same, a topological left
     module over $C^\infty(M_{\theta})\widehat\otimes
     C^\infty(M_{\theta})^{opp}$.
     \begin{proposition}\label{prop5}
     The bimodule $C^\infty(M_{\theta},S)$ is diagonal and
     (topologically) left and right finite projective over
     $C^\infty(M_{\theta})$.
     \end{proposition}
   The proof of this proposition uses the equivalence between the
   category of $\sigma$-equivariant finite projective modules over
   $C^\infty(M)$ (i.e. of $\sigma$-equivariant vector bundles over $M$)
   and the category of finite projective modules over the cross-product
   $C^\infty(M)\rtimes_{\sigma}T^n$, the fact that one has
   $C^\infty(M)\rtimes_{\sigma}T^n\simeq C^\infty(M_{\theta}
)\rtimes_{\sigma}T^n$
   and finally the equivalence between the category of finite
   projective modules over $C^\infty(M_{\theta})\rtimes_{\sigma}T^n$ and
   the category of $\sigma$-equivariant finite projective modules over
   $C^\infty(M_{\theta})$ \cite{pj}.\\

\noindent Let $D$ be a continuous
     $\mathbb C$-linear operator on $C^\infty(M,S)$
     such that
     \begin{equation}
      DV_{s}=V_{s}D
      \label{bd}
      \end{equation}
      for any $s\in T^n$. Then $C^\infty(M_{\theta},S)$ ($\subset
      C^\infty(M,S)\widehat\otimes C^\infty(T^n_{\theta})$) is stable
      by $D\otimes I$ which defines the operator $D_{\theta}$
      ($=D\otimes I\restriction C^\infty(M_{\theta},S)$) on
      $C^\infty(M_{\theta},S)$. If $D$ is a first-order differential
      operator it follows immediately from the definition that
      $D_{\theta}$ is a first-order operator of the bimodule
      $C^\infty(M_{\theta},S)$ over $C^\infty(M_{\theta})$ into itself,
      \cite{connes:03}, \cite{mdv:m}. If $D$ is of order zero i. e. is a
      module homomorphism over $C^\infty(M)$ then it is obvious
      that $D_{\theta}$ is a bimodule homomorphism  over
      $C^\infty(M_{\theta})$.\\

\noindent We already met this construction in the case of $S=\wedge T^\ast
      M$ and $D=d$. There  $D_{\theta}$ is the
      differential $d$ of $\Omega(M_{\theta})$ which is a first-order
      operator on the bimodule $\Omega(M_{\theta})$ over
      $C^\infty(M_{\theta})$. Let $\omega\mapsto \ast\omega$ be the
Hodge operator
      on $\Omega(M)$ corresponding to a $\sigma$-invariant riemannian
      metric on $M$. One has  $\ast\circ \sigma_{s}=\sigma_{s}\circ \ast$
      thus $\ast$ satisfies (\ref{bd}) from which one obtains an
      endomorphism $\ast_{\theta}$ of $\Omega(M_{\theta})$ considered as
      a bimodule over $C^\infty(M_{\theta})$. We shall denote
      $\ast_{\theta}$ simply by $\ast$ in the following. One has
      $\ast\Omega^p(M_{\theta})\subset \Omega^{m-p}(M_{\theta})$.

      \begin{theo}\label{theo4}
       Let the $2n$-sphere $S^{2n}$ be endowed with its usual metric,
       let $\ast$ be defined as above on $\Omega(S^{2n}_{\theta})$ and
       let $e$ be the hermitian projection of  Theorem~\ref{theo1}. Then $e$
       satisfies the self-duality equation $\ast e(de)^n=i^n e(de)^n$.
       \end{theo}
       Indeed using the splitting homomorphism, $e$ identifies with
       \[
       e=\frac{1}{2} (\bbbone + \sum^n_{\mu=1}(u^\mu_{(0)}
       \tilde \Gamma^{\mu\ast}+u^{\mu\ast}_{(0)} \tilde
       \Gamma^\mu)+u\gamma)
       \]
       where $u^\mu_{(0)} ,\cdots,u$ are
       now the classical coordinates of $\mathbb R^{2n+1}$ for
       $S^{2n}\subset \mathbb R^{2n+1}$ and where $\tilde
       \Gamma^{\mu\ast}=\Gamma^{\mu\ast}\otimes U^\mu$,
       $\tilde\Gamma^\mu=\Gamma^\mu\otimes U^{\mu\ast}$ with $\gamma$
       identified with $\gamma\otimes \bbbone \in
       M_{2^n}(C^\infty(T^n_{\theta}))$. Now one verifies easily that
       the $\tilde\Gamma^{\mu\ast}$, $\tilde\Gamma^\nu$ satisfy the
       relations of the usual Clifford algebra of $\mathbb R^{2n}$ so
       $\ast e(de)^n=i^n e(de)^n$ follows from the classical relation
       (\ref{sd}) for $P_{+}=e\mid_{\theta=0}$ and from $\ast=\ast\otimes I$
       where on the right-hand side $\ast$ is the classical one.\\

\noindent Similarily one has $\ast e_{-}(de_{-})^n=-i^n e_{-}(de_{-})^n$.
       Notice that if one replaces the usual metric of $S^{2n}$ by
       another $\sigma$-invariant metric which is conformally
       equivalent, the same result holds but that $\sigma$-invariance is
       a priori necessary for this.

\noindent Let us now compute the Hochschild dimension of $M_{\theta}$.
We first construct a continuous projective resolution of the
       left $C^\infty(M_{\theta})\widehat\otimes
       C^\infty$ $(M_{\theta})^{opp}$-module $C^\infty(M_{\theta})$.

      \begin{lemma}\label{lem3}
       There are continuous homomorphisms  of left modules
       \[
       i_{p}:\Omega^p(M_{\theta})\widehat\otimes C^\infty(M_{\theta})\rightarrow
       \Omega^{p-1}(M_{\theta})\widehat\otimes C^\infty(M_{\theta})
       \]
       over
       $C^\infty(M_{\theta})\widehat\otimes
       C^\infty(M_{\theta})^{opp}$ for $p\in\{1,\cdots,m\}$ such that
the sequence
      \[
       0\rightarrow\Omega^m(M_{\theta})\widehat\otimes
       C^\infty(M_{\theta})\stackrel{i_{m}}{\rightarrow}\cdots
      \stackrel{i_{1}}{\rightarrow}
       C^\infty(M_{\theta})\widehat\otimes
       C^\infty(M_{\theta})\stackrel{\mu}{\rightarrow}
       C^\infty(M_{\theta})\rightarrow 0
       \]
       is exact, where $\mu$ is induced by the product of
       $C^\infty(M_{\theta})$.

       \end{lemma}

\noindent In fact as was shown and used in \cite{connes:02} one has
       continuous projective resolutions of $C^\infty(M)$ and of
       $C^\infty(T^n_{\theta})$ of the form
       \[
       0\rightarrow \Omega^m(M)\widehat\otimes
       C^\infty(M)\stackrel{i^0_{m}}{\rightarrow} \cdots
       \stackrel{i^0_{1}}{\rightarrow} C^\infty(M)\widehat\otimes
       C^\infty(M)\stackrel{\mu}{\rightarrow}C^\infty(M)\rightarrow 0
       \]
        \[
       0\rightarrow \Omega^n(T^n_{\theta})\widehat\otimes
       C^\infty(T^n_{\theta})\stackrel{j_{n}}{\rightarrow} \cdots
       \stackrel{j_{1}}{\rightarrow} C^\infty(T^n_{\theta})\widehat\otimes
   C^\infty(T^n_{\theta})\stackrel{\mu}{\rightarrow}C^\infty(T^n_{\theta}
)\rightarrow 0
       \]
       which combine to give a continuous projective resolution of
     \[
       C^\infty(M)\widehat\otimes
       C^\infty(T^n_{\theta})=C^\infty(M\times T^n_{\theta})
       \]
       of the form
       \[
       \begin{array}{l}
       0\rightarrow \Omega^{m+n}(M\times T^n_{\theta})\widehat \otimes
       C^\infty (M\times T^n_{\theta})\stackrel{\tilde
       i_{m+n}}{\rightarrow}\cdots\\
       \stackrel{\tilde
       \imath_{1}}{\rightarrow} C^\infty(M\times
       T^n_{\theta})\widehat\otimes C^\infty(M\times
       T^n_{\theta})\stackrel{\mu}{\rightarrow}C^\infty(M\times
       T^n_{\theta})\rightarrow 0
       \end{array}
       \]
       where $\Omega^p(M\times T^n_{\theta})=\oplusinf_{p\geq k\geq 0}
       \Omega^k(M)\widehat \otimes \Omega^{p-k}(T^n_{\theta})$ and
       where
       \[
       \tilde\imath_{p}=\sum_{k}(i^0_{k}\otimes I+(-I)^k\otimes
       j_{p-k}).
       \]
       There is some freedom in the choice of the
       $i^0_{k}$, $j_{\ell}$ and one can choose them equivariant (by
       choosing a $\sigma$-invariant metric on $M$, etc.) in such a way
       that the $\tilde\imath_{p}$ restrict as continuous homomorphisms
       \[
       i_{p}:\Omega^p_{B}(M\times T^n_{\theta})\widehat \otimes
       C^\infty (M_{\theta})\rightarrow \Omega^{p-1}_{B}(M\times
T^n_{\theta})\widehat \otimes
       C^\infty (M_{\theta})
       \]
       of $C^\infty(M_{\theta})\widehat\otimes
       C^\infty(M_{\theta})^{opp}$-modules which gives the desired resolution
       of $C^\infty(M_{\theta})$ using Proposition~\ref{prop5}.\\

\noindent This shows that the Hochschild dimension $m_{\theta}$ of $M_{\theta}$
is $\leq m$ where
$m$ is the dimension of $M$.

\noindent Let
        $w \in\Omega^m(M)$ be a non-zero $\sigma$-invariant form of degree $m$
        on $M$ (obtained by a straightforward local averaging). In view of
Proposition~\ref{prop3}, $w
        \otimes \bbbone=w_{\theta}$ is a $\sigma\mathrm{-invariant}$ element
        of $\Omega^m(M_{\theta})$, i.e. $w_{\theta}\in
        \Omega^m(M_{\theta})^\sigma$ which defines canonically
        a non-trivial invariant cycle $v_{\theta}$ in
        $Z_{m}(C^\infty(M_{\theta}), C^\infty(M_{\theta}))$. Thus one
has $m_{\theta}\geq m$
        and therefore the following result.

        \begin{theo}\label{theo5}
         Let $M_{\theta}$ be a $\theta$-deformation of $M$ then one has $
         \dim(M_{\theta})=\dim(M)$, that is the Hochschild dimension
         $m_{\theta}$ of $C^\infty(M_{\theta})$ coincides with the
dimension $m$ of $M$.
       \end{theo}

\noindent Note that the conclusion of the theorem fails for general
deformations by actions of
  $\mathbb R^d$ as described in \cite{rief:3}.
Indeed, in the simplest case of the Moyal deformation
of $\mathbb R^{2n}$ the Hochschild dimension drops down to zero
for non-degenerate values of the deformation parameter.
It is however easy to check that periodic cyclic cohomology
(but not its natural filtration) is unaffected by the $\theta$-deformation.

\section{Metric aspect: The spectral triple}\label{sec11}
\setcounter{equation}{0}
As in the last section we let $M$ be a smooth
$m$-dimensional manifold endowed with a smooth action $s\mapsto
\sigma_{s}$ of $T^n$.
It is well-known and easy to check that
we can average any
    riemannian metric on $M$ under the action of $\sigma$
and obtain one for which the action
   $s\mapsto \sigma_{s}$ of $T^n$ on $M$ is isometric.
   Let us  assume moreover that $M$ is a
  spin manifold.
Let $S$ be the spin bundle over $M$ and let $D$ be the Dirac operator
on $C^\infty(M,S)$. The bundle $S$
is not $\sigma$-equivariant in the sense of the last section but is
equivariant in a slightly generalized sense which we now explain. In
fact the isometric action $\sigma$ of $T^n$ on $M$ does not lift
directly to $S$ but lifts only modulo $\pm I$. More precisely one has
a twofold covering $p:\tilde T^n\rightarrow T^n$ of the group $T^n$,
and a group homomorphism $\tilde s\mapsto V_{\tilde s}$ of
$\tilde T^n$ into the group $\aut(S)$ which covers the action
$s\mapsto \sigma_{s}$ of $T^n$ on $M$. In terms of smooth sections,
(\ref{bc}) generalizes here as
\begin{equation}
   V_{\tilde s}(f\psi)=\sigma_{s}(f) V_{\tilde s}(\psi)
   \label{bcd}
   \end{equation}
   where $f\in C^\infty(M)$ and $\psi\in C^\infty(M,S)$ with $s=p(\tilde
   s)$. The bundle $S$ is also a hermitian vector bundle and one has
\begin{equation}
   (V_{\tilde s}(\psi),V_{\tilde s}(\psi'))=\sigma_{s}((\psi,\psi'))
   \label{be}
   \end{equation}
   for $\psi,\psi'\in C^\infty(M,S)$, $\tilde s\in \tilde T^n$ and
   $s=p(\tilde s)$ where (.,.) denotes
   the hermitian scalar product.  Furthermore, the Dirac operator $D$
   commutes with the $V_{\tilde s}$.\\

\noindent To the projection $p:\tilde T^n\rightarrow T^n$ corresponds an
   injective homomorphism of $C^\infty(T^n)$ into $C^\infty(\tilde
   T^n)$ which identifies $C^\infty(T^n)$ with the
   subalgebra\linebreak[4]
   $C^\infty (\tilde T ^n)^{\ker(p)}$ of $C^\infty(\tilde T^n)$ of
   elements which are invariant by the action of the subgroup
   $\ker(p)\simeq \mathbb Z_{2}$ of $\tilde T^n$. Let $\tilde
   T^n_{\theta}$ be the noncommutative $n$-torus
   $T^n_{\frac{1}{2}\theta}$ and let $\tilde s\mapsto \tilde
   \tau_{\tilde s}$ be the canonical action of the $n$-torus $\tilde
   T^n$ that is the canonical group-homomorphism of $\tilde T^n$ into
   the group $\aut (C^\infty(\tilde T^n_{\theta}))$. The very reason for
   these notations is that $C^\infty(T^n_{\theta})$ identifies with the
   subalgebra $C^\infty(\tilde T^n_{\theta})^{\ker(p)}$ of
   $C^\infty(\tilde T^n_{\theta})$ of elements which are invariant
   by the $\tilde \tau_{\tilde s}$ for $\tilde s\in \ker(p)\simeq
   \mathbb Z_{2}$. Under this identification, one has $\tilde
   \tau_{\tilde s}(f)= \tau_{s}(f)$ for $f\in C^\infty(T^n_{\theta})$ and
   $s=p(\tilde s)\in T^n$.\\

\noindent Define $C^\infty(M_{\theta},S)$ to be the closed subspace of
   $C^\infty(M,S)\widehat \otimes C^\infty(\tilde T_{\theta})$
   consisting of elements $\Psi$ which are invariant by the
   diagonal action $V\times \tilde \tau^{-1}$ of $\tilde T^n$;  this is
   canonically a topological bimodule over $C^\infty(M_{\theta})$. Since
   the Dirac operator commutes with the $V_{\tilde s}$,
   $C^\infty(M_{\theta},S)$ is stable by $D\otimes I$ and we denote
   by $D_{\theta}$ the corresponding operator on
   $C^\infty(M_{\theta},S)$. Again, $D_{\theta}$ is a
   first-order operator of the bimodule $C^\infty(M_{\theta},S)$ over
   $C^\infty(M_{\theta})$ into itself. The space
   $C^\infty(M,S)\widehat \otimes C^\infty(\tilde T^n_{\theta})$ is
   canonically a bimodule over $C^\infty(M)\widehat \otimes C^\infty
   (\tilde T^n_{\theta})$ (and therefore also on $C^\infty(M)\widehat
\otimes C^\infty
   (T^n_{\theta})$). One defines a hermitian structure on
   $C^\infty(M,S)\widehat \otimes C^\infty(\tilde T^n_{\theta})$ for
its right-module structure over
   $C^\infty(M)\widehat\otimes C^\infty(\tilde T^n_{\theta})$
\cite{connes:01} by setting
   \[
   (\psi\otimes t, \psi'\otimes t')= (\psi,\psi')\otimes t^\ast t'
   \]
   for $\psi,\psi'\in C^\infty(M,S)$ and $t,t'\in
   C^\infty(\tilde T^n_{\theta})$. This gives by restriction the hermitian
   structure of $C^\infty(M_{\theta},S)$ considered as a right
   $C^\infty(M_{\theta})$-module; that is one has
   \[
   (\psi f, \psi' f')= f^\ast(\psi,\psi') f'
   \]
   for any $\psi, \psi'\in C^\infty(M_{\theta},S)$ and $f,f'\in
   C^\infty(M_{\theta})$. Notice that when $\mbox{dim}(M)$ is even, one
   has a $\mathbb Z_{2}$-grading $\gamma$ of $C^\infty(M,S)$ as
   hermitian module which induces a $\mathbb Z_{2}$-grading,  again
   denoted by $\gamma$, of $C^\infty(M_{\theta},S)$ as hermitian right
   $C^\infty(M_{\theta})$-module.\\

\noindent Let $J$ denote the charge conjugation of $S$. This is an antilinear
   mapping of $C^\infty(M,S)$ into itself such that
   \begin{eqnarray}
    (J\psi,J\psi) & = & (\psi,\psi)
    \label{bf}\\
    JfJ^{-1} & = & f^\ast
    \label{bg}
    \end{eqnarray}
    for any $\psi\in C^\infty(M,S)$ and for any $f\in C^\infty(M)$,
    $(f^\ast(x)=\overline{f(x)})$. Furthermore one has also
    \begin{equation}
     JV_{\tilde s}=V_{\tilde s} J
\label{bh}
\end{equation}
for any $\tilde s\in \tilde T^n$. Let us define $\tilde J$ to be the unique
antilinear operator on $C^\infty(M,S)\widehat\otimes
C^\infty(\tilde T^n_{\theta})$ satisfying $\tilde J(\psi\otimes
t)=J\psi\otimes t^\ast$ for $\psi\in C^\infty(M,S)$ and $\tilde t\in
C^\infty(\tilde T^n_{\theta})$. The subspace $C^\infty(M_{\theta},S)$ is
stable by $\tilde J$ and we define $J_{\theta}$ to be the induced
antilinear mapping of $C^\infty(M_{\theta},S)$ into itself. It follows
from (\ref{bf}), (\ref{bg}) and from the definition that one has
\begin{eqnarray}
   (J_{\theta}\psi,J_{\theta}\psi)& = &(\psi,\psi) \label{bi}\\
   J_{\theta}f J^{-1}_{\theta}\psi & = & \psi f^\ast
   \label{bj}
   \end{eqnarray}
   for any $\psi\in C^\infty (M_{\theta},S)$ and $f\in
   C^\infty(M_{\theta})$. Thus left multiplication by $J_{\theta}f^\ast
   J_{\theta}^{-1}$ is the same as right multiplication by $f$. Obviously
   $J_{\theta}$ satisfies, in function of $\mbox{dim}(M)$ modulo 8, the
   table of normalizations, commutations with $D_{\theta}$ and with
   $\gamma$ in the even dimensional case which corresponds to the
   reality conditions 7) of \cite{connes:61}.  This follows of course
   from the same properties of $J, D, \gamma$ (i.e. the same properties
   for $\theta=0$). So equipped $C^\infty(M_{\theta},S)$ is in
   particular an involutive bimodule with a right-hermitian structure
   \cite{rief:2}, \cite{mdv:pm2}.\\

\noindent Let us now investigate the symbol of $D_{\theta}$. It is easy to see
   that the left universal symbol $\sigma_{L}(D_{\theta})$ of
   $D_{\theta}$ (as defined in \cite{mdv:m}) factorizes through a
   homomorphism
   \[
   \hat\sigma_{L}(D_{\theta}):\Omega^1(M_{\theta})\otimesinf_{C^\infty(M_
{\theta})}
     C^\infty(M_{\theta},S)\rightarrow C^\infty(M_{\theta},S)
     \]
     of bimodules over $C^\infty(M_{\theta})$. By definition, one has
     \[
     [D_{\theta},f]\psi = \hat\sigma_{L}(D_{\theta}) (df\otimes \psi)
     \]
     for $f\in C^\infty(M_{\theta})$ and
     $\psi\in C^\infty(M_{\theta},S)$ and $df\mapsto [D_{\theta},f]$
     extends as an injective linear mapping of $\Omega^1(M_{\theta})$
     into the continuous linear endomorphisms of
     $C^\infty(M_{\theta},S)$.

     \begin{lemma} \label{lem4}
      Let $f_{i}, g_{i}$ be a finite family of elements of
      $C^\infty(M_{\theta})$ such that $\sum_{i}f_{i}[D_{\theta},g_{i}]=0$.
      Then the endomorphism
$\sum_{i}[D_{\theta},f_{i}][D_{\theta},g_{i}]$ is the left
      multiplication in $C^\infty(M_{\theta},S)$ by an element of
      $C^\infty(M_{\theta})$.
      \end{lemma}
      When no confusion arises, we shall summarize this statement by
      writing \break
      $\sum_{i}[D_{\theta}, f_{i}][D_{\theta},g_{i}]\in
      C^\infty(M_{\theta})$ whenever $\sum_{i}
      f_{i}[D_{\theta},g_{i}]=0$. Indeed, using the fact that
      $D_{\theta}$ is the restriction of $D\otimes I$ where $D$ is the
      classical Dirac operator on $M$ one shows that
      \[
     \sum_{i}[D_{\theta},f_{i}][D_{\theta},g_{i}]+\sum_{i}
     f_{i}\Delta_{\theta}(g_{i})=[D_{\theta},\sum
     f_{i}[D_{\theta},g_{i}]]=0
     \]
     where $\Delta_{\theta}$ is the restriction of $\Delta\otimes I$ to
     $C^\infty (M_{\theta})$ with $\Delta$ being the ordinary Laplace
     operator on $M$ which is $\sigma$-invariant. This implies that
     $\sum_{i }f_{i }\Delta_{\theta}(g_{i})$ is in
     $C^\infty(M_{\theta})$ and therefore the result.\\

\noindent Concerning the particular case $M=\mathbb R^{2n}$ one shows the
     following result using the splitting homomorphism.

     \begin{proposition}\label{prop6}
     Let $z^\mu,\bar z^\nu\in C^\infty(\mathbb R^{2n}_{\theta})$ be as
     in Section \ref{sec2}. Then the
     $\hat\Gamma^\mu=[D_{\theta},z^\mu]$, $\hat{\bar
     \Gamma}^\nu=[D_{\theta},\bar z^\nu]$ satisfy the relations
     \[
     \hat \Gamma^\mu \hat\Gamma^\nu +\lambda^{\mu\nu} \hat \Gamma^\nu
     \hat\Gamma^\mu  =  0,\>
     \hat{\bar\Gamma}^\mu \hat{\bar \Gamma}^\nu + \lambda^{\mu\nu}
     \hat{\bar\Gamma}^\nu \hat{\bar \Gamma}^\mu   = 0,\>
     \hat{\bar\Gamma}^\mu \hat\Gamma^\nu + \lambda^{\nu\mu}
     \hat\Gamma^\nu \hat{\bar \Gamma}^\mu   = \delta^{\mu\nu}\bbbone
     \]
     where $\bbbone$ is the identity mapping of $C^\infty(\mathbb
     R^{2n}_{\theta},S)$ onto itself.
     \end{proposition}
     This $\theta$-twisted version of the generators of the Clifford
     algebra connected with the symbol of $D_{\theta}$ differs from the
one introduced in Section \ref{sec2} by
     the replacement $\lambda^{\mu\nu}\mapsto \lambda^{\nu\mu}$ and is
     the version associated with the $\theta$-twisted version
     $\wedge_{c}\mathbb R^{2n}_{\theta}$ of the exterior algebra which
     is itself behind the differential calculus $\Omega(\mathbb
     R^{2n}_{\theta})$. This is a counterpart for this example of the
     fact that $\Omega_{D_{\theta}}=\Omega(M_{\theta})$.\\

\noindent We now make contact with the axiomatic framework of \cite{connes:61}.
   To simplify the discussion we shall assume now that $M$ is a compact
   oriented $m\mbox{-}$dimensional riemannian spin manifold endowed with an
   isometric action of $T^n$, (i.e. we add compactness).
   One defines a positive definite scalar product on $C^\infty(M,S)$ by
   setting
   \[
   \langle \psi,\psi'\rangle = \int_{M}(\psi,\psi') \vol
   \]
   where $\vol$ is the riemannian volume $m$-form which is
   $\sigma$-invariant and we denote by $\calh=L^2(M,S)$ the Hilbert
   space obtained by completion. As an unbounded operator in $\calh$,
   the Dirac operator $D:C^\infty(M,S)\rightarrow C^\infty(M,S)$ is
   essentially self-adjoint on $C^\infty(M,S)$. We identify $D$ with its
   closure that is with the corresponding self-adjoint operator in
   $\calh$. The spectral triple $(C^\infty(M),\calh,D)$ together with
   the real structure $J$ satisfy the axioms of \cite{connes:61}. The
   homomorphism $\tilde s\mapsto V_{\tilde s}$ uniquely extends as a unitary
   representation of the group $\tilde T^n$ in $\calh$ which will be still
   denoted by $\tilde s\mapsto V_{\tilde s}$. On the other hand the action
$\tilde s\mapsto \tilde \tau_{\tilde s}$ of $\tilde T^n$ on
$C^\infty(\tilde T^n_{\theta})$ extends as a
unitary action again denoted by $\tilde s\mapsto \tilde \tau_{\tilde s}$
of $\tilde T^n$ on the
Hilbert space $L^2(\tilde T^n_{\theta})$ which is obtained from
$C^\infty(\tilde T^n_{\theta})$ by completion for the Hilbert norm $f\mapsto
\parallel f\parallel = \mbox{tr} (f^\ast f)^{1/2}$ where  tr is the
usual normalized trace of $C^\infty(\tilde T^n_{\theta})=
C^\infty(T^n_{\frac{1}{2}\theta})$ .  We now define the spectral triple
$(C^\infty(M_{\theta}),\calh_{\theta}, D_{\theta})$ to be the
following one. The Hilbert space $\calh_{\theta}$ is the subspace of
the Hilbert tensor product $\calh \widehat\otimes L^2(\tilde T^n_{\theta})$
which consists of elements $\Psi$ which are invariant by the
diagonal action of $\tilde T^n$, that is which satisfy $V_{\tilde s}\otimes
\tilde \tau_{-\tilde s}(\Psi)=\Psi$, $\forall \tilde s\in \tilde
T^n$.  The operator $D_{\theta}$
identifies with an unbounded operator in $\calh_{\theta}$ which is
essentially self-adjoint on the dense subspace
$C^\infty(M_{\theta},S)$. We also identify $D_{\theta}$ with its
closure that is with the self-adjoint operator which is also the
restriction to $\calh_{\theta}$ of $D\otimes I$. The antilinear
operator $J_{\theta}$ canonically extends as  anti-unitary operator
in $\calh_{\theta}$ (again denoted by $J_{\theta}$).
\begin{theo}\label{theo6}
   The spectral
triple $(C^\infty(M_{\theta}),\calh_{\theta}, D_{\theta})$ together
with the real structure $J_{\theta}$ satisfy all axioms of
noncommutative geometry of \cite{connes:61}.
\end{theo}
Notice that axiom 4) of orientability is
directly connected to the $\sigma$-invariance of the $m$-form $\vol$ on $M$.
  Consequently this form defines a
$\sigma$-invariant $m$-form on $M_{\theta}$ in view of
Proposition~\ref{prop3}
which corresponds to a $\sigma$-invariant Hochschild cycle in
$Z_{m}(A,A)$ for both $A=C^\infty(M)$ and $A=C^\infty(M_{\theta})$.
The argument for Poincar\'e duality is the same as in
\cite{connes:08}.
Finally, the isospectral nature of the deformation
$(C^\infty(M),\calh, D, J)\mapsto (C^\infty(M_{\theta}),
\calh_{\theta}, D_{\theta}, J_{\theta})$ follows immediately from the
fact that $D_{\theta}=D\otimes I$.\\

\noindent Coming back to the notations of sections 4 and 5, we can 
then return to
the noncommutative geometry of  $S^{m}_{\theta}$. \\
This geometry (with variable metric) is entirely specified by the
projection $e$, the matrix algebra (which together generate the
algebra of coordinates) and the Dirac operator which fulfill a
polynomial equation of degree $m$.

\begin{theo} \label{theo7}
Let $g$ be any $T^n$-invariant Riemannian metric on $S^m$, $m=2n$
or $m=2n-1$, whose volume form
is the same as for the round metric.
   $(i)$ Let $e\in M_{2^n}(C^\infty(S^{2n}_{\theta}))$ be the
projection of Theorem~\ref{theo1} . Then the Dirac operator $D_{\theta}$
of $S^{2n}_{\theta}$ associated to the metric $g$ satisfies
\[
\langle(e-\frac{1}{2})[D_{\theta},e]^{2n}\rangle=\gamma
\]
where $\langle\>\rangle$ is the projection on the commutant of
$M_{2^n}(\mathbb C)$.\\
\phantom{THEOREM 8 }  $(ii)$ Let $U\in
M_{2^{n-1}}(C^\infty(S^{2n-1}_{\theta}))$ be the unitary of
Theorem~\ref{theo1}.
Then the Dirac operator  $D_{\theta}$ of $S^{2n-1}_{\theta}$
associated to the metric $g$
satisfies
\[
\langle
U[D_{\theta},U^\ast]([D_{\theta},U][D_{\theta},U^\ast])^{n-1}\rangle =1
\]
where $\langle\>\rangle$ is the projection on the commutant of
$M_{2^{n-1}}(\mathbb C)$.
\end{theo}
Using the splitting homomorphism
as for Theorem~\ref{theo4}
it is enough to show that this holds
for the classical case $\theta=0$, i. e. when
$D$ is the classical Dirac operator associated to the metric $g$.\\

\noindent This result is of course a
straightforward extension of results of \cite{connes:07},
\cite{connes:08}. Since the deformed algebra $C^\infty(S^{m}_{\theta})$
is highly nonabelian the inner fluctuations of the noncommutative
metric (\cite{connes:61})
generate non-trivial internal gauge fields which compensate for the loss of
gravitational degrees of freedom imposed by the $T^n$-invariance of
the metric $g$.

\section{ Further prospect}\label{seconclu}
\setcounter{equation}{0}
We have shown that the basic $K$-theoretic equation defining
spherical manifolds admits a complete solution in dimension 3
and that for generic values of the deformation
parameters the obtained algebras of polynomials on the
deformed $\mathbb R^4_{\mathbf{u}}$ only depend on two parameters and 
are isomorphic to the algebras introduced
by Sklyanin in connection with the Yang-Baxter equation.
The spheres themselves do depend on the three initial parameters and 
we postpone their analysis to part II.

\noindent We did concentrate here on the critical values of the 
deformation parameters i.e. on the subclass of ${\theta}$-deformations
  and identified as $m$-dimensional noncommutative
spherical manifolds the noncommutative $m$-sphere $S^m_{\theta}$ for
any $m\in \mathbb N$. For this class we completed the path
from the crudest level of the algebra $\calg(S)$ of polynomial functions on $S$
  to the full-fledged structure
of noncommutative geometry \cite{connes:61}, as exemplified in theorem 9.
We showed that the basic polynomial equation fulfilled by the Dirac operator
held unaltered in the noncommutative case.
We also obtained the noncommutative analogue of the self-duality equations
and described concretely the quantum symmetry groups.\\

\noindent Needless to say
  our goal in part II will be to analyse general spherical
3-manifolds
including their smooth
structure, their differential calculus and metric aspect.
For these non-critical generic values the scale invariance inherited 
from criticality in the above examples will no longer hold. This will 
generate very interesting new phenomena. The analysis of the 
corresponding noncommutative spaces $S^3_{\mathbf{u}}$ is much more 
involved as we shall see in part II.

\section{Appendix : Relations in the noncommutative Grassmannian}
\setcounter{equation}{0}
Let ${\cal A}$ be the universal
Grassmannian   generated by the $2^2$ elements $\alpha, \beta,
\gamma, \delta$  with the relations,
\begin{equation}
U \, U^*=U^* \, U = 1,  \qquad U = \left(\begin{array} {cc}
\alpha & \beta\\
\gamma & \delta
\end{array}
\right)
\label{eq12.6}
\end{equation}
In this appendix we shall show that the intersection $\cal J$ of the
kernels of the representations $\rho$ of ${\cal A}$ such that
${\rm ch}_{\frac{1}{2}} (\rho (U)) = 0 $ is a non-trivial two sided ideal
of ${\cal A}$. Thus the odd Grassmanian $\cal B$
which was introduced in \cite{connes:08} is a nontrivial quotient of 
${\cal A}$.

\noindent Given an algebra ${\cal A}$ and elements $x_j \in {\cal A}$ we let,
\begin{equation}
\label{eq1}
[x_1 , \ldots , x_n] = \Sigma \ \varepsilon (\sigma) \, x_{\sigma
(1)} \ldots x_{\sigma (n)}
\end{equation}
where the sum is over all permutations and $\varepsilon (\sigma)$ is
the signature of the permutation.

\noindent  With the above notations, let
$\mu = [\alpha, \beta, \gamma, \delta]$.
We shall check that,
\begin{lemma}
In any representation $\rho$ of ${\cal A}$ for which
${\rm ch}_{\frac{1}{2}} (\rho (U)) = 0 $
  one has,
$ \rho ([\mu , \mu^*]) = 0  $. Moreover
$[\mu , \mu^*] \ne 0 \, $ in ${\cal A}$.
\end{lemma}
\noindent{\it Proof.}
For $y_i = \lambda_i^j \, x_j$, one has $[y_1 , \ldots , y_n] = \det
\lambda \, [x_1 , \ldots , x_n]$. This allows to
  extend the map $a \otimes b \otimes c \otimes d
\rightarrow [a,b,c,d]$ to a linear map $c$,
\begin{equation}
\label{eq32}
c:  \, \wedge^4 \, {\cal A} \rightarrow {\cal A} \, .
\end{equation}
Let us now show that, for any representation $\rho$ of ${\cal A}$ for which
${\rm ch}_{\frac{1}{2}} (\rho (U)) = 0 $, the following relation 
fulfilled by the matrix elements
$\tilde\alpha=\rho(\alpha),\dots, \tilde\delta=\rho(\delta)$,
\begin{equation}
\label{eq30}
\tilde\alpha \otimes \tilde\alpha^* + \tilde\beta \otimes 
\tilde\beta^* + \tilde\gamma \otimes
\tilde\gamma^* + \tilde\delta \otimes \tilde\delta^* = \tilde\alpha^* 
\otimes \tilde\alpha + \tilde\beta^*
\otimes \tilde\beta + \tilde\gamma^* \otimes \tilde\gamma + 
\tilde\delta^* \otimes \tilde\delta
\end{equation}
implies,
\begin{equation}
\label{eq31}
[\tilde\alpha , \tilde\beta , \tilde\gamma , \tilde\delta] \, 
[\tilde\alpha , \tilde\beta , \tilde\gamma ,
\tilde\delta]^* = [\tilde\alpha , \tilde\beta , \tilde\gamma , 
\tilde\delta]^* \, [\tilde\alpha , \tilde\beta ,
\tilde\gamma , \tilde\delta] \, .
\end{equation}
It follows from (\ref{eq30}) that,
\begin{equation}
\label{eq33}
(\tilde\alpha \wedge \tilde\beta \wedge \tilde\gamma \wedge 
\tilde\delta) \otimes (\tilde\alpha^*
\wedge \tilde\beta^* \wedge \tilde\gamma^* \wedge \tilde\delta^*) = 
(\tilde\alpha^* \wedge
\tilde\beta^* \wedge \tilde\gamma^* \wedge \tilde\delta^*) \otimes 
(\tilde\alpha \wedge \tilde\beta
\wedge\tilde \gamma \wedge\tilde \delta) \, .
\end{equation}
Indeed we view $\tilde{\cal A}=\rho(\cala)$ as a linear space and 
consider the tensor product
of exterior algebras,
\begin{equation}
\label{eq34}
\wedge \tilde{\cal A} \otimes \wedge \tilde{\cal A} \qquad \hbox{(ungraded
tensor product)}.
\end{equation}
We then take the $4^{\rm th}$ power of (\ref{eq30}) and get,
\begin{equation}
\label{eq35}
24 \, (\tilde\alpha \wedge \tilde\beta \wedge \tilde\gamma \wedge 
\tilde\delta) \otimes
(\tilde\alpha^* \wedge \tilde\beta^* \wedge \tilde\gamma^* \wedge 
\tilde\delta^*) = 24 \,
(\tilde\alpha^* \wedge \tilde\beta^* \wedge \tilde\gamma^* \wedge 
\tilde\delta^*) \otimes
(\tilde\alpha \wedge \tilde\beta \wedge \tilde\gamma \wedge \tilde\delta) \, .
\end{equation}
We can then apply $c \otimes c$ on both sides and compose
with $m : \tilde{\cal A} \otimes \tilde{\cal A} \rightarrow \tilde{\cal
A}$, the product, to get (\ref{eq31}), that is
\begin{equation}
\label{eq36}
\rho([\mu,\mu^\ast])=0
\end{equation}

\noindent  It remains to check that $[\mu , \mu^*] \ne 0 \, $ in ${\cal A}$.

\noindent  One has
\begin{equation}
\label{eq4}
M_2 (\Cb) * \, \Cb \, \Zb = M_2 ({\cal A})
\end{equation}
where the free product in the left hand side is the free algebra 
generated by $M_2 (\Cb)$ and a
unitary $U$, $U^* \, U = U \, U^* = 1$.
As above ${\cal A}$ is generated by the matrix elements of $U$,
\begin{equation}
\label{eq5}
U = \left(\begin{array}{cc}
\alpha & \beta\\
\gamma & \delta
\end{array}\right)\> \> ; \> \> \alpha,\beta,\gamma, \delta \in \cala
\end{equation}
As a linear basis of $M_2 (\Cb)$ we use the Pauli spin matrices,
which we view as a projective representation of $\Gamma =
(\Zb/2)^2$,
\begin{equation}
\label{eq7}
(0,0) \build\longrightarrow_{}^{\sigma} 1 \ , \quad (0,1)
\build\longrightarrow_{}^{\sigma} \sigma_1 \ , \quad (1,0)
\build\longrightarrow_{}^{\sigma} \sigma_2 \ , \quad (1,1)
\build\longrightarrow_{}^{\sigma} \sigma_3
\end{equation}
with $\sigma \, (a+b) = c \, (a,b) \ \sigma (a) \ \sigma (b)
\qquad \forall \, a,b \in (\Zb/2)^2$.

\noindent Since we are dealing with a free product, we have a
natural basis of $M_2 ({\cal A})$ given by the monomials,
\begin{equation}
\label{eq8}
\sigma_{i_1} \, U^{j_1} \ \sigma_{i_2} \, U^{j_2} \ldots
\sigma_{i_k} \, U^{j_k}
\end{equation}
where $i_1$ and $j_k$ can be $0$ but all other $i_{\ell} , j_{\ell}$
are $\ne 0$. The projection to ${\cal A}$ is given by,
\begin{equation}
\label{eq9}
P (T) = \frac{1}{4} \ \sum_{\Gamma} \ \sigma (a) \, T \
\sigma^{-1} (a) \, .
\end{equation}
In particular the matrix components $\alpha , \beta , \gamma ,
\delta$, of $U$ are linear combinations of the four elements,
\begin{equation}
\label{eq10}
x_a = P \, (\sigma (a) \, U) \ , \qquad a \in \Gamma = (\Zb/2)^2
\, .
\end{equation}
We want to compute $[\alpha , \beta , \gamma , \delta]$ or equivalently
  $[x_0 , x_1 , x_2 , x_3]$.

\noindent Let us first rewrite the product $x_{a_1} \, x_{a_2} \,
x_{a_3} \, x_{a_4}$, which is up to an overall coefficient $4^{-4}$,
\begin{eqnarray}
\label{eq11}
&\build\sum_{b_i}^{} \ \sigma (b_1) \, \sigma (a_1) \, U \, \sigma
(b_1)^{-1} \, \sigma (b_2) \, \sigma (a_2) \, U \, \sigma (b_2)^{-1}
\, \sigma (b_3) \, \sigma (a_3) \, U \, \sigma (b_3)^{-1} \nonumber \\
&\sigma (b_4) \, \sigma (a_4) \, U \, \sigma (b_4)^{-1}
\end{eqnarray}
as a sum of terms of the form,
\begin{eqnarray}
\label{eq13}
&\sigma (c_1) \, U \, \sigma (c_2) \, U \ldots U \, \sigma (c_4) \,
U \, \sigma (c_4)^{-1} \, \sigma (c_3)^{-1} \, \sigma (c_2)^{-1} \,
\sigma (c_1)^{-1} \nonumber \\
&\lambda \, (c_1 , \ldots , c_4) \,
\sigma (a_1) \, \sigma (a_2) \, \sigma (a_3) \, \sigma (a_4) \, .
\end{eqnarray}
where $c_1 = b_1 + a_1$, $c_2 = b_2 - b_1 + a_2$, $c_3 =
b_3 - b_2 + a_3$, $c_4 = b_4 - b_3 + a_4$ vary independently
in $\Gamma$ and
$\lambda \, (c_1 , \ldots , c_4) \in U(1)$ can be computed
using the trivial
representation, $U \rightarrow 1$ by,
\begin{eqnarray}
\label{eq133}
&\sigma
(b_1) \, \sigma (a_1) \, \sigma (b_1)^{-1} \, \sigma (b_2) \,
\sigma (a_2) \sigma (b_2)^{-1} \, \sigma (b_2) \, \sigma (a_3) \, \sigma
(b_3)^{-1} \, \sigma (b_4) \, \sigma (a_4) \, \sigma (b_4)^{-1} \nonumber \\
&= \, \lambda \, (c_1 , \ldots , c_4) \, \sigma (a_1) \, \sigma (a_2) 
\, \sigma (a_3) \, \sigma (a_4) .
\end{eqnarray}

\noindent Each term in the reduced expansion of $[x_0 , x_1 , x_2 , x_3]$
  is the sum of the above expressions multiplied by $\varepsilon 
(a)=\delta_{0 \, 1 \, 2 \, 3}^{a_1 \, a_2 \,
a_3 \, a_4}$
the signature of the permutation $\{ 0, 1, 2, 3 \} \rightarrow \{
a_1, a_2, a_3, a_4 \}$.

\noindent To see that $[x_0 , x_1 , x_2 , x_3] \ne 0$ we compute the
terms in
\begin{equation}
\label{eq19}
U^3 \, \sigma (c) \, U \, \sigma (c)^{-1} \, .
\end{equation}
Fixing $c$ there is one contribution for each of the permutation of
$\{0,1,2,3\}$ and in (\ref{eq11}) we have
\begin{equation}
\label{eq20}
b_1 = a_1 \, , \ b_2 = a_1 + a_2 \, , \ b_3 = a_1 + a_2 + a_3 \, , \
b_4 = c+ a_1 + a_2 + a_3 + a_4 \, .
\end{equation}
In $(\Zb / 2)^2 = \Gamma$ one has $a_1 + a_2 + a_3 + a_4 = 0$ so that
$b_4 = c$, $b_3 = a_4$. Since $\sigma (x)^2 = 1$ one can thus write
(\ref{eq11}) as
\begin{equation}
\label{eq21}
U \, \sigma (a_1) \, \sigma (a_1 + a_2) \, \sigma (a_2) \, U \, \sigma
(a_1 + a_2) \, \sigma (a_4) \, \sigma (a_3) \, U \, \sigma (a_4) \,
\sigma (c) \, \sigma (a_4) \, U \, \sigma (c)
\end{equation}
which we should multiply by $\varepsilon (a)$ and sum over $a$.

\noindent It is clear here that $\sigma (a_1) \, \sigma (a_1 + a_2) \,
\sigma (a_2)$ and $\sigma (a_1 + a_2) \, \sigma (a_4) \, \sigma (a_3)$
are {\em scalar} and thus commute with $U$ which allows to write
(\ref{eq21}) as follows,
\begin{equation}
\label{eq22}
U^3 \, \sigma (a_1) \, \sigma (a_1 + a_2) \, \sigma (a_2) \, \sigma
(a_1 + a_2) \, \sigma (a_4) \, \sigma (a_3) \, \sigma (a_4) \, \sigma
(c) \, \sigma (a_4) \, U \, \sigma (c) \, .
\end{equation}
One has,
\begin{equation}
\label{eq17}
\sigma (a) \, \sigma (a') \, \sigma (a)^{-1} \, \sigma (a')^{-1} =
(-1)^{\langle a,a' \rangle} \qquad \forall \, a,a' \in \Gamma
\end{equation}
using the bilinear form
with $\Zb/2$ values on $\Gamma$ given by,
\begin{equation}
\label{eq16}
\langle a , a' \rangle = \alpha \, \beta' - \alpha' \, \beta \quad
\hbox{for} \quad a = (\alpha , \beta) \, , \ a' = (\alpha' , \beta')
\in (\Zb/2)^2 \, .
\end{equation}
Permuting $\sigma (a_1 + a_2)$ with $\sigma (a_2)$ and $\sigma (a_4)$
with $\sigma (a_3)$ introduces terms in $(-1)^n$ with $n = \langle a_1
+ a_2 , a_2 \rangle + \langle a_3 , a_4 \rangle = \langle a_1 , a_2
\rangle + \langle a_3 , a_4 \rangle$. One has $0 \in \{ a_1 , a_2 \}$
or $0 \in \{ a_3 , a_4 \}$. In the first case $\langle a_1 , a_2
\rangle = 0$ and $\langle a_3 , a_4 \rangle = 1$ since they are
distinct $\ne 0$. Similarly if $0 \in \langle a_3 , a_4 \rangle$ we
get $\langle a_1 , a_2 \rangle + \langle a_3 , a_4 \rangle = 1$ in all
cases. We can thus replace (\ref{eq22}) by
\begin{equation}
\label{eq23}
- \, U^3 \, \sigma (a_1) \, \sigma (a_2) \, \sigma (a_3) \, \sigma (c)
\, \sigma (a_4) \, U \, \sigma (c) \, .
\end{equation}
Permuting $c$ with $a_4$ gives a $(-1)^{\langle c , a_4 \rangle}$. We
have,
\begin{equation}
\label{eq24}
  \sigma (a_1)
\, \sigma (a_2) \, \sigma (a_3) \, \sigma (a_4) = (-1)^s \, \sigma_1
\, \sigma_2 \, \sigma_3 \, = \, i \, (-1)^s.
\end{equation}
where $(-1)^s$ is
  the signature of the permutation of $\{1,2,3\}$ given by the {\em
non zero} $a_j$'s.
The coefficient of $U^3 \, \sigma (c) \, U \, \sigma (c)^{-1}$ is thus,
\begin{equation}
\label{eq25}
- 4^{-4}  \Sigma \,  i \, \varepsilon (a) \, (-1)^s \, (-1)^{\langle 
c , a_4 \rangle} \, .
\end{equation}
Taking $c = (1,0)$ we find 16 $-$ signs and 8 $+$ signs so that we get
the term,
\begin{equation}
\label{eq26}
4^{-4} (- (-16+8) i) \, U^3 \, \sigma_2 \, U \, \sigma_2 =
\frac{i}{32} \ U^3 \, \sigma_2 \, U \, \sigma_2 \, .
\end{equation}
Taking $c = (0,1)$ we also find 16 $-$ signs and 8 $+$ signs which
gives
\begin{equation}
\label{eq27}
\frac{i}{32} \ U^3 \, \sigma_1 \, U \, \sigma_1 \, .
\end{equation}
Thus if we let $\mu = [x_0 , x_1 , x_2 , x_3]$ and compute $\mu \mu^*$
we get terms of the form,
\begin{equation}
\label{eq28}
\frac{-1}{(32)^2} \ U^3 \, \sigma_1 \, U \, \sigma_1 \, \sigma_2 \,
U^{-1} \, \sigma_2 \, U^{-3}
\end{equation}
which cannot be simplified and do not appear in the product $\mu^*
\mu$ where we always have negative powers for the first $U$'s on the
left followed by positive powers.

\noindent Thus we conclude that in the universal algebra $\cal A$
one has
\begin{equation}
\label{eq29}
[\mu , \mu^*] \ne 0 \, .
\end{equation}

   \end{document}